\documentclass[reqno, 12pt]{report}
\usepackage{amssymb, latexsym, amsmath, amscd, amsmath, amsthm, thesis, epsf}
\newtheorem{theorem}{Theorem}[section]
\newtheorem{corollary}[theorem]{Corollary}
\newtheorem{lemma}[theorem]{Lemma}
\newtheorem{proposition}[theorem]{Proposition}

\theoremstyle{definition}
\newtheorem{definition}[theorem]{Definition}
\newtheorem{example}[theorem]{Example}

\theoremstyle{remark}
\newtheorem*{remark}{Remark}
\newtheorem*{note}{Note}

\everymath{\displaystyle}

\newcommand{\Syl}{\mbox{\rm{Syl}}}
\newcommand{\Hol}{\mbox{\rm{Hol}}}
\newcommand{\scriptHol}{\mbox{\scriptsize\rm{Hol}}}
\newcommand{\Aut}{\mbox{\rm{Aut}}}
\newcommand{\scriptAut}{\mbox{\scriptsize\rm{Aut}}}
\newcommand{\Hom}{\mbox{\rm{Hom}}}
\newcommand{\Z}{{\mathbb Z}}
\newcommand{\F}{{\mathbb F}}
\newcommand{\N}{{\mathbb N}}

\newcommand{\Ker}{\mbox{\rm{Ker}}}
\newcommand{\moplus}{\mathop{\oplus}}
\newcommand{\free}{\mathop{\ast}}
\newcommand{\msplit}{\buildrel \curvearrowleft \over \rightarrow}
\newcommand{\C}{{\mathfrak C}}
\newcommand{\M}{{\mathfrak M}}
\newcommand{\Freemodule}{{\mathfrak F}}
\newcommand{\Ob}{\mbox{\rm{Ob}}}
\newcommand{\Mor}{\mbox{\rm{Mor}}} 

\begin{document}  

\thispagestyle{empty}

\vfill

\centerline{\large On the Holomorph of a Discrete Group}

\vfill

\centerline{\large by}

\vfill

\centerline{\large Maria S. Voloshina}

\vfill

\vfill 

\centerline{\large Submitted in Partial Fulfillment}

\vfill

\centerline{\large of the}

\vfill

\centerline{\large Requirements for the Degree}

\vfill

\centerline{\large Doctor of Philosophy}

\vfill

\vfill

\centerline{\large Supervised by}

\vfill

\centerline{\large Professor Frederick R. Cohen}

\vfill

\centerline{\large Department of Mathematics}
\centerline{\large The College}
\centerline{\large Arts \& Sciences}

\vfill

\centerline{\large University of Rochester}
\centerline{\large Rochester, New York}

\vfill

\centerline{\large 2003} 

\newpage

\pagestyle{kisa}
\pagenumbering{roman}
\setcounter{page}{2}

\centerline{\large Curriculum Vitae}

\vspace{1cm}

The author was born in Moscow, Russia on 08 September, 1975. She attended Moscow State University from 1992 to 1998,
and graduated with Diploma from the Faculty of Mechanics and Mathematics in 1998. She came to the University of
Rochester in the Fall of 1998 and began graduate studies in Mathematics. She received the Master of Arts degree from
the University of Rochester in 2000. She received Dean's Teaching Fellowship in the Spring and Fall of 2001, Spring of
2002, and Spring of 2003. She pursued her research in the area of cohomology of groups under the direction of
Professor Frederick R. Cohen, and was supported by an NSF grant in the Fall of 2002.

\newpage

\centerline{\large Acknowledgements}

\vspace{1cm}

I would like to express my deep gratitude to my advisor Frederick Cohen for introducing me to the
subject, suggesting this topic for my dissertation, his enormous help, many suggestions on improvement of
this text, his incredible patience, and encouragement. His aid can not be overestimated.

I am very thankful to Jonathan Pakianathan for pointing out the congruence subgroups, helping me to go
through his thesis and other papers of which we make a heavy use in sections \ref{congruence_subgroups}
and \ref{minh_simonds}, and many useful conversations.

I would also like to thank Frederick Cohen, Douglas Ravenel, Saul Lubkin, and Naomi Jochnowitz for their excellent 
courses in topology and algebra.

Finally, I thank my friends Zokhrab Moustafaev, Raluca Felea, Jonathan Pakianathan, Michael Knapp, Iryna Labachova, 
Inga Johnson, and Sung Eun Kim for bringing fun and diversity to my days at the University of Rochester.

\newpage

\centerline{\large Abstract}

\vspace{1cm}

The holomorph of a discrete group $G$ is the universal semi-direct product of $G$. In chapter 1 we describe why it is an 
interesting object and state main results.
 
In chapter \ref{general-chapter} we recall the classical definition of the holomorph as well as this universal property, and 
give some group theoretic properties and examples of holomorphs. In particular, we give a necessary and sufficient condition 
for the existence of a map of split extensions for holomorphs of two groups. 

In chapter \ref{cyclic-chapter} we construct a resolution for  Hol$({\mathbb Z}_{p^r})$ for every prime $p$, where
$\Z_m$ denotes a cyclic group of order $m$, and use it to compute the integer homology and mod $p$ cohomology ring of 
Hol$({\mathbb Z}_{p^r})$. 

In chapter \ref{matrices-chapter} we study the holomorph of the direct sum of several copies of $Z_{p^r}$. 
We identify this holomorph as a nice subgroup of $GL(n+1, Z_{p^r})$, thus its cohomology informs on the 
cohomology of the general linear group which has been of interest in the subject. 
Using a map of split extensions from chapter \ref{general-chapter}, we show that the Lyndon-Hochschild-Serre spectral sequence 
for $\displaystyle H^*\left(\Hol\left(\moplus_n{\mathbb Z}_{p^r}\right); \F_p \right)$ does not collapse at the $E_2$ stage 
for $p^r\ge 8$. Also, we compute mod $p$ cohomology and the first Bockstein homomorphisms of the congruence subgroups given by
$\mbox{Ker} \left(\displaystyle \Hol\left(\moplus_n{\mathbb Z}_{p^r}\right) \rightarrow 
\Hol\left(\moplus_n{\mathbb Z}_{p}\right) \right) . $ 

In chapter \ref{wreath-chapter} we recall wreath products and permutative categories, and their connections
with holomorphs. We show that the families $ \{ \Aut(G^n), n\ge 0 \}$ and $ \{ \Hol(G^n), n\ge 0 \}$ form  permutative 
categories.

In chapter \ref{injectives-chapter} we give a short proof of the well-known fact due to S. Eilenberg and J. C. Moore that 
the only injective object in the category of groups is the trivial group.

\thispagestyle{kisa}

\tableofcontents

\thispagestyle{kisa}

\thispagestyle{kisa}

\listoffigures

\newpage

\pagenumbering{arabic}
\setcounter{page}{1}

\chapter{Introduction}

\thispagestyle{kisa}

\section{What is the holomorph of a group and why to study it}

The holomorph of a group $K$, denoted by $\Hol(K)$ and defined below, is a semi-direct product of 
$K$ and its automorphism group $\Aut(K)$. It is the universal semi-direct product of $K$ in the sense described below. 
A semi-direct product of groups $K$ and $H$, denoted by $K \ltimes H$, is a variation of a product: 
elements of $K \ltimes H$ are pairs $(k,h)$ where $k$ is in $K$ and $h$ is in $H$, but the multiplication is twisted, i.e. 
$(k,h)(l,m) \not= (kl,hm)$. 

\begin{figure}
\centerline{\epsfbox{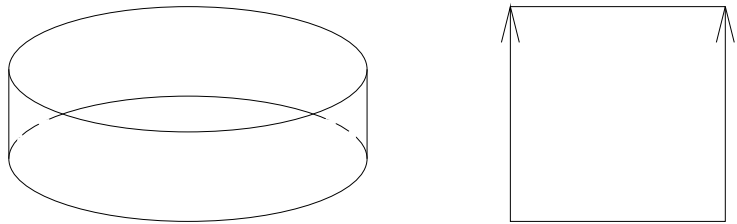}}
\caption{Untwisted band}
\label{untwisted-band}
\end{figure}

\begin{figure}
\centerline{\epsfbox{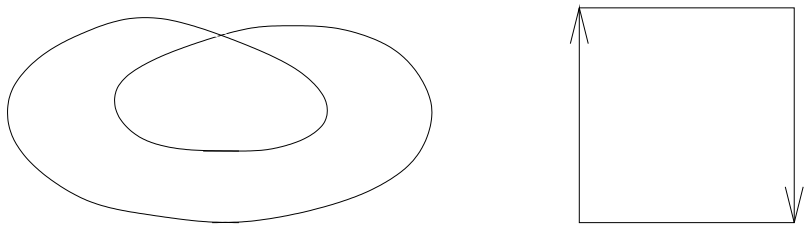}}
\caption{M\"obius band}
\label{mobius-band}
\end{figure}

The difference between direct products  and semi-direct products of groups is like the difference between Cartesian 
products and 
``twisted'' products of spaces. For example, an untwisted band shown in figure \ref{untwisted-band} is the Cartesian product of a circle 
and an interval, while a M\"obius band shown in figure \ref{mobius-band} is a ``twisted'' product of a circle and an interval. Both can 
be obtained from a square by gluing one pair of opposite edges, but for the M\"obius band the orientation of one side is reversed. These 
two spaces are not homeomorphic. However, they are homotopy equivalent (being both homotopy equivalent to a circle).

\begin{figure}
\centerline{\epsfbox{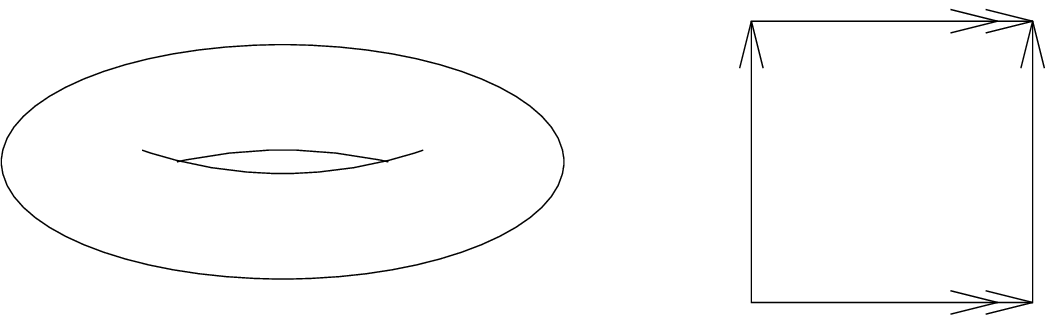}}
\caption{Torus}
\label{torus}
\end{figure}

\begin{figure}
\centerline{\epsfbox{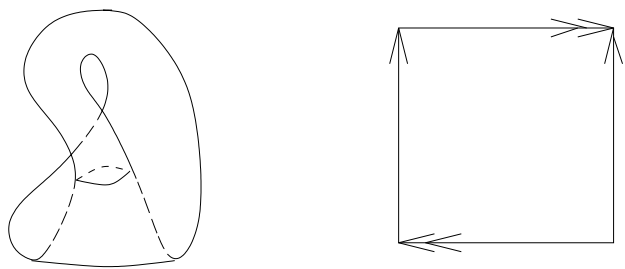}}
\caption{Klein bottle}
\label{klein}
\end{figure}

A more interesting example is a torus (surface of a donut), the Cartesian product of two circles, and a Klein bottle, 
a ``twisted'' product of two circles. A torus can be obtained from a square by gluing the opposite edges as shown in 
figure \ref{torus}. A Klein bottle is also obtained from a square by gluing the opposite edges, but now one pair 
should be glued in a twisted way: see figure \ref{klein}.

A torus and a Klein bottle are not homotopy equivalent. This can be shown for example by computing their fundamental groups and 
showing that they are not isomorphic. The fundamental group of a torus is
$$<a,b \ | \ aba^{-1}b^{-1}=1> \ \cong \ <a,b \ | \ ab=ba> \ \cong \ \Z \times \Z,$$
the direct sum of two copies of $\Z$. The fundamental group of a Klein bottle is
$$<a,b \ | \ abab^{-1}=1> \ \cong \ <a,b \ | \ ab=ba^{-1}> \ \cong \ \Z \ltimes \Z,$$
a semi-direct product of two copies of $\Z$. These two groups
are not isomorphic e.g. because the first one is commutative while the second is not.

Semi-direct products occur frequently in algebra and topology, and have been studied a lot. Thus the holomorph of a group, 
being the 
 universal semi-direct product in the sense described below, is an interesting object to study. It is constructed as a 
subgroup 
of the ``symmetric group on the set of elements of the group''. More precisely, let $S(K)$ be the group of set automorphisms of 
the set of elements of $K$ with multiplication given by composition. If $K$ is finite, then $S(K)$ is the symmetric group on 
$|K|$ letters. Let 
$$\begin{CD}
K @>i_K>> S(K)
\end{CD}$$
\noindent be the Cayley map. Also, consider the action of $\Aut(K)$ on the set of elements of $K$ given by 
$f:a\mapsto f(a)$ for $f\in \Aut(K)$ and $a\in K$. This gives an inclusion
$$\begin{CD}
\Aut(K) @>{i_{\scriptsize \Aut(K)}}>> S(K).
\end{CD}$$

\noindent {\bf Definition \ref{hol-definition}.} The holomorph of a group $K$, denoted by $\Hol(K)$, is the subgroup of $S(K)$ generated 
by the images of $i_K$ and $i_{\scriptsize \Aut(K)}$. 

Below we list some facts about $\Hol(K)$. 

\noindent {\bf Facts.} 

\bigskip

\noindent {\bf 1.} \ There is a split short exact sequence of groups
$$1 \rightarrow K \rightarrow \mbox{Hol}(K) \msplit \mbox{Aut}(K) \rightarrow 1.$$

\bigskip  

\noindent {\bf 2.} \ The group $\Hol(K)$ is the universal semi-direct product of $K$, namely  
\begin{enumerate}
\item[\bf (a)] For any group homomorphism $\phi: H \rightarrow \Aut(K)$ there is a map of split extensions
$$\begin{CD}
1 @>>> K    @>i_1>> \xi_{\phi, j_2}  @>\buildrel s_1 \over \curvearrowleft>j_1> H          @>>> 
1 \\
@.     @|           @VVfV                                                       @VV{\phi}V      @. \\
1 @>>> K    @>i_2>> \Hol(K)          @>\buildrel s_2 \over \curvearrowleft>j_2> \Aut(K)    @>>> 1
\end{CD}$$
where $\xi_{\phi, j_2}$ is the pullback. 
\item[\bf (b)] For any split extension
$1 \rightarrow K \rightarrow G \buildrel \curvearrowleft \over \rightarrow H \rightarrow 1$
there is a unique homomorphism
$H \rightarrow \Aut(K)$ such that the pullback in the diagram above is isomorphic to $G.$
\item[\bf (c)] That is, there is a one-to-one correspondence between the set of split extensions of $K$ and group 
homomorphisms with target $\Aut(K)$. 
\end{enumerate}

\noindent {\bf 3.} \ Let $P$ be a permutation group on $q$ letters (i.e. a subgroup of $S_q$). The wreath product $P \wr 
G$ is isomorphic to  the pullback $\xi$ in the diagram
$$\begin{CD}
\xi  @>>>      P        \\
@VVV           @VVV     \\
\Hol(G^q) @>>> \Aut(G^q) 
\end{CD}$$
where $P \rightarrow \Aut(G^q)$ is the composite $P \hookrightarrow S_q \hookrightarrow \Aut(G^q).$
  
\noindent {\bf 4.} \ The group $\Hol\left( \moplus_n \Z_2 \right)$ is isomorphic to a maximal subgroup of 
$GL(n+1, \Z_2)$. 

\bigskip

\noindent {\bf 5.} \ Let $B_G$ denote a classifying space for a group $G$, and let $E_G$ denote a principal $G$-bundle. 
The action of $\Aut(G)$ 
on $G$ induces an action of $\Aut(G)$ on $B_G$, and we have the associated bundle 
$$ B_G \rightarrow B_G \times_{\scriptAut(G)} E_{\scriptAut(G)} \rightarrow B_{\scriptAut(G)},$$
where $B_G \times_{\scriptAut(G)} E_{\scriptAut(G)}$ is a classifying space for $\Hol(G)$. Thus if $G$ is discrete, then 
$B_G \times_{\scriptAut(G)} E_{\scriptAut(G)}$ is an Eilenberg - Mac Lane space $K(\Hol(G),1)$.

Let us mention a related work by M. Furusawa, M. Tezuka, and N. Yagita. Since $S^1\times S^1$ is a $B_{\Z\oplus \Z}$, 
taking $B$ of $\Z\oplus \Z$ twice gives the bundle 
$$ B_{S^1\times S^1} \rightarrow  B_{S^1\times S^1} \times_{GL(2,\Z)} E_{GL(2,\Z)} \rightarrow B_{GL(2,\Z)}.$$
Let $\Gamma$ be a subgroup of $SL(2,\Z)\hookrightarrow GL(2,\Z)$, and take the pullback 
$$\begin{CD}
\xi                                               @>>>  \Gamma \\
@VVV                                                    @VVV   \\
B_{S^1\times S^1} \times_{GL(2,\Z)} E_{GL(2,\Z)}  @>>>  B_{GL(2,\Z)}
\end{CD}$$
The cohomology of such pullbacks were studied in \cite{FTY}. 

\bigskip

The purpose of this work is to study some general properties of holomorphs, and compute the cohomology of holomorphs for 
some natural and interesting examples of groups. 
For instance, if $\displaystyle K=\moplus_n R$ where $R$ is any commutative ring with a unit, e.g. $\Z$ or $\Z /m \Z$, then 
$\Aut(K)\cong GL(n,R)$, the general linear group of $n$ by $n$ matrices with entries in $R$ (the general linear group is the group of 
invertible matrices). Thus $\Hol\left(\moplus_n R\right)$ is the semi-direct product of $\moplus_n R$ and $GL(n,R)$ via 
the natural 
action of $GL(n,R)$ on $\moplus_n R$. The general linear group is of great interest to topologists. We show that $GL(n,R)$ and 
$\Hol\left(\moplus_n R\right)$ inform on each other.  

Homology and cohomology of groups are an interplay between algebraic topology and homological algebra. They can be defined purely 
algebraically or topologically (see \cite{AM}, \cite{B}, \cite{M}). Low dimension homology and cohomology groups have simple 
algebraic meanings, e.g. the first homology, $H_1(G)$, is the abelianization of the group $G$; if $K$ is abelian, then 
$H^2(H;K)$ classifies central extensions $1\rightarrow K\rightarrow G\rightarrow H\rightarrow 1$; and $H^3(H;K)$ gives the 
obstruction for existence of such an extension. Higher dimension homology and cohomology groups have further interpretations in terms of  
extensions.

Main results of this work are listed in the next section. 

\section{Main results}

In chapter \ref{general-chapter} we state some group theoretic properties of the holomorph. In particular, in section  
\ref{hol-map-section}, we study homomorphisms between holomorphs of two groups.
  
\noindent {\bf Definition \ref{compatible-definition}.} Let $G$ and $H$ be groups, and let $F: G \rightarrow H$ and $F': \Aut (G) 
\rightarrow
\Aut (H)$ be group homomorphisms. Then $F$ and $F'$ are called compatible if for every $x \in G$ and every
$g \in \Aut (G)$
$$F'(g)(F(x)) = F(g(x)),$$
i.e. the diagram
$$\begin{CD}
G     @>g>>     G     \\
@VVFV           @VVFV \\
H     @>F'(g)>> H
\end{CD}$$
is commutative for every $g \in \Aut (G)$.

\noindent {\bf Theorem \ref{compatibility}.} {\it 
The map of sets $F'': \Hol (G) \rightarrow \Hol (H)$ defined by\\
$ \phantom{move} \hspace{4cm} F''((f,x))=(F'(f),F(x)) \hfill (*)$ \\
is a group homomorphism if and only if $F$ and $F'$ are compatible. In this case we have a map of split extensions \\
$ \phantom{move} \hspace{2cm} 
\begin{CD}
1  @>>> G     @>i_1>> \Hol(G) @>{\buildrel s_1 \over \curvearrowleft}>j_1> \Aut(G) @>>> 1  \\
@.      @VVFV         @VVF''V                                              @VVF'V       @. \\
1  @>>> H     @>i_2>> \Hol(H) @>{\buildrel s_2 \over \curvearrowleft}>j_2> \Aut(H) @>>> 1
\end{CD} \hfill (**) 
$ \\
Conversely, for any map of split extensions ($**$), $F$ and $F'$ must be compatible and $F''$ is defined by 
($*$). } 

For example, if $\{ \Gamma^i(G) \}$ is the descending central series for a group $G$, then as shown in example 
\ref{characteristic-example-Gamma}, the natural homomorphisms  
$G \rightarrow G/\Gamma^i(G) $ and $\Aut(G) \rightarrow \Aut(G/\Gamma^i(G)) $ are compatible, and hence give a map of 
split extensions 
$\Hol(G) \rightarrow \Hol(G/\Gamma^i(G)) $. Some other natural examples of compatible homomorphisms are given in 
\ref{hol-map-section}. 

In section \ref{modules-section} we study modules over the group ring of a split extension. In particular, we show that 
any free $\Z[K]$-module can be turned into a $\Z[\Hol(K)]$-module using the natural actions of $K$ and $\Aut(K)$ on 
$\Z[K]$. This fact generalizes to the case of an arbitrary split extension. 

\noindent {\bf Theorem \ref{hol-module-theorem}.} {\it 
Let $\M_G$, $\M_K$, and $\Freemodule_K$ denote the categories of $\Z[G]$-modules, $\Z[K]$-modules,
and free $\Z[K]$-modules respectively. We have the  inclusion $\Freemodule_K \buildrel i \over \hookrightarrow \M_K$ and
forgetful functor $\M_G \buildrel f \over \rightarrow \M_K$. Let $\xi_{i,f}$ be the pullback, i.e. the category of
$\Z[G]$-modules that are free over $\Z[K]$. Then the top map in the diagram
$$\begin{CD}
\xi_{i,f} @>>> \Freemodule_K \\
@VVV          @VViV          \\
\M_G     @>f>> \M_K
\end{CD}$$
is a split epimorphism.   
}

In chapter \ref{cyclic-chapter} we compute the integer homology and mod $p$ cohomology ring of Hol$({\mathbb Z}_{p^r})$ 
for any prime number $p$, whose description in terms of generators and relations is given in section 
\ref{description-holZ2r} for $p=2$ and in section \ref{description-HolZpr} for $p$ odd. In case $p$ is odd 
$\Hol(\Z_{p^r})$ is metacyclic, and mod $p$ cohomology 
rings of metacyclic groups were computed by Huebschmann in \cite{H}. However, the author has learned about the work of 
Huebschmann only after this work was done. Also, we will need to know the behavior of the LHS spectral sequence in chapter 
\ref{matrices-chapter}, so we will state some intermediate results and show the spectral sequence, but omit some details of 
computations. In addition, we compute the Bockstein homomorphisms.

We will use the method described in \cite{W} to construct an explicit resolution of Hol$(K).$ Namely, in
\cite{W} Wall showed how, given a short exact sequence  
$$1\rightarrow K \rightarrow G \rightarrow H \rightarrow 1$$
where $K$ is a normal subgroup of $G$, and given resolutions for $K$ and $H$, to construct a free
resolution for $G$. He applied his method to find an explicit resolution for a split extension of
cyclic by cyclic groups, and used this to compute its integer homology.

If $r\ge 3,$ $H=\Aut(\Z_{2^r})$ is the direct sum of two cyclic groups, so we generalize Wall's resolution to cover this case.
The explicit resolution is given in section \ref{generalization-section}. 

The cohomology rings  $\Hol(\Z_{p^r})$ are as follows. 

\noindent {\bf Theorem \ref{coh-ring-p2}.} {\it 
\begin{enumerate}
\item If $r>3$, then
$H^*(\Hol(\Z_{2^r};\F_2)) \cong $\\
$\cong \Lambda(x) \otimes \F_2[a,b,c,y,z] / a^2=ax+ay, az=0, ab=by, b^2=acx+bxz+acy,$ \\
where $|x|=|y|=|a|=1,$ $|z|=2$, $|b|=3$, and $|c|=4$.
\item  $H^*(\Hol(\Z_8;\F_2)) \cong $\\
$\cong \F_2[a,b,c,x, y] / a^2=ax+ay, ax^2=0, ab=by, b^2=acx+bx^3+acy+cx^2,$ \\
where $|x|=|y|=|a|=1,$ $|b|=3$, and $|c|=4$.
\end{enumerate}
}

\noindent {\bf Theorem \ref{coh-ring-podd}.} {\it 
$$H^*(\mbox{\Hol}({\mathbb Z}_{p^r}); {\mathbb F}_p) \cong \Lambda (d_0, d_1, \ldots, d_{p-2}, e, x) \otimes {\mathbb
F}_p[f,z]/
d_i d_j = d_i e = d_i z = 0$$
where $|d_i|=2(p-1)(i+1)-1,$ $|e|=2p(p-1)-1,$ $|f|=2p(p-1),$ $|x|=1,$ $|z|=2.$
}

The last theorem also follows from theorem A, part (2), of \cite{H}.

We also compute the continuous cohomology of $\lim_{\leftarrow} ( \Hol(\Z_{p^n}) )$. 

\noindent {\bf Theorem \ref{my-continuous-theorem}.} {\it 
$H_{cont}^* (\lim_{\leftarrow} ( \Hol(\Z_{p^n}) ); \F_p) = $

$ \hfill = \left\{ \begin{array}{@{}ll@{\,}}
\Lambda(x)\otimes \F_2[u,y]/u^2=ux+uy, & |x|=|u|=|y|=1, \ \hfill p=2, \\
\Lambda(x),                  & |x|=1,         \hfill p \mbox{ is odd.}
\end{array}\right.$ \\
The Poincar\'e series of $H_{cont}^* (\lim_{\leftarrow} ( \Hol(\Z_{2^n}) ); \F_2)$ is $\frac{(1+t)^2}{1-t}$.
}

In chapter \ref{matrices-chapter} we study the holomorph of the direct sum of several copies of $\Z_{p^r}$.
This holomorph is identified as a subgroup of $GL(n+1, \Z_{p^r})$ with certain properties listed in section 
\ref{matrices-section}. For example, one property is that if $p^r=2$, then $\Hol\left( \moplus_n \Z_2 \right) $ is a 
maximal subgroup of $GL(n+1, \Z_2)$. 
If $r=1$, $\Syl_p\left(\Hol\left(\moplus_n \Z_{p}\right)\right) \cong \Syl_p\left(GL(n+1, \Z_{p})\right)$. Thus the 
cohomology of such a holomorph informs on the cohomology of the general linear group which has been of interest in the subject.
Using a map of split extensions from chapter \ref{general-chapter} and an analogue of the Dickson algebra, we show that the 
Lyndon-Hochschild-Serre (LHS) spectral sequence for the cohomology of $\displaystyle \Hol\left(\moplus_n{\mathbb 
Z}_{p^r}\right)$ does not collapse at the $E_2$ stage: 

\noindent {\bf Theorem \ref{not-collapse-theorem}.} {\it 
The LHS spectral sequence for the extension
$$1 \rightarrow \moplus_n \Z_{p^r} \rightarrow \Hol\left( \moplus_n \Z_{p^r} \right) \rightarrow \Aut\left( \moplus_n \Z_{p^r} 
\right)
\rightarrow 1$$
does not collapse at the $E_2$ stage if $p=2$ and $r\ge 3$, or $p$ is odd and $r\ge 2$. 
} 

\noindent {\bf Corollary \ref{nonzero-cohomology-corollary}.} {\it 
Let $n\ge 2$. Then
$$\displaystyle H^2\left( GL(n,\Z_{p^r}); H^{2p^n-3}\left( \moplus_n \Z_{p^r} ; \F_p \right)\right)\ne 0$$
if $p=2$ and $r\ge 3$, or if $p$ is odd and $r\ge 2$, where the local coefficient system is that in 
\ref{not-collapse-theorem}. 
}

In section \ref{congruence_subgroups} we show that for $p$ odd, the kernels of reduction maps 
$$\displaystyle \Hol\left(\moplus_n{\mathbb Z}_{p^r}\right) \rightarrow \Hol\left(\moplus_n{\mathbb Z}_{p}\right)$$
form an infinite uniform tower (the definition of a uniform tower is recalled in \ref{congruence_subgroups}), and compute 
their mod $p$ cohomology, as well as the first Bockstein homomorphisms in case $p\ge 5$: 

\noindent {\bf Theorem \ref{congruence-ring-Bock}.} {\it For any prime number $p$, 
$$H^* \left( \Gamma_{n,k}^{\scriptHol}; \F_p \right) \cong \Lambda(x_{ij}, x_{l}) \otimes \F_p[s_{ij}, s_{l}] 
\hspace{1cm} 1 \le i,j,k, \le n$$
where $|x_{ij}|=|x_l|=1$ and $|s_{ij}|=|s_l|=2$. Also, the first Bocksteins are given by
$$\begin{array}{rll}
\beta(x_{tu}) & = & \left\{
        \begin{array}{lll}
        \displaystyle - \sum_{i=t}^n x_{ti}x_{iu} + \sum_{i=1}^{t-1} x_{iu}x_{ti} & \mbox{if} & t<u \\
        \displaystyle - \sum_{i=t+1}^n x_{ti}x_{iu} + \sum_{i=1}^{t-1} x_{iu}x_{ti} & \mbox{if} & t=u \\
        \displaystyle - \sum_{i=t+1}^n x_{ti}x_{iu} + \sum_{i=1}^{t} x_{iu}x_{ti} & \mbox{if} & t>u
        \end{array} \right. \\
\beta(x_t)    & = & \displaystyle - \sum_{i=1}^n x_{ti} x_i \\
\beta(s_{tu}) & = & \displaystyle \sum_{i=1}^n (s_{ti}x_{iu}-s_{iu}x_{ti}) \\
\beta(s_t)    & = & \displaystyle \sum_{i=1}^n (s_{ti} x_i - s_i x_{ti})
\end{array}$$
}

In chapter \ref{wreath-chapter} we recall wreath products and their connections with holomorphs. In particular, we 
show that for any group $G$, $S_q \wr \Aut(G)$ and $S_q \wr \Hol(G)$ are natural subgroups of $\Aut(G^q)$ and $\Hol(G^q)$ respectively. 

Also, we show how the families of groups  
$$ \left\{ \Aut(G^n), n\ge 0 \right\},$$ 
$$ \left\{ \Hol(G^n), n\ge 0 \right\},$$
and 
$$ \left\{ \Aut\left(\free_n G\right), n\ge 0 \right\}$$
admit additional structure. 

\noindent {\bf Theorem \ref{permutative-theorem}.} {\it 
The above families form  permutative categories.
}

Chapter \ref{injectives-chapter} is independent from the preceding ones. There we give a short proof of the 
following well-known fact proved in \cite{EM}.

\noindent {\bf Theorem \ref{injectives-theorem}} {\it (S. Eilenberg and J. C. Moore) The only injective object in the category 
of groups is the trivial group. }

\chapter{Holomorph of a group} 

\label{general-chapter}

\thispagestyle{kisa}

\section{Classical definition of $\Hol(K)$}

\label{classical-def-section}

Classically (see e.g. \cite[p.90]{K}), $\Hol(K)$ arises as follows. If $K$ is a normal subgroup of a group $G$, then every inner 
automorphism of $G$ induces an
automorphism of $K$. Since inner automorphisms are easier to study than automorphisms in general, a natural question is, given a group 
$K$, does there exist a group $G$ containing $K$ as a normal subgroup, and such that every automorphism of $K$ is induced by some inner 
automorphism of $G$? The answer is ``yes'', and an example of such a group $G$ is $\Hol(K)$. 

Let $S(K)$ be the group of set automorphisms of the set of elements of $K$ with multiplication
given by composition. If $K$ is finite, then $S(K)$ is the symmetric group on $|K|$ letters.
The group $K$ acts on the set of its elements by left translation: 

\centerline{if $x,a\in K$, then $x:a\mapsto xa$.} This action gives the inclusion 

\begin{equation}\label{Cayley-map} 
\begin{CD}
K @>i_K>> S(K)
\end{CD}
\end{equation}
known as Cayley map. 
Also, the group of automorphisms of $K$, $\Aut(K)$, acts on the set of elements of $K$ by applying the automorphism: 

\centerline{if $f\in \Aut(K)$ and $a\in K$, then $f:a\mapsto f(a)$.} 
\noindent This gives an inclusion 

\begin{equation}\label{Aut-inclusion} 
\begin{CD}
\Aut(K) @>{i_{\scriptsize \Aut(K)}}>> S(K). 
\end{CD}\end{equation}

\begin{definition} \label{hol-definition}
The holomorph of a group $K$, denoted by $\Hol(K)$, is the subgroup of $S(K)$ generated by the images of 
\eqref{Cayley-map} and \eqref{Aut-inclusion}. 
\end{definition}

Equivalently, $\Hol(K)$ is the normalizer of the $K$ in $S(K)$ (see \cite{K}).

\section{$\Hol(K)$ as a split extension}

\label{hol-as-split-section}

For $f\in\Aut(K)$ and $x\in K$, 
$$(i_K (x)\circ i_{\scriptsize \Aut(K)}(f)) (a) = xf(a)=f(f^{-1}(x)a) = 
(i_{\scriptsize \Aut(K)}(f) \circ i_K (f^{-1}(x)))(a),$$
therefore every element of $\Hol(K)$ can be written as a composition $i_{\scriptsize \Aut(K)}(f)\circ i_K(x)$ for some 
$f\in\Aut(K)$ and $x\in K$. 
Now let us denote this composition by $(f,x)$. Then $i_K(x)=(1,x)$ and $i_{\scriptsize \Aut(K)}(f)=(f,1)$.  

Thus elements of $\Hol(K)$ can be written as pairs $(f,x)$ where $f\in\Aut(K)$ and $x\in K$, and multiplication is given by 
$$(f,x)(g,y)=(f,1)(1,x)(g,1)(1,y)=(f,1)(g,g^{-1}(x))(1,y)=(fg,g^{-1}(x)y).$$

A trivial but a useful fact is that the inverse of $(f,x)$ is $$(f,x)^{-1}=(f^{-1},f(x^{-1})).$$

From now on, we will identify $K$ with its image in $\Hol(K)$. 

\begin{lemma}
$K$ is normal in $\Hol(K)$, and every automorphism of $K$ is induced by an inner automorphism of $\Hol(K)$.
\end{lemma}

\begin{proof}
$(f,x)(1,y)(f,x)^{-1} = (f,x)(1,y)(f^{-1},f(x^{-1})) = (1, f(xyx^{-1}))$. This proves the first statement of the lemma.

In particular, for $x=1$ we have $(f,1)(1,y)(f,1)^{-1}=(1, f(y))$ which means that the automorphism $f$ of $K$ is induced by the 
conjugation by the element $(f,1)$ in $\Hol(K)$. 

\end{proof}
 
Clearly, the quotient group $\Hol(K)/K$ is isomorphic to $\Aut(K)$. Thus we say that $\Hol(K)$ is the extension of $K$ by $\Aut(K)$ via 
the natural action of $\Aut(K)$ on $K$. 

It is well-known that the short exact sequence
$$1 \rightarrow K \rightarrow \mbox{Hol}(K) \rightarrow \mbox{Aut}(K) \rightarrow 1.$$
is always split, namely there is a section $\Aut(K) \buildrel s \over \rightarrow \Hol(K)$ defined by $s(f)=(f,1),$
and that Hol$(K)$ is the universal semi-direct product of $K$ in the sense of the theorem \ref{universal_theorem} below.
For $K$ abelian, this theorem is proved in \cite{BB}. The author does not know a reference for the general 
case. 

First recall the definition of a pullback. 

\begin{definition} Let $H$, $G$, and $L$ be groups. Given homomorphisms $H\buildrel \phi \over \rightarrow L$ and 
$G \buildrel \theta \over \rightarrow L$, the pullback $\xi_{\phi, \theta}$ is the subgroup of $H\times G$ consisting of 
all pairs $(h,g)$ such that $\phi(h)=\theta(g)$. We then have the commutative diagram 
$$\begin{CD}
\xi_{\phi, \theta} @>p_1>>    H   \\
@VVp_2V                       @VV{\phi}V \\
G                  @>{\theta}>> L
\end{CD}$$
where $p_1$ and $p_2$ are projection onto the first and second coordinates respectively.
\end{definition}

\begin{theorem} \label{universal_theorem}

(1) For any group homomorphism $\phi: H \rightarrow \Aut(K)$ there is a map of split extensions

\begin{equation}\label{universal_eq}
\begin{CD}
1 @>>> K    @>i_1>> \xi_{\phi,j_2}   @>\buildrel s_1 \over \curvearrowleft>j_1> H         @>>> 1 \\
@.     @|           @VVfV                                                       @VV{\phi}V      @. \\
1 @>>> K    @>i_2>> \Hol(K)          @>\buildrel s_2 \over \curvearrowleft>j_2> \Aut(K)   @>>> 1 
\end{CD}
\end{equation}

(2) For any split extension
$1 \rightarrow K \rightarrow G \buildrel \curvearrowleft \over \rightarrow H \rightarrow 1$
there is a unique map
$H \rightarrow \Aut(K)$ such that the pullback in the diagram above is isomorphic to $G.$

\end{theorem}

\begin{proof}

(1) An element of the pullback has the form $(h,(a,x))$ where $h\in H,$ $(a,x)\in \mbox{Hol}(K)$ s.t. $\phi(h)=j_2((a,x))=a.$
Define $s_1$ and $i_1$ by $s_1(h)=(h,(\phi(h),1_K))$ and $i_1(x)=(1_H, (1_{\scriptsize \Aut(K)}, x)).$ Then
$$1 \longrightarrow K \buildrel i_1 \over \longrightarrow \xi_{\phi,j_2} \longrightarrow H \longrightarrow 1 $$
is a split extension and \eqref{universal_eq} is clearly commutative.

(2) Given
$\ 1 \longrightarrow K \buildrel i_3 \over \longrightarrow G$ 
$\renewcommand{\baselinestretch}{0.5}\normalsize 
{\begin{array}{c} {\scriptsize s_3} \\ \curvearrowleft \\ \longrightarrow \\ {\scriptsize j_3} \\  \end{array} }$ 
$ H \longrightarrow 1 $, \
define $\ \phi: H \rightarrow \mbox{Aut}(K) \ $ by
\newline $\phi(h)(x)=i_3^{-1}(s_3(h)^{-1}xs_3(h)),$ and obtain the diagram \eqref{universal_eq}.

Define $\theta: \xi_{\phi,j_2} \rightarrow G$ by $\theta((h,(a,x)))=s_3(h)i_3(x).$ It is easy to see that

$$\begin{CD}
1  @>>> K  @>i_1>> \mbox{pullback} @>{\buildrel s_1 \over \curvearrowleft}>j_1> H  @>>> 1 \\
@.      @|         @VV{\theta}V                                                 @|        \\
1  @>>> K  @>i_3>> G               @>{\buildrel s_3 \over \curvearrowleft}>j_3> H  @>>> 1 
\end{CD}$$
is commutative, therefore $\theta$ is an isomorphism by the 5-lemma.
\end{proof}

\begin{example}
Let $G=K \times H$. Then the homomorphism $H \buildrel \phi \over \longrightarrow \Aut(K)$ in \eqref{universal_eq} is
the trivial homomorphism.
\end{example}

\section{Other properties and examples of holomorphs}

\label{properties-examples-section}

\begin{example}
If $K\cong\Z_3$, $\Aut(K)\cong\Z_2$. It is easy to check that $\Hol(\Z_3)$ is not abelian and has order 6, thus 
$\Hol(\Z_3)\cong S_3$. 
\end{example}

In fact, $\Hol(K)$ is abelian only in two cases as shown in the lemma below. 

\begin{lemma}
The group $\Hol(K)$ is not abelian unless $K$ is trivial or $K \cong \Z_2$.
\end{lemma}

\begin{proof}

$(f,x)(g,y)=(g,y)(f,x)$ implies that $fg=gf$ and

\begin{equation} \label{non_abelian}
g^{-1}(x)y=f^{-1}(y)x.
\end{equation}

In particular, for $y=1$, we have $g^{-1}(x)=x$, thus $\Aut(K)$ is trivial.
Then \eqref{non_abelian} implies $xy=yx$, i.e. $K$ is abelian.
If $K$ has a cyclic subgroup of order at least 3, then $\Aut(K)$ is not trivial because the automorphism group of a
cyclic group of order $n\ge 3$ is non-trivial (it consists of multiplications by numbers between 1 and $n$, and
relatively prime to $n$).
$\Aut(\Z_2 \oplus \Z_2)$ is also non-trivial, thus $K$ is either a trivial group (in this case $\Hol(K)$
is also trivial), or isomorphic to $\Z_2$ (and so is $\Hol(K)$).
\end{proof}

Holomorphs of cyclic groups are studied in chapter \ref{cyclic-chapter}. 

\begin{example}
If $K=\Z_2 \oplus \Z_2$, obviously every permutation of the 3 nonidentity elements is an automorphism of $K$, so 
$\Aut(K)\cong S_3$, thus $\Hol(\Z_2 \oplus \Z_2)$ is an extension of $\Z_2 \oplus \Z_2$ by $S_3$. We know that $\Hol(K)$ is a 
subgroup of $S(K)$, and $|\Hol(\Z_2 \oplus \Z_2)|=|\Z_2 \oplus \Z_2|\cdot |S_3| = 4 \cdot 6 = 24 = |S_4|$. Therefore 
$\Hol(\Z_2 \oplus \Z_2) \cong S_4$. The extension 
$1\rightarrow \Z_2 \oplus \Z_2 \rightarrow S_4 \rightarrow S_3 \rightarrow 1$ is an exceptional one appearing in work of 
Dickson. This is one of only four examples when $\Hol(K) \cong S(K)$ as proved in the following lemma.
\end{example}

\begin{lemma} 
The inclusion $\Hol(K) \hookrightarrow S_{|K|}$ is an isomorphism if and only if $K$ is isomorphic to the trivial group, $\Z_2$, $\Z_3$, 
or $\Z_2 \oplus \Z_2$.
\end{lemma}

\begin{proof}
If $\Hol(K) \cong S_{|K|}$ then $|\Hol(K)| = |S_{|K|}| = |K|!$. Therefore
$\displaystyle |\Aut(K)| = \frac{|\Hol(K)|}{|K|} = (|K|-1)!$.
Every automorphism preserves the identity element, so the maximum possible order of $\Aut(K)$ is $(|K|-1)!$.
Thus in order for  $\Hol(K)$ and $S_{|K|}$ to be isomorphic, every permutation preserving $1$ has to be an automorphism.
For $|K| \ge 4$, let $x$ and $y$ be distinct nonidentity elements which are not inverses of each other. For any 
automorphism $f$ of $G$, $f(xy) = f(x)f(y)$, i.e. $f(x)$ and $f(y)$ determine $f(xy)$, therefore $|K| \le 4$.

If $K$ is isomorphic to the trivial group, $\Z_2, \ \Z_3,$ or $\Z_2 \oplus \Z_2$, it is easy to check that $|\Hol(K)| =
|K|!$, so $\Hol(K) \cong S_{|K|}$, and if $K \cong \Z_4$, $|\Hol(K)| = 8$, so $\Hol(K) \not\cong S_{|K|}$.
\end{proof}

\begin{example}
More generally, if $\displaystyle K=\moplus_n \Z_m$, then $\Aut(K)\cong GL(n,\Z_m)$. The structure and cohomology of 
$\displaystyle \Hol\left(\moplus_n \Z_m\right)$ are studied in chapter \ref{cyclic-chapter}. In particular, we will show in 
section \ref{matrices-section} that $\displaystyle \Hol\left(\moplus_n \Z_m\right)$ is isomorphic to the group of 
$GL(n+1,\Z_m)$ fixing the first basis vector.
\end{example}

Homology of holomorphs of free groups have been studied by Craig Jensen \cite{J1} and \cite{J2}.

\section{Maps of holomorphs}

\label{hol-map-section}

\begin{definition} \label{compatible-definition}
Let $G$ and $H$ be groups, and let $F: G \rightarrow H$ and $F': \Aut (G) \rightarrow \Aut (H)$ be group
homomorphisms. Then $F$ and $F'$ are called compatible if for every $x \in G$ and every $g \in \Aut (G)$
$$F'(g)(F(x)) = F(g(x)),$$
i.e. the diagram
\begin{equation}\label{compatible_diagram}\begin{CD}
G     @>g>>     G     \\
@VVFV           @VVFV \\
H     @>F'(g)>> H     
\end{CD}\end{equation}
is commutative for every $g \in \Aut (G)$.
\end{definition}

\begin{theorem} \label{compatibility}
The map of sets $F'': \Hol (G) \rightarrow \Hol (H)$ defined by
\begin{equation}\label{hol_map_def}
F''((f,x))=(F'(f),F(x))
\end{equation}
is a group homomorphism if and only if $F$ and $F'$ are compatible. In this case we have a map of split extensions
\begin{equation}\label{map_of_extensions}
\begin{CD}
1  @>>> G     @>i_1>> \Hol(G) @>{\buildrel s_1 \over \curvearrowleft}>j_1> \Aut(G) @>>> 1  \\
@.      @VVFV         @VVF''V                                              @VVF'V       @. \\
1  @>>> H     @>i_2>> \Hol(H) @>{\buildrel s_2 \over \curvearrowleft}>j_2> \Aut(H) @>>> 1   
\end{CD}
\end{equation}
Conversely, for any map of split extensions \eqref{map_of_extensions}, $F$ and $F'$ must be compatible and $F''$ is defined by
\eqref{hol_map_def}.
\end{theorem}

\begin{proof}
The map $F''$ is a group homomorphism if and only if for every $(f,x)$ and $(g,y)$ in $\Hol (G)$ we have
\begin{equation}\label{compatibility_condition} F''((f,x)) F''((g,y)) = F''((f,x)(g,y)) \end{equation}
By our definitions, \\
$F''((f,x))F''((g,y)) = (F'(f),F(x)) (F'(g),F(y)) = $ \\
$= (F'(f)F'(g), (F'(g))^{-1}(F(x))F(y)) = (F'(fg), F'(g^{-1})(F(x))F(y))$ \\
and \\
$F''((f,x)(g,y)) = F''(fg,g^{-1}(x)y) = (F'(fg),F(g^{-1}(x)y)) = $ \\
$= (F'(fg), F(g^{-1}(x))F(y)), $ \\
thus \eqref{compatibility_condition} holds if and only if $F'(g^{-1})(F(x)) = F(g^{-1}(x)).$
Clearly in this case \eqref{map_of_extensions} is commutative.

Conversely, if we have a map of split extensions \eqref{map_of_extensions}, then \\
$F''((f,1))=F''(s_1(f))=s_2(F'(f))=(F'(f),1)$ and \\
$F''((1,x))=F''(i_1(x))=i_2(F(x))=(1,F(x)).$ \\
Since $(f,x)=(f,1)(1,x)$, \eqref{hol_map_def} is true for every element of $\Hol(G)$.
As shown above, in this case $F$ and $F'$ are compatible.
\end{proof}

\begin{remark}
It is not true that every homomorphism $F'': \Hol(G) \rightarrow \Hol(H)$ is defined by \eqref{hol_map_def} for some compatible $F$ and 
$F'$ as shown in the next example.
\end{remark}

\begin{example} \label{not_map_of_split_ext}
Let $G \cong \Z_3$, $\Aut(\Z_3) \cong \Z_2$, $\Hol(\Z_3) \cong S_3$, and $H \cong \Z_2$, $\Aut(\Z_2)$ is trivial, 
$\Hol(\Z_2) \cong \Z_2$. The only possible maps \mbox{$F: \Z_3 \rightarrow \Z_2$} and $F': \Z_2 \rightarrow 1$ are trivial 
maps, however there exists a nontrivial map \mbox{$F'': S_3 \rightarrow  \Z_2$.}
\end{example}

Note that in the above example we do have a map of extensions

$$\begin{CD}
1  @>>> \Z_3 @>>> S_3  @>>> \Z_2 @>>> 1  \\
@.      @VVV      @VVV      @VVV      @. \\  
1  @>>> \Z_2 @>>> \Z_2 @>>> 1    @>>> 1
\end{CD}$$
but it is not a map of $split$ extensions. Below is an example of a map $F''$ that does not fit into any map of
extensions.

\begin{example}
Let $G \cong Z_3$ as in example \ref{not_map_of_split_ext}, let $H \cong \Z_2 \oplus \Z_2$,
$\Aut(\Z_2 \oplus \Z_2) \cong S_3$, $\Hol(\Z_2 \oplus \Z_2) \cong S_4$, and let
$F'':S_3 \rightarrow S_4$ be an inclusion. Since the only possible homomorphism 
\mbox{$F:Z_3 \rightarrow \Z_2 \oplus \Z_2$} is the trivial one, $F''$ doesn't fit into any map of extensions.
\end{example}

\begin{example}
For any $F \in \Aut(G)$ there exists a unique $F':\Aut(G) \rightarrow \Aut(G)$ compatible with $F$. Namely,
$F'$ is defined by
\begin{equation} \label{iso_compatible_def}
F'(g)(x) = FgF^{-1}(x)
\end{equation}
for any $ g\in \Aut(G)$ and $x \in G$. This follows from the fact that the diagram
$$\begin{CD}
G     @>g>>     G     \\
@VVFV           @VVFV \\
G     @>F'(g)>> G
\end{CD}$$
must commute for every $ g\in \Aut(G)$. (Note: In particular, if $\Aut(G)$ is abelian, then for every automorphism $F$ of $G$ the
compatible $F'$ is the identity homomorphism.) Thus any automorphism of $G$ induces a map $F'': \Hol(G) \rightarrow \Hol(G)$, and hence
a map of extensions.
\end{example}

\begin{example} \label{product_example}
If $H\cong G_1 \times \ldots \times G_q$, then $G_i \hookrightarrow H$ and $\Aut(G_i) \hookrightarrow \Aut(H)$ are
compatible, thus  $\Hol(G_i) \hookrightarrow \Hol(H)$, and there is a map of split extensions.
Moreover, $\displaystyle \prod_{i=1}^q \Aut(G_i) \hookrightarrow \Aut(H)$,
$\displaystyle \prod_{i=1}^q \Hol(G_i) \hookrightarrow \Hol(H)$, and

$$\begin{CD}
1 @>>>  \prod_{i=1}^q G_i @>>> \prod_{i=1}^q \Hol(G_i)             @>>> \prod_{i=1}^q \Aut(G_i)            @>>> 1 \\
@.      @|                     @VVV                                     @VVV                               @.\\
1 @>>>  \prod_{i=1}^q G_i @>>> \Hol\left(\prod_{i=1}^q G_i\right)  @>>> \Aut\left(\prod_{i=1}^q G_i\right) @>>> 1
\end{CD}$$
\end{example}

\begin{proposition} \label{characteristic_proposition}
If $K$ is a characteristic subgroup of $G$ (i.e. it is preserved by every automorphism of $G$), then for every automorphism $f$ of $G$, 
and $a$ in $K$, $f(aK)=f(a)K$, thus $f$ induces an automorphism $f'$ of the quotient group $G/K$. Thus we have a map $\Aut(G) 
\rightarrow \Aut(G/K)$. This map is obviously a group homomorphism, and it is compatible with the quotient map $G \rightarrow G/K$: 
commutativity of 
$$\begin{CD}
G         @>f>>  G         \\
@VV{\pi}V        @VV{\pi}V \\
G/K       @>f'>> G/K
\end{CD}$$
follows from the definition of $f'$ (i.e. that $f'$ is induced by $f$). 
Thus we have a map of split extensions $\Hol(G) \rightarrow \Hol(G/K)$. 
\end{proposition}

There are at least two important examples of characteristic subgroups: 
\begin{example} \label{characteristic-example-Z}
For $k>l$, $p^l \Z_{p^k} \cong \Z_{p^{k-l}}$ is a characteristic subgroup of $\Z_{p^k}$ with quotient group isomorphic to 
$\Z_{p^l}$. Therefore there is a reduction homomorphism $\Hol(\Z_{p^k})\rightarrow\Hol(\Z_{p^l})$, and a map of split extensions
$$\begin{CD}
1  @>>> \Z_{p^k} @>>> \Hol(\Z_{p^k}) @>{\curvearrowleft}>> \Aut(\Z_{p^k}) @>>> 1  \\
@.      @VVV          @VVV                                 @VVV                @. \\
1  @>>> \Z_{p^l} @>>> \Hol(\Z_{p^l}) @>{\curvearrowleft}>> \Aut(\Z_{p^l}) @>>> 1
\end{CD}$$
\end{example}
This map will be used in chapters \ref{cyclic-chapter} and \ref{matrices-chapter}. 

\begin{example} \label{characteristic-example-Gamma}
For any group $G$, groups in the descending central series of $G$, $\Gamma^i(G)$ (i.e. $\Gamma^i(G)$ is generated by commutators 
of length i; e.g. $\Gamma^2(G)=[G,G]$ is the commutator subgroup of $G$), are characteristic in $G$. Thus there are maps of split 
extensions $\Hol(G) \rightarrow \Hol(G/\Gamma^i(G))$. $i=2$ is an important special case: $\Hol(G) \rightarrow \Hol(G_{ab})$ where 
$G_{ab}=G/[G,G]$ is the abelianization of $G$. 
\end{example}

\section{Modules over the group ring of a split extension}

\label{modules-section}

Let $K$ be a discrete group.

\begin{theorem} \label{hol-module-theorem} If $M$ is a left $\Z[K]$-module and a left $\Z[\Aut(K)]$-module, then
\begin{equation}\label{action-def}
(f,x)(m)=f(x(m))
\end{equation}
for $f\in \Aut(K),$ $x\in K,$ and $m\in M$, is a well-defined action of $\Hol(K)$ on $M$ (i.e. $M$ is a left
$\Z[\Hol(K)]$-module with this action) if and only if for every $f\in \Aut(K),$ $x\in K,$ and $m\in M$,
\begin{equation} \label{action-cond}  
f(x(m))=f(x)(f(m)).
\end{equation}
\end{theorem}

\begin{proof} If \eqref{action-cond} is satisfied, then \\
$(f,x)((g,y)(m)) = (f,x)(g(y(m))) = f(xg(y)(g(m))) = f(xg(y))(fg(m))$, and \\
$((f,x)(g,y))(m) = (fg, g^{-1}(x)y)(m) = fg(g^{-1}(x)y(m)) = fg(g^{-1}(x)y)(fg(m)) =$\\
$= f(xg(y))(fg(m))$, so the action is well-defined.

Conversely, if \eqref{action-def} is a well-defined action, then \\
$f(x(m)) = (f,x)(m) = (1,f(x))(f,1)(m) = (1,f(x))(f(m)) = f(x)(f(m))$.
\end{proof}

\begin{example}
For $M=\Z[K]$, consider the actions of $K$ and $\Aut(K)$ on $M$ given by \hfill 
$ x \left( \sum_{x_i\in K} n_i x_i \right) = \sum_{x_i\in K} n_i (x x_i)$ \hfill and \hfill  
$ f \left( \sum_{x_i\in K} n_i x_i \right) = \sum_{x_i\in K} n_i f(x_i).$ \\
Obviously, these actions satisfy \eqref{action-cond}, thus $\Z[K]$ can be regarded as a $\Z[\Hol(K)]$-module. Similarly, any 
free $\Z[K]$-module can be regarded as a $\Z[\Hol(K)]$-module. However, it may not be free over $\Z[\Hol(K)]$. 
\end{example}

More generally, let $1\rightarrow K \rightarrow G \buildrel \curvearrowleft \over \rightarrow H \rightarrow 1 $ be a split extension,
and \mbox{$\phi : H\rightarrow \Aut(K)$} as in theorem \ref{universal_theorem}. Let $M$ be a left $\Z[K]$-module and a left 
$\Z[H]$-module.
Then $(h,x)(m)=h(x(m))$ for $h\in H$, $x\in K$, and $m\in M$, is a well-defined action of $G$ on $M$ (i.e. $M$ is a left
$\Z[G]$-module with this action) if and only if for every $h\in H$, $x\in K$, and $m\in M$, 
\begin{equation}\label{action-cond-general}
h(x(m))=((\phi(h))(k))(h(m)).
\end{equation} 
The proof mimics the proof of theorem \ref{hol-module-theorem}. Also, the actions \\ 
$ x \left( \sum_{x_i\in K} n_i x_i \right) = \sum_{x_i\in K} n_i (x x_i)$ \hfill and \hfill  
$ h \left( \sum_{x_i\in K} n_i x_i \right) = \sum_{x_i\in K} n_i ( (\phi(h))(x_i) )$ \\
of $K$ and $H$ on $\Z[K]$ satisfy 
\eqref{action-cond-general}, thus any free $\Z[K]$-module can be turned into a $\Z[G]$-module, although not necessarily free 
over $\Z[G]$. Thus we have the following theorem. 

\begin{theorem} Let $1\rightarrow K \rightarrow G \buildrel \curvearrowleft \over \rightarrow H \rightarrow 1 $ be a split extension.  
Let $\M_G$, $\M_K$, and $\Freemodule_K$ denote the categories of $\Z[G]$-modules, $\Z[K]$-modules, and free 
$\Z[K]$-modules respectively. We have the  inclusion $\Freemodule_K \buildrel i \over \hookrightarrow \M_K$ and forgetful 
functor $\M_G \buildrel f \over \rightarrow \M_K$. Let $\xi_{i,f}$ be the pullback, i.e. the category of $\Z[G]$-modules 
that are free over $\Z[K]$. Then the top map in the diagram 
$$\begin{CD}
\xi_{i,f} @>>> \Freemodule_K \\
@VVV          @VViV          \\
\M_G     @>f>> \M_K
\end{CD}$$
is a split epimorphism.
\end{theorem}

Recall that for any ring $R$, a left $R$-module $M$ can be turned into a right $R$-module with the right action defined by
$(m)r = r^{-1}(m)$. This action is well-defined since 
$((m)r_1)r_2 = (r_1^{-1}(m))r_2 = r_2^{-1}(r_1^{-1}(m)) = (r_1r_2)^{-1}(m) = $ \\
$= (m)(r_1r_2)$.

\begin{theorem} Let $S$ be a commutative ring with a unit. 
Let $L$, $M$, and $N$ be left $S[\Hol(K)]$-modules with $L$ a trivial $S[K]$-module. Each of these can also be turned into a right
$S[\Hol(K)]$-module as above. Then there is a natural isomorphism between
$\left( L\otimes_{S}M\right) \otimes_{S[\scriptHol(K)]} N$ and $L \otimes_{S[\scriptAut(K)]} \left( M \otimes_{S[K]} N \right)$
where the action of $\Hol(K)$ on $L\otimes_{S} M$ and the action of $\Aut(K)$ on $M\otimes_{S[K]} N$ are diagonal.
More precisely, these abelian groups are obtained from $L \otimes_{S} M \otimes_{S} N$ by the same identifications, i.e. the diagram
\begin{equation}\label{modules-commute}\begin{array}{r@{}c@{}l}
                                              & L\otimes_{S} M \otimes_{S} N &           \\
                                     \swarrow &                                & \searrow  \\
(L\otimes_{S} M) \otimes_{S[K]} N           &                                & L\otimes_{S} (M \otimes_{S[K]} N) \\
\downarrow                       \hspace{1cm} &                                & \hspace{1cm} \downarrow \\
(L\otimes_{S} M) \otimes_{S[\scriptHol(K)]} N   &        \longrightarrow         & L\otimes_{S[\scriptAut(K)]} (M \otimes_{S[K]} N)
\end{array}
\end{equation}
commutes.
\end{theorem}

\begin{proof} $(L\otimes_{S} M) \otimes_{S[\scriptHol(K)]} N \cong $\\
$\cong L\otimes_{S} M \otimes_{S} N /
(l\otimes m)x\otimes n = l\otimes m \otimes x(n), (l\otimes m)f \otimes n = l\otimes m \otimes f(n) \cong $\\
$\cong L\otimes_{S} M \otimes_{S} N /
l \otimes x^{-1}(m) \otimes n = l\otimes m \otimes x(n), f^{-1}(l)\otimes f^{-1}(m)\otimes n = l\otimes m \otimes f(n),$
and $L\otimes_{S[\scriptAut(K)]} (M \otimes_{S[K]} N) \cong $\\
$\cong L\otimes_{S} M \otimes_{S} N /
l\otimes (m)x \otimes n = l\otimes m \otimes x(n), (l)f\otimes m\otimes n = l\otimes f(m \otimes n) \cong$\\
$\cong L\otimes_{S} M \otimes_{S} N /
l\otimes x^{-1}(m) \otimes n = l\otimes m \otimes x(n), f^{-1}(l)\otimes m\otimes n = l\otimes f(m)\otimes f(n) \cong$\\
$\cong L\otimes_{S} M \otimes_{S} N /
l\otimes x^{-1}(m) \otimes n = l\otimes m \otimes x(n), f^{-1}(l)\otimes f^{-1}(m')\otimes n = l\otimes m' \otimes f(n)$.

Now defining the map
$\left( L\otimes_{S}M\right) \otimes_{S[\scriptHol(K)]} N \rightarrow  L \otimes_{S[\scriptAut(K)]} \left( M \otimes_{S[K]} N 
\right)$
by sending the class of $l \otimes m \otimes n$ in $\left( L\otimes_{S}M\right) \otimes_{S[\scriptHol(K)]} N$ to the class of
$l \otimes m \otimes n$ in
$L \otimes_{S[\scriptAut(K)]} \left( M \otimes_{S[K]} N \right)$ gives an isomorphism such that \eqref{modules-commute} commutes.
\end{proof}

Again, this can be generalized to the case of an arbitrary  split extension
$$1\rightarrow K \rightarrow G \buildrel \curvearrowleft \over \rightarrow H \rightarrow 1.$$
Namely, if $L,$ $M$, and $N$ are left $S[G]$-modules with $L$ being a trivial $S[K]$-module,
then the abelian groups $\left( L\otimes_{S} M\right) \otimes_{S[G]} N$ and \mbox{$L\otimes_{S[H]} \left(M \otimes_{S[K]} 
N\right)$},
where the action of $G$ on $L\otimes_{S} M$ and the action of $H$ on $M \otimes_{S[K]} N$ are diagonal,
are obtained from $L\otimes_{S} M \otimes_{S} N$ by the same identifications, i.e. the diagram
$$\begin{array}{r@{}c@{}l}
                                                & L\otimes_{S} M \otimes_{S} N &          \\
                                       \swarrow &                                & \searrow \\
\left(L\otimes_{S} M \right) \otimes_{S[G]} N &      \longrightarrow           & L\otimes_{S[H]} \left( M \otimes_{S[K]} N \right)
\end{array} $$
commutes.

\chapter{Holomorphs of cyclic groups}

\thispagestyle{kisa}

\label{cyclic-chapter}

\section{Introduction}

Let $m=p_1^{r_1}\cdot \ldots \cdot p_n^{r^n}$ where $p_1,$ $\ldots,$ $p_n$ are distinct prime numbers. By example 
\ref{product_example}, there is a map of split extensions 
$$\begin{CD}
1 @>>> \prod_{i=1}^n \Z_{p_i^{r_i}} @>>> \prod_{i=1}^n \Hol(\Z_{p_i^{r_i}}) @>>> \prod_{i=1}^n \Aut(\Z_{p_i^{r_i}}) @>>> 1 \\
@.     @|                                @VVV                                    @VVV                               @.\\
1 @>>> \Z_m                         @>>> \Hol(\Z_m)                         @>>> \Aut(\Z_m) @>>> 1
\end{CD}$$
Since $p_i$'s are distinct, $\Aut(\Z_m) \cong \prod_{i=1}^n \Aut(\Z_{p_i^{r_i}})$ (see \cite{AM}), therefore 
$$\Hol(\Z_m) \cong \prod_{i=1}^n \Hol(\Z_{p_i^{r_i}}).$$ 

In this chapter, we compute, for any prime number $p$, the integer homology and mod $p$ cohomology of the group 
$\Hol(\Z_{p^r})$, whose description in terms of generators and relations is given in section \ref{description-holZ2r} for 
$p=2$ and in section \ref{description-HolZpr} for $p$ odd. In case $p$ is odd $\Hol(\Z_{p^r})$ is metacyclic, and mod $p$ 
cohomology rings of metacyclic groups were computed by Huebschmann in \cite{H}. However, the author has learned about the 
work of Huebschmann only after this work was done. Also, this work provided additional information: we will need to know 
the behavior of the LHS spectral sequence in chapter \ref{matrices-chapter}, so we will state some intermediate results 
and show the spectral sequence, but omit some details of computations. In addition, we compute the Bockstein homomorphisms. 

We will use the method described in \cite{W} to construct an explicit resolution of Hol$(K).$ Namely, in
\cite{W} Wall showed how, given a short exact sequence
$$1\rightarrow K \rightarrow G \rightarrow H \rightarrow 1$$
where $K$ is a normal subgroup of $G$, and given resolutions for $K$ and $H$, to construct a free
resolution for $G$. He applied his method to find an explicit resolution for a split extension of
cyclic by cyclic groups, and used this to compute its integer homology.

If $p$ is an odd prime, Aut$({\mathbb Z}_{p^r})$ is cyclic, therefore the integer homology of $\Hol(\Z_{p^r})$ 
can be computed using Wall's application. The only difficulty here is that we don't know any formula for a generator 
of Aut$({\mathbb Z}_{p^r})$ (see \cite{AM}). But it turns out that it is not necessary to know that.

The group Aut$({\mathbb Z}_2)$ is trivial, therefore Hol$({\mathbb Z}_2)\cong {\mathbb Z}_2$ whose homology 
is well-known; Aut$({\mathbb Z}_4)\cong {\mathbb Z}_2$, therefore the homology of
Hol$({\mathbb Z}_4)\cong D_8$ (dihedral group of order 8) can be computed easily using Wall's application (and is also well-known). 
So we will only work out the answer for $r\ge 3.$ In this case $H=\Aut(\Z_{2^r})$ is the direct sum of two 
cyclic groups, but we will generalize Wall's resolution to cover this case.

As in \cite{W}, we will choose the usual resolution \cite[p.251]{CE} for each cyclic group
${\mathbb Z}_n = \{ g | g^n=1 \} :$

$$\begin{array}{lc}
R_k \mbox{ is free on one generator } r_k & (k \ge 0),\\
dr_{2k+1}=(g-1)r_{2k}                   & (k \ge 0),\\
dr_{2k}=N_g r_{2k-1}                    & (k \ge 1),
\end{array}$$

where $\displaystyle N_g=\sum_{i=0}^{n-1}{g^i}.$

Also, we will make use of the following obvious lemma.

\begin{lemma} \label{obvious} Suppose we have two chain complexes
$$ \ldots \leftarrow {\mathbb Z}_{p^r} \buildrel G \over \longleftarrow {\mathbb Z}_{p^r} \buildrel F \over \longleftarrow 
{\mathbb Z}_{p^r} \leftarrow \ldots$$  
and
$$ \ldots \leftarrow {\mathbb Z}_{p^r} \buildrel G' \over \longleftarrow {\mathbb Z}_{p^r} \buildrel F' \over \longleftarrow 
{\mathbb Z}_{p^r} \leftarrow \ldots,$$     
where $G$, $F$, $G'$, and $F'$ are multiplications by $p^a,$ $p^b,$ $p^a c,$ and $p^b d$ respectively,
where $c$ and $d$ are relatively prime to $p$. Then $\mbox{Ker\,}G=\mbox{Ker\,}G'$ and $\mbox{Im\,}F = \mbox{Im\,}F'.$

Therefore, to compute the homology of a complex
$$ \ldots \leftarrow {\mathbb Z}_{p^r} \buildrel {\times m} \over \longleftarrow {\mathbb Z}_{p^r} \buildrel {\times n} \over 
\longleftarrow \Z_{p^r} \leftarrow \ldots$$
it suffices to know the largest powers of $p$ that divide $m$ and $n$.

\end{lemma}

Let $\nu_p(n)$ denote the highest exponent of $p$ that divides $n$. 

\section{Description of $\Hol(\Z_{2^r})$ for $r\ge 3$}

\label{description-holZ2r}

\begin{lemma} \label{divisible-2}
\begin{enumerate}
\item If $n$ is odd, then $\nu_2(3^n-1)=1$.
\item If $m=2^q$ where $q \ge 1$, then $\nu_2(3^m-1) = q+2.$
\item Moreover, if $m=2^q \cdot n$ where $n$ is odd, then $\nu_2(3^m-1) = q+2.$
\end{enumerate}

\end{lemma}

\begin{proof}
\begin{enumerate}
\item $\nu_2(3^n-1) = \nu_2((3-1)(3^{n-1} + \ldots + 1)) = 1.$
\item Induction on $q$.    
$q=1$: $m=2$, $ 3^2-1=2^{1+2}.$ \\
$\nu_2(3^{2^{q+1}}-1) = \nu_2(3^{2^q \cdot 2}-1) = \nu_2((3^{2^q}-1) \cdot (3^{2^q}+1)) = (q+2)+1=q+3.$
\item $\nu_2(3^m-1) = \nu_2(3^{2^q \cdot n}-1) = \nu_2((3^{2^q}-1)(3^{2^q \cdot (n-1)} + \ldots + 1)) = 
\nu_2(3^{2^q}-1) = $ 
$ = q+2.$
\end{enumerate}
\end{proof}

\begin{corollary} \label{mod-2}
For $r\ge 4$, $3^{2^{r-3}} \equiv 1+2^{r-1} (\mbox{mod } 2^r).$
\end{corollary}

Therefore, if $r\ge 4$, then multiplication by 3 generates a cyclic subgroup of order $2^{r-2}$ in
$\mbox{Aut}({\mathbb Z}_{2^r}).$ This subgroup contains multiplication by $2^{r-1}+1,$ an element
of order 2. There are two other elements of order 2 in $\mbox{Aut}({\mathbb Z}_{2^r}),$ namely
multiplication by
$2^r-1$ and multiplication by $2^{r-1}-1.$ Since the order of $\mbox{Aut}({\mathbb Z}_{2^r})$ is
$2^{r-1},$ multiplication by 3 and multiplication by $2^r-1$ generate all of $\mbox{Aut}({\mathbb
Z}_{2^r}).$\\
If $r=3,$ $\mbox{Aut}({\mathbb Z}_8) \cong {\mathbb Z}_2 \oplus {\mathbb Z}_2,$ and we also can choose
multiplication by 3 and multiplication by $2^r-1=7$ as generators of the ${\mathbb Z}_2$'s.\\
Thus we have

\begin{corollary} \label{generated}
If $r\ge 3$, then 
$ \Aut({\mathbb Z}_{2^r}) \cong {\mathbb Z}_{2^{r-2}} \oplus {\mathbb Z}_2$, where ${\mathbb Z}_{2^{r-2}}$
is generated by multiplication by 3, and ${\mathbb Z}_2$ is generated by multiplication by $2^r-1$.

\end{corollary}

\begin{corollary} \label{exact-sequence} The split short exact sequence
\begin{equation} \label{split-short-seq-hol2r}
1 \rightarrow \Z_{2^r} \rightarrow \Hol(\Z_{2^r}) \buildrel \curvearrowleft \over \rightarrow \Aut(\Z_{2^r}) \rightarrow 1 
\end{equation}
implies that for $r\ge 3$
$$\Hol({\mathbb Z}_{2^r}) \cong \{x, \ y, \ z\ |\ x^{2^{r-2}}=y^2=z^{2^r}=1, \ xy=yx, \ zx=xz^3, \
zy=yz^{2^r-1}\}.$$

\end{corollary}

\section{Application of Wall's method (Generalization of his resolution \cite{W})}
\label{generalization-section}

Let $G$ be a split extension of a cyclic group by the direct sum of two cyclic groups:
$$1 \rightarrow {\mathbb Z}_q \rightarrow G \buildrel \curvearrowleft \over \rightarrow {\mathbb Z}_{s_1} \oplus {\mathbb Z}_{s_2}
\rightarrow 1.$$
More precisely, let $G$ be given by generators and relations:
$$G\cong \{ x, \ y, \ z \ |\ x^{s_1}=y^{s_2}=z^q=1, \ xy=yx, \ zx=xz^{t_1}, \ zy=yz^{t_2} \}$$
where $t_1^{s_1} \equiv 1(\mbox{mod }q), \ \ t_2^{s_2} \equiv 1(\mbox{mod }q).$

We will construct a free resolution for Hol$(G)$ with underlying module $A = \{A_{n,m}\}$ graded by
deg$(A_{n,m})=n+m$.

Let $B_m,$ $R_i',$ and $R_j''$ be the resolutions as above for the cyclic groups ${\mathbb Z}_q,$
${\mathbb Z}_{s_1},$ and ${\mathbb Z}_{s_2}$ on $z, \ x, \ y$ respectively.

Let the resolution for $H = {\mathbb Z}_{s_1} \oplus {\mathbb Z}_{s_2}$ be the tensor product of the
resolutions for ${\mathbb Z}_{s_1}$ and ${\mathbb Z}_{s_2}:$
$$C_n = \moplus_{i+j=n}{(R_i' \otimes_{\Z [H]} R_j'')} = \moplus_{i=0}^n{\Z [H]},$$
$$d_{n,i}=d_i' \otimes 1 \oplus (-1)^i 1 \otimes d_{n-i}''.$$
Then $C_n$ is free on $n+1$ generators. Define $\displaystyle A_{n,m}=\moplus_{n+1}\left( \Z[G]\otimes_{\Z[K]} B_m \right)$ 
(see \cite{W}).

\begin{lemma}\label{crazy_differentials_lemma}
Let $A_{n,m}\ (n \ge 0, m \ge 0)$ be free on $n+1$ generators, $a_{n,m,i},$ \mbox{$0\le i\le n$}, and let
the differentials be given by
$$\begin{array}{@{}rcl@{}}
d_0a_{n,2m+1,i} &=& (z-1)a_{n,2m,i},\\
d_0a_{n,2m,i} &=& N_z a_{n,2m,i},\\
d_1a_{2n+1,2m,2i} &=& \displaystyle (yL_{t_2,1}^m-1)a_{2n,2m,2i} +
\sum_{j=0}^{s_1-1}{x^jL_{t_1,j}^m}a_{2n,2m,2i-1},\\ 
\end{array}$$

$$\begin{array}{@{}rcl@{}}
d_1a_{2n+1,2m,2i+1} &=& \displaystyle -\sum_{j=0}^{s_2-1}{y^jL_{t_2,j}^m}a_{2n,2m,2i+1} +
(xL_{t_1,1}^m-1)a_{2n,2m,2i}, \\
d_1a_{2n+1,2m-1,2i} &=& \displaystyle -(yL_{t_2,1}^m-1)a_{2n,2m-1,2i} -
\sum_{j=0}^{s_1-1}{x^jL_{t_1,j}^m}a_{2n,2m-1,2i-1},\\
d_1a_{2n+1,2m-1,2i+1} &=& \displaystyle \sum_{j=0}^{s_2-1}{y^jL_{t_2,j}^m}a_{2n,2m-1,2i+1} -
(xL_{t_1,1}^m-1)a_{2n,2m-1,2i},\\
d_1a_{2n,2m,2i} &=& \displaystyle \sum_{j=0}^{s_2-1}{y^jL_{t_2,j}^m}a_{2n-1,2m,2i} +
\sum_{j=0}^{s_1-1}{x^jL_{t_1,j}^m}a_{2n-1,2m,2i-1},\\
d_1a_{2n,2m,2i+1} &=& -(yL_{t_2,1}^m-1)a_{2n-1,2m,2i+1} + (xL_{t_1,1}^m-1)a_{2n-1,2m,2i},\\
d_1a_{2n,2m-1,2i} &=& \displaystyle -\sum_{j=0}^{s_2-1}{y^jL_{t_2,j}^m}a_{2n-1,2m-1,2i} -
\sum_{j=0}^{s_1-1}{x^jL_{t_1,j}^m}a_{2n-1,2m-1,2i-1},\\
d_1a_{2n,2m-1,2i+1} &=& (yL_{t_2,1}^m-1)a_{2n-1,2m-1,2i+1} - (xL_{t_1,1}^m-1)a_{2n-1,2m-1,2i},\\
d_2a_{n,2m,i} &=& 0,\\
d_2a_{n,2m-1,i} &=& -{\frac{1}{q}}(t_2^{ms_2}-1)a_{n-2,2m,i} -
{\frac{1}{q}}(t_1^{ms_2}-1)a_{n-2,2m,i-2},
\end{array}$$
where $\displaystyle L_{f,j}=\sum_{i=0}^{f^j-1}{z^i}.$ Let $\displaystyle d=\sum_{k=0}^2 d_k$. Then $(A,d)$ is acyclic. 
\end{lemma}

\begin{proof}
Let us ``place'' this resolution in a three-dimensional lattice, with coordinates $x, \ y,$ and $z$, as follows.
We put each generator $a_{n,m,i}$ at the lattice point $(i,n-i,m)$, i.e. $i$ is its $x$-coordinate, $n$ is its total 
horizontal degree, 
and $m$ is its vertical degree. It is easy to see from our formulas that every nonzero differential is a linear 
combination of a differential lying  in a plane parallel to the $xz$-plane and a differential lying in a plane parallel 
to the $yz$-plane. The $xz$-plane contains Wall's resolution for 
the extension $1 \rightarrow \Z_q \rightarrow G^x \rightarrow \Z_{s_1} \rightarrow 1$ where 
$G^x\cong \{ x, \ z \ |\ x^{s_1}=z^q=1, zx=xz^{t_1}\}$ is the subgroup of $G$ generated by $x$ and $z$, 
and the $yz$-plane contains Wall's resolution for the extension 
$1 \rightarrow \Z_q \rightarrow G^y \rightarrow \Z_{s_2} \rightarrow 1$ where 
$G^y\cong \{ y, \ z \ |\ y^{s_2}=z^q=1, zy=yz^{t_2}\}$ is the subgroup of $G$ generated by $y$ and $z$. Moreover, all the 
differentials in every plane $y=c$ coincide with the corresponding differentials in the $xz$-plane, and the differentials 
in every plane $x=c$ coincide with the corresponding differentials in the $yz$-plane except for $d_1$ is multiplied by 
$(-1)^i$ where $i$ is the 
$x$-coordinate of the source of $d_1$. Notice that the vertical differentials $d_0$ are the same for both extensions. 
So let's write $d_1 = d_1^x + (-1)^i d_1^y$ and $d_2 = d_2^x + d_2^y$ where $d_j^x$ and $d_j^y$ are Wall's differentials for 
the above extensions. To show acyclicity, we have to show that for each $k$, $\displaystyle \sum_{j=0}^k d_j d_{k-j}=0$.

First of all, notice that all $d_2$'s are multiplications by numbers, all $d_1^x$'s and $d_1^y$'s are multiplications by 
expressions of $x$ and $y$ respectively, thus $d_2$ commute with everything, and $d_1^x$ commutes with $d_1^y$ as $x$ and 
$y$ commute in our group $G$.   

\begin{itemize}

\item For $k=0$ we have $d_0d_0=0$ because $B_m$ is a resolution for $\Z_q$. 

\item For $k=1$, since $d_0$ does not change the horizontal bidegree, 
$$\begin{array}{@{}l} 
d_0d_1 + d_1d_0 \\
		= d_0(d_1^x+(-1)^id_1^y) + (d_1^x+(-1)^id_1^y)d_0 \\
                = (d_0d_1^x+d_1^xd_0) + (-1)^i(d_0d_1^y+d_1^yd_0) \\
                = 0.
\end{array}$$ 
by \cite{W}. 

\item For $k=2$, since $d_1^x$ lowers $i$ by 1 and $d_1^y$ does not change $i$, 
$$\begin{array}{@{}l} 
d_0d_2 + d_1d_1 + d_2d_0 \\
    = d_0(d_2^x+d_2^y) + d_1(d_1^x+(-1)^id_1^y) + (d_2^x + d_2^y)d_0 \\
    = d_0d_2^x + d_0d_2^y + (d_1^x+(-1)^{i-1}d_1^y)d_1^x + (d_1^x+(-1)^id_1^y)(-1)^id_1^y + d_2^xd_0 + d_2^yd_0 \\ 
    = (d_0d_2^x + d_1^xd_1^x + d_2^xd_0) + (d_0d_2^y + d_1^yd_1^y + d_2^yd_0) + ((-1)^{i-1}d_1^yd_1^x + (-1)^id_1^xd_1^y) \\ 
    = 0.
\end{array}$$
The sums in the first two parenthesis here are 0 by \cite{W}, and the sum in the third parenthesis is identically 0. 

\item For $k=3$, since $d_j=0$ for $j\ge 3$, and $d_2^x$ lowers $i$ by 2 while $d_2^y$ does not change $i$, we have 
$$\begin{array}{@{}l}
d_0d_3 + d_1d_2 + d_2d_1 + d_3d_0  \\
=  d_1d_2 + d_2d_1 \\
= (d_1^x+(-1)^id_1^y)(d_2^x+d_2^y) + (d_2^x+d_2^y)(d_1^x+(-1)^id_1^y) \\
= (d_1^xd_2^x + d_2^xd_1^x)+(-1)^i(d_1^yd_2^y + d_2^yd_1^y) + (d_1^xd_2^y + d_2^yd_1^x) + (-1)^i(d_1^yd_2^x + d_2^xd_1^y)\\
= 0.
\end{array}$$
The sums in the first two parenthesis are again 0 by \cite{W}. The sum in the third parenthesis is 0 because of the 
following. $d_2$ raises 
the vertical degree by 1; if the vertical degree of the source is even, then $d_2=0$, while if the vertical degree of the 
source is odd, then $d_2^y$'s in this sum are multiplication by the same expression of $y$, and $d_1^x$'s are 
multiplication by expressions that differ by sign. The sum in the fourth parenthesis is 0 for a similar reason. 

\item For $k=4$ we have $d_2d_2=0$ because, as said above, $d_2=0$ raises the vertical degree by 1, and it is 0 whenever the 
vertical degree of the source is even. 

\item For $k\ge 5$, $\displaystyle \sum_{j=0}^k d_j d_{k-j}=0$ because $d_j=0$ for $j\ge 3$. 
\end{itemize}	
\end{proof} 

As in Wall's application (where $H$ is cyclic), on forming the tensor product with ${\mathbb Z}$, $z-1$
becomes zero, and the graded submodules
$$\begin{array}{rcl}
A^0 & = & \displaystyle \sum_{n\ge 0}{A_{n,0} \otimes_{\Z[G]} \Z},\\
A^m & = & \displaystyle \sum_{n\ge 0}{(A_{n,2m-1}+A_{n,2m}) \otimes_{\Z[G]} \Z}
\end{array}$$
are invariant under $d$. Thus $A \otimes_{\Z[G]} \Z$ is a direct sum of chain complexes. Hence it is sufficient to 
compute the homology of $d$ on each summand separately.

To compute the homology of $A^m$ for positive m, consider the spectral sequence of the filtration
$\displaystyle F_k(A^m) = \sum_{s\le k}A_{s,*}.$ The differential $d_0$ is induced by $N_z$, thus is multiplication by
$q.$ Hence the $E_1$-stage collapses to the \mbox{$(2m-1)$-row.}

\begin{remark}
The above construction can be generalized to the case of an arbitrary finite abelian group $H$ in the split short exact 
sequence
$$1 \rightarrow {\mathbb Z}_q \rightarrow G \msplit H \rightarrow 1.$$ 
If $\displaystyle H\cong \moplus_{i=1}^k{\mathbb Z}_{s_i}, $ let $A_{n,m}$ be
free on $\binom{n+k-1}{k-1}$ generators, $a_{n,m,i_1,\ldots i_k},$ $0\le i_j\le n,$
$\displaystyle \sum_{j=1}^k{i_j}=n.$ Place each generator $a_{n,m,i_1,\ldots i_k}$ at the lattice point 
$(i_1, i_2, \ldots, i_k, m)$. The differentials are as follows. Construct Wall's resolution for 
\begin{equation}\label{i_th_resolution}
1 \rightarrow {\mathbb Z}_q \rightarrow G_i \msplit \Z_{s_i} \rightarrow 1
\end{equation}
for each $i$ from 1 to $k$, and put it in each plane parallel to the $x_i z$-plane in the $(k+1)$-dimensional 
$x_1 x_2 \ldots x_k z$-space. Let $d_1^{x_i}$ and $d_2^{x_i}$ be differentials in Wall's resolution for 
\eqref{i_th_resolution}. Then define $d_1$ and $d_2$ by 
$$\begin{array}{rcl}
d_1 & = & d_1^{x_1} + (-1)^{i_1}d_1^{x_2} + (-1)^{i_1+i_2}d_1^{x_3} + \ldots + (-1)^{i_1+i_2+\ldots+i_{k-1}}d_1^{x_k}, \\
d_2 & = & d_2^{x_1} + d_2^{x_2} + \ldots + d_2^{x_k}, 
\end{array}$$
where $i_j$ is the $j$-th coordinate of the source of the differential. 
Define $d_0$ as in each resolution \eqref{i_th_resolution}. Let 
$$d=d_0 + d_1 + d_2.$$
Then $(A,d)$ is a resolution (the proof is a slight generalization of the proof of lemma \ref{crazy_differentials_lemma}). 
The graded submodules 
$$\begin{array}{rcl}
A^0 & = & \displaystyle \sum_{i_1+\ldots+i_k=n\ge 0}{A_{n,0,i_1,\ldots,i_k} \otimes_{\Z[G]} \Z},\\
A^m & = & \displaystyle \sum_{i_1+\ldots+i_k=n\ge 0}{(A_{n,2m-1,i_1,\ldots,i_k}+A_{n,2m,i_1,\ldots,i_k}) \otimes_{\Z[G]} \Z}
\end{array}$$
are still invariant under $d,$ however, the larger $k$ is, the more complicated is the 
calculation of $H_*(A^m)$.
\end{remark}

\section{$H_*(\Hol(\Z_{2^r}); \Z)$}

The application of the method described above to computation of
$H_*(\Hol({\mathbb Z}_{2^r}); {\mathbb Z})$
is straightforward: according to corollary \ref{exact-sequence}, we set $q=2^r,$  $s_1=2^{r-2},$
$s_2=2,$  $t_1=3,$ and $t_2=2^r-1.$

Let $\{ f_{n,i}, \ 0\le i\le n\}$ be the basis of the $E^1$-stage of the spectral sequence for $A^m$
deduced from the basis $\{ a_{n,2m-1,i}, 0\le i\le n\}.$ Then for
$m=0,$
$$\begin{array}{rcl}
d_1f_{2n+1,2i} &=& 2^{r-2}f_{2n,2i-1},\\
d_1f_{2n+1,2i+1} &=& -2f_{2n,2i+1},\\   
d_1f_{2n,2i} &=& 2f_{2n-1,2i} + 2^{r-2}f_{2n-1,2i-1},\\
d_1f_{2n,2i+1} &=& 0.
\end{array}$$

Then
$$\begin{array}{rcl}
H_{2n+1}(A^0 ; {\mathbb Z}) &\cong & \displaystyle (\moplus_{n-1}{{\mathbb Z}_2})\moplus {\mathbb Z}_{2^{r-2}} \\
H_{2n}(A^0 ; {\mathbb Z}) &\cong & \displaystyle \moplus_n{{\mathbb Z}_2}, \ n\ge 1.
\end{array}$$

For $m\ge 1,$
\begin{equation} \label{2-differentials} \begin{array}{rcl}
d_1 f_{2n+1,2i} &=& \displaystyle -((-1)^m-1)f_{2n,2i} - {\frac{3^{2^{r-2}m}-1}{3^m-1}}f_{2n,2i-1},\\
d_1 f_{2n+1,2i+1} &=& \displaystyle (1+(-1)^m)f_{2n,2i+1} - (3^m-1)f_{2n,2i},\\
d_1 f_{2n,2i} &=& \displaystyle -(1+(-1)^m)f_{2n-1,2i} - {\frac{3^{2^{r-2}m}-1}{3^m-1}}f_{2n-1,2i-1},\\
d_1 f_{2n,2i+1} &=& ((-1)^m-1)f_{2n-1,2i+1} - (3^m-1)f_{2n-1,2i}.
\end{array} \end{equation}

By lemma \ref{obvious}, to compute $H_*(A^m; {\mathbb Z}),$ it suffices to know the largest powers of 2 that
divide the coefficients in \eqref{2-differentials}.

By lemma \ref{divisible-2},
\begin{itemize}
\item
if $m$ is odd, then $\nu_2(3^{2^{r-2}\cdot m}-1) = r$ and $\nu_2(3^m-1) = 1$, therefore \\ 
$\nu_2 \left( {\frac{3^{2^{r-2} m}-1}{3^m-1}} \right) = r-1,$
\item
if \,$\nu_2(m)=q \ge 1$, \,then \,$\nu_2(3^{2^{r-2} m}-1) = r+q$ \,and \,
$\nu_2(3^m-1) = q+2$, \,therefore  $\nu_2 \left( {\frac{3^{2^{r-2}\cdot m}-1}{3^m-1}} \right) = r-2.$
\end{itemize}
It follows that:
\begin{itemize}
\item
if $m$ is odd, then $H_{n}(A^m ; {\mathbb Z}) = {\mathbb Z}_2$ for any $n\ge 0,$
\item
if $V_2(m)=q,$ $1\le q\le r-2$, then
$$\begin{array}{rcl}
H_{0}(A^m ; {\mathbb Z}) &\cong & {\mathbb Z}_{2^{q+2}} \\
H_{n}(A^m ; {\mathbb Z}) &\cong & \displaystyle \left( \moplus_n{\mathbb Z}_2 \right) \moplus {\mathbb Z}_{2^q} \ \ (n\ge 1)
\end{array}$$
\item
If $q\ge r-2$, then $2^{q+2}\equiv 0 (\mbox{mod}2^r)$, therefore the homology of $A^{2^q\cdot l}$ does not
depend on $q$, and equals that of $A^{2^{r-2}}.$
\end{itemize}

\begin{figure}
\hspace{0.9cm} \input{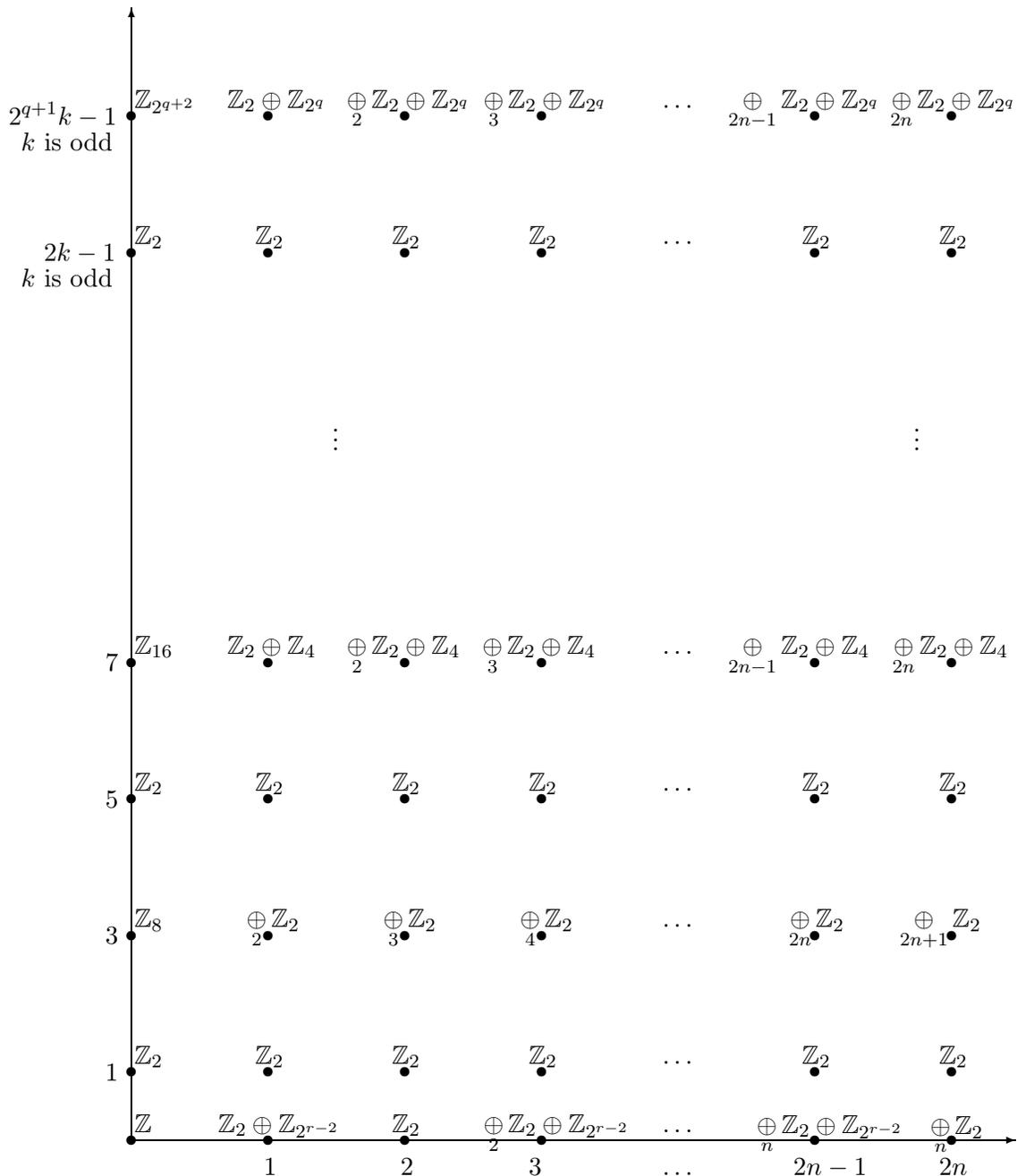}
\caption{The first $2^{r-1}$ rows of the $E^2=E^{\infty}$ stage of the spectral sequence for $H_*(\Hol(\Z_{2^r}); \Z)$ induced by the 
filtration $F$ of $A$, $r\ge 3$}
\label{E-infty_stage}
\end{figure}

The $E^{\infty}$ stage of this spectral sequence is shown in figure \ref{E-infty_stage}. 

Finally, $H_*(\mbox{Hol}({\mathbb Z}_{2^r}); {\mathbb Z})$ is the direct sum of the above for all $m.$ We have
$$H_q(\mbox{Hol}({\mathbb Z}_8); {\mathbb Z}) \cong \moplus_{i=1}^3{ \left( \moplus_{n_i}{\mathbb Z}_{2^i} \right) },$$
where

\bigskip

$\begin{array}{rcl}
n_1 & = & \displaystyle {\frac{q^2}{8}}+{\frac{3q}{4}} + \left\{
   \begin{array}{rl}
   0    & \mbox{if } q\equiv 0(\mbox{mod }4)\\
   17/8  & \mbox{if } q\equiv 1(\mbox{mod }4)\\
   -1/4  & \mbox{if } q\equiv 2(\mbox{mod }4)\\
   5/8 & \mbox{if } q\equiv 3(\mbox{mod }4)
   \end{array} \right.\\
 & & \\
 &=& \left\{
   \begin{array}{ll}
   2n^2+3n  & \mbox{if } q=4n\\
   2n^2+4n+3  & \mbox{if } q=4n+1\\
   2n^2+5n+2 & \mbox{if } q=4n+2\\
   2n^2+6n+4    & \mbox{if } q=4n+3
   \end{array} \right.\\
n_2 & = & 0\\
n_3 & = & \left\{
   \begin{array}{rl}
   1 & \mbox{if } q\equiv 3(\mbox{mod }4)\\
   0 & \mbox{otherwise}
   \end{array} \right.\\
\end{array}$

$$H_q(\mbox{Hol}({\mathbb Z}_{16}); {\mathbb Z}) \cong \moplus_{i=1}^4{ \left( \moplus_{n_i}{\mathbb Z}_{2^i} \right) },$$
where

$\begin{array}{@{}rcl@{}}
n_1 & = & \displaystyle {\frac{q^2}{8}}-{\frac{q}{4}}+ \left[ {\frac{q}{4}}\right] -\left[
{\frac{q}{8}}\right] +\left[ \frac{q+3}{4} \right] +\left[ {\frac{q+1}{2}}\right] +
\left\{
   \begin{array}{@{}rl@{}}
   1/8  & \mbox{if } q\equiv 1(\mbox{mod }4)\\
   -3/8 & \mbox{if } q\equiv 3(\mbox{mod }4)\\
   0    & \mbox{otherwise}
   \end{array} \right.\\
 & & \\
n_2 & = & \displaystyle \left[ {\frac{q}{8}}\right] + \left\{
   \begin{array}{rl}
   1 & \mbox{if } q\equiv 1(\mbox{mod }2)\\
   0 & \mbox{otherwise}
   \end{array} \right.\\
 & & \\ 
n_3 & = & \left\{
   \begin{array}{rl}   
   1 & \mbox{if } q\equiv 3(\mbox{mod }8)\\
   0 & \mbox{otherwise}
   \end{array} \right.\\
 & & \\
n_4 & = & \left\{
   \begin{array}{rl}
   1 & \mbox{if } q\equiv 7(\mbox{mod }8)\\
   0 & \mbox{otherwise}
   \end{array} \right.\\  
\end{array}$

\vspace{1cm}

If $r\ge 5,$
$$H_q(\mbox{Hol}({\mathbb Z}_{2^r}); {\mathbb Z}) \cong \moplus_{i=1}^r{ \left( \moplus_{n_i}{\mathbb Z}_{2^i} \right) },$$
where
   
$\begin{array}{@{}l@{}cl@{}}
n_1 & = & \displaystyle {\frac{q^2}{8}}-{\frac{q}{4}}+\left[ {\frac{q}{4}}\right] -\left[
{\frac{q}{8}}\right] +\left[ {\frac{q+3}{4}}\right] +\left[ {\frac{q+1}{2}}\right] +
\left\{
   \begin{array}{@{}rl@{}}
   1/8  & \mbox{if } q\equiv 1(\mbox{mod }4)\\
   -3/8 & \mbox{if } q\equiv 3(\mbox{mod }4)\\
   0    & \mbox{otherwise}
   \end{array} \right.\\
 & & \\
n_2 & = & \displaystyle \left[ {\frac{q}{8}}\right] -\left[ {\frac{q}{16}}\right] \\
 & & \\
n_3 & = & \displaystyle \left[ {\frac{q}{16}}\right] -\left[ {\frac{q}{32}}\right] + \left\{
   \begin{array}{rl}
   1  & \mbox{if } q\equiv 3(\mbox{mod }8)\\
   0  & \mbox{otherwise}
   \end{array} \right.\\  
    &   & \ldots\\
\end{array}$

$\begin{array}{@{}l@{}cl@{}} 
n_i & = & \displaystyle \left[ {\frac{q}{2^{i+1}}}\right] -\left[ 
{\frac{q}{2^{i+2}}}\right] + \left\{
   \begin{array}{rl}
   1  & \mbox{if } q\equiv 2^{i-1}-1(\mbox{mod }2^i)\\
   0  & \mbox{otherwise}
   \end{array} \right.\\
    &   & (3\le i \le r-3)\\
    &   & \ldots\\
n_{r-2} & = & \displaystyle \left[ {\frac{q}{2^{r-1}}}\right] + \left\{
   \begin{array}{rl}
   1 & \mbox{if } q\equiv 1(\mbox{mod }2)\\
   1 & \mbox{if } q\equiv 2^{r-3}-1(\mbox{mod }2^{r-2})\\
   0 & \mbox{otherwise}   
   \end{array} \right.\\
 & & \\
n_{r-1} & = & \left\{
   \begin{array}{rl}
   1 & \mbox{if } q\equiv 2^{r-2}-1(\mbox{mod }2^{r-1})\\
   0 & \mbox{otherwise}
   \end{array} \right. \\
 & & \\
n_r & = & \left\{   
   \begin{array}{rl}
   1 & \mbox{if } q\equiv -1(\mbox{mod }2^{r-1})\\
   0 & \mbox{otherwise}   
   \end{array} \right. \\
\end{array}$

\section{$H^*(\Hol(\Z_{2^r});\F_2)$}

First we use the universal coefficient theorem to compute the mod 2 cohomology of Hol$({\mathbb Z}_{2^r})$
additively. 

\begin{proposition} If $r\ge 3,$ then the ranks of
$\displaystyle H^{*}(\Hol({\mathbb Z}_{2^r}); {\mathbb F}_2)$ do not depend on $r:$

$$\begin{array}{rccl}
H^{4n}(\Hol({\mathbb Z}_{2^r}); {\mathbb F}_2)   & \cong & \displaystyle \moplus_{4n^2+5n+1}  & {{\mathbb F}_2}, \\
H^{4n+1}(\Hol({\mathbb Z}_{2^r}); {\mathbb F}_2) & \cong & \displaystyle \moplus_{4n^2+7n+3}  & {{\mathbb F}_2}, \\
H^{4n+2}(\Hol({\mathbb Z}_{2^r}); {\mathbb F}_2) & \cong & \displaystyle \moplus_{4n^2+9n+5}  & {{\mathbb F}_2}, \\
H^{4n+3}(\Hol({\mathbb Z}_{2^r}); {\mathbb F}_2) & \cong & \displaystyle \moplus_{4n^2+11n+7} & {{\mathbb F}_2}
\end{array}$$
\end{proposition}

Next we find the ring structure of $H^*(\mbox{Hol}({\mathbb Z}_{2^r}); {\mathbb F}_2)$ using LHS spectral
sequence. Consider the following three subgroups of
$\mbox{Aut}({\mathbb Z}_{2^r}) \cong {\mathbb Z}_{2^{r-2}} \oplus {\mathbb Z}_2:$
$$\begin{array}{llccc}   
1. & A_1 = & {\mathbb Z}_{2^{r-2}} & \hookrightarrow & {\mathbb Z}_{2^{r-2}} \oplus {\mathbb Z}_2 \\
   &       & 1                     & \mapsto         & (1,0) \\
2. & A_2 = & {\mathbb Z}_2         & \hookrightarrow & {\mathbb Z}_{2^{r-2}} \oplus {\mathbb Z}_2 \\
   &       & 1                     & \mapsto         & (0,1) \\
3. & A_3 = & {\mathbb Z}_2         & \hookrightarrow & {\mathbb Z}_{2^{r-2}} \oplus {\mathbb Z}_2 \\
   &       & 1                     & \mapsto         & (2^{r-3}, 1)
\end{array}$$
and the corresponding pullbacks:
$$\begin{CD}
1  @>>> \Z_{2^r} @>>> G_i            @>>> A_i                     @>>> 1  \\
@.      @VVV          @VVV                @VVV                         @. \\
1  @>>> \Z_{2^r} @>>> \Hol(\Z_{2^r}) @>>> \Z_{2^{r-2}} \oplus\Z_2 @>>> 1   
\end{CD}$$

Then $G_i$ is a split extension of cyclic by cyclic groups, so we can use Wall's method to compute its integer
homology, and then use the universal coefficients theorem to compute its mod 2 cohomology.

\begin{remark}
$G_2$ and $G_3$ are the dihedral and semi-dihedral \cite{MP} groups of order $2^{r+1}$.
\end{remark}

Let $A_1'$ be the subgroup of $A_1$ of order 2. Then the corresponding pullback is the quasi-dihedral group of order 
$2^{r+1}$:
$$ \{ z, s | z^{2^r}=1, s^2=1, zs = sz^{2^{r-1}+1} \}$$

\begin{figure}
\centerline{\setlength{\unitlength}{3400sp}%
\begingroup\makeatletter\ifx\SetFigFont\undefined%
\gdef\SetFigFont#1#2#3#4#5{%
  \reset@font\fontsize{#1}{#2pt}%
  \fontfamily{#3}\fontseries{#4}\fontshape{#5}%
  \selectfont}%
\fi\endgroup%
\begin{picture}(5400,3205)(250,-3500)
\thinlines
\special{ps: gsave 0 0 0 setrgbcolor}\put(1201,-2761){\vector( 1, 0){2400}}
\special{ps: grestore}\special{ps: gsave 0 0 0 setrgbcolor}\put(1201,-2761){\vector( 0, 1){2400}}
\special{ps: grestore}\special{ps: gsave 0 0 0 setrgbcolor}\put(1126,-2461){\line( 1, 0){ 75}}
\special{ps: grestore}\special{ps: gsave 0 0 0 setrgbcolor}\put(1126,-2161){\line( 1, 0){ 75}}
\special{ps: grestore}\special{ps: gsave 0 0 0 setrgbcolor}\put(1126,-1861){\line( 1, 0){ 75}}
\special{ps: grestore}\special{ps: gsave 0 0 0 setrgbcolor}\put(1126,-1561){\line( 1, 0){ 75}}
\special{ps: grestore}\special{ps: gsave 0 0 0 setrgbcolor}\put(1126,-1261){\line( 1, 0){ 75}}
\special{ps: grestore}\special{ps: gsave 0 0 0 setrgbcolor}\put(1126,-961){\line( 1, 0){ 75}}
\special{ps: grestore}\special{ps: gsave 0 0 0 setrgbcolor}\put(1651,-2761){\line( 0,-1){ 75}}
\special{ps: grestore}\special{ps: gsave 0 0 0 setrgbcolor}\put(2101,-2761){\line( 0,-1){ 75}}
\special{ps: grestore}\special{ps: gsave 0 0 0 setrgbcolor}\put(2551,-2761){\line( 0,-1){ 75}}
\special{ps: grestore}\special{ps: gsave 0 0 0 setrgbcolor}\put(3001,-2761){\line( 0,-1){ 75}}
\special{ps: grestore}\special{ps: gsave 0 0 0 setrgbcolor}\put(1201,-2161){\vector( 3,-1){900}}
\special{ps: grestore}\special{ps: gsave 0 0 0 setrgbcolor}\put(1651,-2161){\vector( 3,-1){900}}
\special{ps: grestore}\special{ps: gsave 0 0 0 setrgbcolor}\put(2101,-2161){\vector( 3,-1){900}}
\special{ps: grestore}\special{ps: gsave 0 0 0 setrgbcolor}\put(1201,-961){\vector( 3,-1){900}}
\special{ps: grestore}\special{ps: gsave 0 0 0 setrgbcolor}\put(1651,-961){\vector( 3,-1){900}}
\special{ps: grestore}\special{ps: gsave 0 0 0 setrgbcolor}\put(2101,-961){\vector( 3,-1){900}}
\special{ps: grestore}\special{ps: gsave 0 0 0 setrgbcolor}\put(1201,-661){\line(-1, 0){ 75}}
\special{ps: 
grestore}\put(5000,-1900){\makebox(0,0)[lb]{\smash{\SetFigFont{11}{13.2}{\rmdefault}{\mddefault}{\updefault}\special{ps: 
gsave 0 0 0 setrgbcolor}$r>3$\special{ps: grestore}}}}
\put(1096,-2516){\makebox(0,0)[rb]{\smash{\SetFigFont{11}{13.2}{\rmdefault}{\mddefault}{\updefault}\special{ps: 
gsave 0 0 0 setrgbcolor}{\small $a_1=u_1$}\special{ps: grestore}}}}
\put(1096,-2216){\makebox(0,0)[rb]{\smash{\SetFigFont{11}{13.2}{\rmdefault}{\mddefault}{\updefault}\special{ps: 
gsave 0 0 0 setrgbcolor}{\small $v_1$}\special{ps: grestore}}}}
\put(1096,-1916){\makebox(0,0)[rb]{\smash{\SetFigFont{11}{13.2}{\rmdefault}{\mddefault}{\updefault}\special{ps: 
gsave 0 0 0 setrgbcolor}{\small $b_1=u_1v_1$}\special{ps: grestore}}}}
\put(1096,-1616){\makebox(0,0)[rb]{\smash{\SetFigFont{11}{13.2}{\rmdefault}{\mddefault}{\updefault}\special{ps: 
gsave 0 0 0 setrgbcolor}{\small $c_1=v_1^2$}\special{ps: grestore}}}}
\put(1096,-1316){\makebox(0,0)[rb]{\smash{\SetFigFont{11}{13.2}{\rmdefault}{\mddefault}{\updefault}\special{ps: 
gsave 0 0 0 setrgbcolor}{\small $u_1v_1^2$}\special{ps: grestore}}}}
\put(1096,-1016){\makebox(0,0)[rb]{\smash{\SetFigFont{11}{13.2}{\rmdefault}{\mddefault}{\updefault}\special{ps: 
gsave 0 0 0 setrgbcolor}{\small $v_1^3$}\special{ps: grestore}}}}
\put(1096,-716){\makebox(0,0)[rb]{\smash{\SetFigFont{11}{13.2}{\rmdefault}{\mddefault}{\updefault}\special{ps: 
gsave 0 0 0 setrgbcolor}{\small $u_1v_1^3$}\special{ps: grestore}}}}
\put(1576,-2986){\makebox(0,0)[lb]{\smash{\SetFigFont{11}{13.2}{\rmdefault}{\mddefault}{\updefault}\special{ps: 
gsave 0 0 0 setrgbcolor}{\small $x_1$}\special{ps: grestore}}}}
\put(2026,-2986){\makebox(0,0)[lb]{\smash{\SetFigFont{11}{13.2}{\rmdefault}{\mddefault}{\updefault}\special{ps: 
gsave 0 0 0 setrgbcolor}{\small $z_1$}\special{ps: grestore}}}}
\put(2401,-2986){\makebox(0,0)[lb]{\smash{\SetFigFont{11}{13.2}{\rmdefault}{\mddefault}{\updefault}\special{ps: 
gsave 0 0 0 setrgbcolor}{\small $x_1z_1$}\special{ps: grestore}}}}
\put(2926,-2986){\makebox(0,0)[lb]{\smash{\SetFigFont{11}{13.2}{\rmdefault}{\mddefault}{\updefault}\special{ps: 
gsave 0 0 0 setrgbcolor}{\small $z_1^2$}\special{ps: grestore}}}}
\put(1276,-451){\makebox(0,0)[lb]{\smash{\SetFigFont{11}{13.2}{\rmdefault}{\mddefault}{\updefault}\special{ps: 
gsave 0 0 0 setrgbcolor}{\small $H^*({\mathbb Z}_{2^r};{\mathbb F}_2)$}\special{ps: grestore}}}}
\put(3490,-2686){\makebox(0,0)[lb]{\smash{\SetFigFont{11}{13.2}{\rmdefault}{\mddefault}{\updefault}\special{ps: 
gsave 0 0 0 setrgbcolor}{\small $H^*({\mathbb Z}_{2^{r-2}};{\mathbb F}_2)$}\special{ps: grestore}}}}
\put(2551,-2236){\makebox(0,0)[lb]{\smash{\SetFigFont{11}{13.2}{\rmdefault}{\mddefault}{\updefault}\special{ps: 
gsave 0 0 0 setrgbcolor}{\small $d_2$}\special{ps: grestore}}}}
\end{picture}}
\centerline{\setlength{\unitlength}{3400sp}%
\begingroup\makeatletter\ifx\SetFigFont\undefined%
\gdef\SetFigFont#1#2#3#4#5{%
  \reset@font\fontsize{#1}{#2pt}%
  \fontfamily{#3}\fontseries{#4}\fontshape{#5}%
  \selectfont}%
\fi\endgroup%
\begin{picture}(5400,3205)(4840,-3500)
\thinlines
\special{ps: gsave 0 0 0 setrgbcolor}\put(5801,-2761){\vector( 0, 1){2400}}
\special{ps: grestore}\special{ps: gsave 0 0 0 setrgbcolor}\put(5801,-2761){\vector( 1, 0){2400}}
\special{ps: grestore}\special{ps: gsave 0 0 0 setrgbcolor}\put(6251,-2761){\line( 0,-1){ 75}}
\special{ps: grestore}\special{ps: gsave 0 0 0 setrgbcolor}\put(6701,-2761){\line( 0,-1){ 75}}
\special{ps: grestore}\special{ps: gsave 0 0 0 setrgbcolor}\put(7151,-2761){\line( 0,-1){ 75}}
\special{ps: grestore}\special{ps: gsave 0 0 0 setrgbcolor}\put(7601,-2761){\line( 0,-1){ 75}}
\special{ps: grestore}\special{ps: gsave 0 0 0 setrgbcolor}\put(5801,-2461){\line(-1, 0){ 75}}
\special{ps: grestore}\special{ps: gsave 0 0 0 setrgbcolor}\put(5801,-2161){\line(-1, 0){ 75}}
\special{ps: grestore}\special{ps: gsave 0 0 0 setrgbcolor}\put(5801,-1861){\line(-1, 0){ 75}}
\special{ps: grestore}\special{ps: gsave 0 0 0 setrgbcolor}\put(5801,-1561){\line(-1, 0){ 75}}
\special{ps: grestore}\special{ps: gsave 0 0 0 setrgbcolor}\put(5801,-1261){\line(-1, 0){ 75}}
\special{ps: grestore}\special{ps: gsave 0 0 0 setrgbcolor}\put(5801,-961){\line(-1, 0){ 75}}
\special{ps: grestore}\special{ps: gsave 0 0 0 setrgbcolor}\put(5801,-661){\line(-1, 0){ 75}}
\special{ps: grestore}\special{ps: gsave 0 0 0 setrgbcolor}\put(5801,-2161){\vector( 3,-1){900}}
\special{ps: grestore}\special{ps: gsave 0 0 0 setrgbcolor}\put(6251,-2161){\vector( 3,-1){900}}
\special{ps: grestore}\special{ps: gsave 0 0 0 setrgbcolor}\put(6701,-2161){\vector( 3,-1){900}}
\special{ps: grestore}\special{ps: gsave 0 0 0 setrgbcolor}\put(5801,-961){\vector( 3,-1){900}}
\special{ps: grestore}\special{ps: gsave 0 0 0 setrgbcolor}\put(6251,-961){\vector( 3,-1){900}}
\special{ps: grestore}\special{ps: gsave 0 0 0 setrgbcolor}\put(6701,-961){\vector( 3,-1){900}}
\special{ps: grestore}
\put(9600,-1900){\makebox(0,0)[lb]{\smash{\SetFigFont{11}{13.2}{\rmdefault}{\mddefault}{\updefault}\special{ps: 
gsave 0 0 0 setrgbcolor}$r=3$\special{ps: grestore}}}}
\put(5696,-716){\makebox(0,0)[rb]{\smash{\SetFigFont{11}{13.2}{\rmdefault}{\mddefault}{\updefault}\special{ps: 
gsave 0 0 0 setrgbcolor}{\small $u_1v_1^3$}\special{ps: grestore}}}}
\put(5696,-1316){\makebox(0,0)[rb]{\smash{\SetFigFont{11}{13.2}{\rmdefault}{\mddefault}{\updefault}\special{ps: 
gsave 0 0 0 setrgbcolor}{\small $u_1v_1^2$}\special{ps: grestore}}}}
\put(5696,-1016){\makebox(0,0)[rb]{\smash{\SetFigFont{11}{13.2}{\rmdefault}{\mddefault}{\updefault}\special{ps: 
gsave 0 0 0 setrgbcolor}{\small $v_1^3$}\special{ps: grestore}}}}
\put(5696,-1616){\makebox(0,0)[rb]{\smash{\SetFigFont{11}{13.2}{\rmdefault}{\mddefault}{\updefault}\special{ps: 
gsave 0 0 0 setrgbcolor}{\small $c_1=v_1^2$}\special{ps: grestore}}}}
\put(5696,-1916){\makebox(0,0)[rb]{\smash{\SetFigFont{11}{13.2}{\rmdefault}{\mddefault}{\updefault}\special{ps: 
gsave 0 0 0 setrgbcolor}{\small $b_1=u_1v_1$}\special{ps: grestore}}}}
\put(5696,-2216){\makebox(0,0)[rb]{\smash{\SetFigFont{11}{13.2}{\rmdefault}{\mddefault}{\updefault}\special{ps: 
gsave 0 0 0 setrgbcolor}{\small $v_1$}\special{ps: grestore}}}}
\put(5696,-2516){\makebox(0,0)[rb]{\smash{\SetFigFont{11}{13.2}{\rmdefault}{\mddefault}{\updefault}\special{ps: 
gsave 0 0 0 setrgbcolor}{\small $a_1=u_1$}\special{ps: grestore}}}}
\put(6176,-2986){\makebox(0,0)[lb]{\smash{\SetFigFont{11}{13.2}{\rmdefault}{\mddefault}{\updefault}\special{ps: 
gsave 0 0 0 setrgbcolor}{\small $x_1$}\special{ps: grestore}}}}
\put(6626,-2986){\makebox(0,0)[lb]{\smash{\SetFigFont{11}{13.2}{\rmdefault}{\mddefault}{\updefault}\special{ps: 
gsave 0 0 0 setrgbcolor}{\small $x_1^2$}\special{ps: grestore}}}}
\put(7076,-2986){\makebox(0,0)[lb]{\smash{\SetFigFont{11}{13.2}{\rmdefault}{\mddefault}{\updefault}\special{ps: 
gsave 0 0 0 setrgbcolor}{\small $x_1^3$}\special{ps: grestore}}}}
\put(7526,-2986){\makebox(0,0)[lb]{\smash{\SetFigFont{11}{13.2}{\rmdefault}{\mddefault}{\updefault}\special{ps: 
gsave 0 0 0 setrgbcolor}{\small $x_1^4$}\special{ps: grestore}}}}
\put(5876,-451){\makebox(0,0)[lb]{\smash{\SetFigFont{11}{13.2}{\rmdefault}{\mddefault}{\updefault}\special{ps: 
gsave 0 0 0 setrgbcolor}{\small $H^*({\mathbb Z}_8;{\mathbb F}_2)$}\special{ps: grestore}}}}
\put(8090,-2686){\makebox(0,0)[lb]{\smash{\SetFigFont{11}{13.2}{\rmdefault}{\mddefault}{\updefault}\special{ps: 
gsave 0 0 0 setrgbcolor}{\small $H^*({\mathbb Z}_2;{\mathbb F}_2)$}\special{ps: grestore}}}}
\put(7151,-2236){\makebox(0,0)[lb]{\smash{\SetFigFont{11}{13.2}{\rmdefault}{\mddefault}{\updefault}\special{ps: 
gsave 0 0 0 setrgbcolor}{\small $d_2$}\special{ps: grestore}}}}
\end{picture}}
\caption{$E_2$ stage of the LHS spectral sequence for $H^*(G_1; {\mathbb F}_2)$}
\label{Lyndon_G_1}
\end{figure}
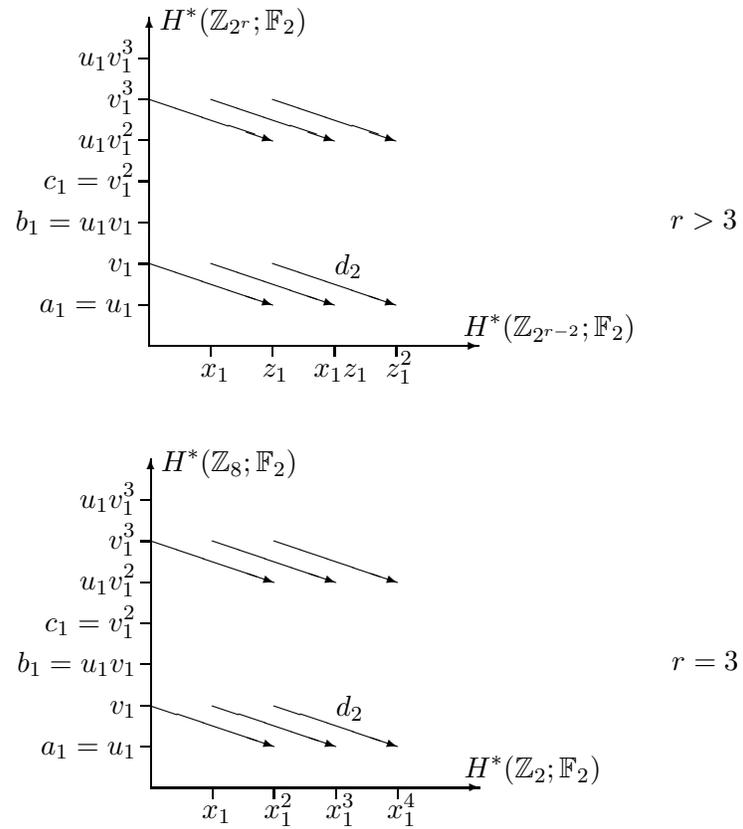

Now consider the LHS spectral sequence for $H^*(G_1; \F_2)$. Its $E_2$ stage is shown in Figure \ref{Lyndon_G_1}. 

Using Wall's formulas we have $H_1(G_1;\Z)\cong \Z_{2^{r-2}}\oplus \Z_2$ and $H_2(G_1;\Z)=0$. 
Then $H^2(G_1;\Z_2)\cong \Z_2\oplus \Z_2$. It follows  that
$$d_2(v_1)=u_1z_1 \mbox{ \ \ if } r>3,$$
$$d_2(v_1)=u_1x_1^2 \mbox{ \ \ if } r=3.$$
Then
$$d_2(v_1^{2n+1} x_1^i z_1^j)=u_1 v_1 ^{2n} x_1^i z_1^{j+1} \mbox{ \ \ if } r>3,$$
$$d_2(v_1^{2n+1} x_1^i)=u_1 v_1^{2n} x_1^{i+2} \mbox{ \ \ if } r=3,$$
and all other differentials are zero.

The spectral sequence for $H^*(\Hol(\Z_{2^r}); {\mathbb F}_2)$ shown in figure \ref{Lyndon_Hol_Z2r}
maps into the spectral sequence for $H^*(G_1; \F_2)$ shown in figure \ref{Lyndon_G_1}; $x$ is mapped to $x_1$, $z$
to $z_1$, $y$ to 0, $u$ to $u_1$, and $v$ to $v_1$. Therefore,
$$d_2(v)=uz+ \mbox{possibly some other term(s) \ \ \ if } r>3,$$
$$d_2(v)=ux^2+ \mbox{possibly some other term(s) \ \ \ if } r>3.$$

\begin{figure}
\centerline{\setlength{\unitlength}{3400sp}%
\begingroup\makeatletter\ifx\SetFigFont\undefined%
\gdef\SetFigFont#1#2#3#4#5{%
  \reset@font\fontsize{#1}{#2pt}%
  \fontfamily{#3}\fontseries{#4}\fontshape{#5}%
  \selectfont}%
\fi\endgroup%
\begin{picture}(6000,3905)(600,-4236)
\thinlines
\special{ps: gsave 0 0 0 setrgbcolor}\put(1201,-2761){\vector( 1, 0){2400}}
\special{ps: grestore}\special{ps: gsave 0 0 0 setrgbcolor}\put(1201,-2761){\vector( 0, 1){2400}}
\special{ps: grestore}\special{ps: gsave 0 0 0 setrgbcolor}\put(1126,-2461){\line( 1, 0){ 75}}
\special{ps: grestore}\special{ps: gsave 0 0 0 setrgbcolor}\put(1126,-2161){\line( 1, 0){ 75}}
\special{ps: grestore}\special{ps: gsave 0 0 0 setrgbcolor}\put(1126,-1861){\line( 1, 0){ 75}}
\special{ps: grestore}\special{ps: gsave 0 0 0 setrgbcolor}\put(1126,-1561){\line( 1, 0){ 75}}
\special{ps: grestore}\special{ps: gsave 0 0 0 setrgbcolor}\put(1126,-1261){\line( 1, 0){ 75}}
\special{ps: grestore}\special{ps: gsave 0 0 0 setrgbcolor}\put(1126,-961){\line( 1, 0){ 75}}
\special{ps: grestore}\special{ps: gsave 0 0 0 setrgbcolor}\put(1651,-2761){\line( 0,-1){ 75}}
\special{ps: grestore}\special{ps: gsave 0 0 0 setrgbcolor}\put(2101,-2761){\line( 0,-1){ 75}}
\special{ps: grestore}\special{ps: gsave 0 0 0 setrgbcolor}\put(2551,-2761){\line( 0,-1){ 75}}
\special{ps: grestore}\special{ps: gsave 0 0 0 setrgbcolor}\put(3001,-2761){\line( 0,-1){ 75}}
\special{ps: grestore}\special{ps: gsave 0 0 0 setrgbcolor}\put(1201,-2161){\vector( 3,-1){900}}
\special{ps: grestore}\special{ps: gsave 0 0 0 setrgbcolor}\put(1651,-2161){\vector( 3,-1){900}}
\special{ps: grestore}\special{ps: gsave 0 0 0 setrgbcolor}\put(2101,-2161){\vector( 3,-1){900}}
\special{ps: grestore}\special{ps: gsave 0 0 0 setrgbcolor}\put(1201,-961){\vector( 3,-1){900}}
\special{ps: grestore}\special{ps: gsave 0 0 0 setrgbcolor}\put(1651,-961){\vector( 3,-1){900}}
\special{ps: grestore}\special{ps: gsave 0 0 0 setrgbcolor}\put(2101,-961){\vector( 3,-1){900}}
\special{ps: grestore}\special{ps: gsave 0 0 0 setrgbcolor}\put(1201,-661){\line(-1, 0){ 75}}
\special{ps: grestore}
\put(5951,-2236){\makebox(0,0)[lb]{\smash{\SetFigFont{11}{13.2}{\rmdefault}{\mddefault}{\updefault}\special{ps: 
gsave 0 0 0 setrgbcolor}$r>3$\special{ps: grestore}}}}
\put(1096,-2516){\makebox(0,0)[rb]{\smash{\SetFigFont{11}{13.2}{\rmdefault}{\mddefault}{\updefault}\special{ps: 
gsave 0 0 0 setrgbcolor}{\small $a=u$}\special{ps: grestore}}}}
\put(1096,-2216){\makebox(0,0)[rb]{\smash{\SetFigFont{11}{13.2}{\rmdefault}{\mddefault}{\updefault}\special{ps: 
gsave 0 0 0 setrgbcolor}{\small $v$}\special{ps: grestore}}}}
\put(1096,-1916){\makebox(0,0)[rb]{\smash{\SetFigFont{11}{13.2}{\rmdefault}{\mddefault}{\updefault}\special{ps: 
gsave 0 0 0 setrgbcolor}{\small $b=uv$}\special{ps: grestore}}}}
\put(1096,-1616){\makebox(0,0)[rb]{\smash{\SetFigFont{11}{13.2}{\rmdefault}{\mddefault}{\updefault}\special{ps: 
gsave 0 0 0 setrgbcolor}{\small $c=v^2$}\special{ps: grestore}}}}
\put(1096,-1316){\makebox(0,0)[rb]{\smash{\SetFigFont{11}{13.2}{\rmdefault}{\mddefault}{\updefault}\special{ps: 
gsave 0 0 0 setrgbcolor}{\small $uv^2$}\special{ps: grestore}}}}
\put(1096,-1016){\makebox(0,0)[rb]{\smash{\SetFigFont{11}{13.2}{\rmdefault}{\mddefault}{\updefault}\special{ps: 
gsave 0 0 0 setrgbcolor}{\small $v^3$}\special{ps: grestore}}}}
\put(1096,-716){\makebox(0,0)[rb]{\smash{\SetFigFont{11}{13.2}{\rmdefault}{\mddefault}{\updefault}\special{ps: 
gsave 0 0 0 setrgbcolor}{\small $uv^3$}\special{ps: grestore}}}}
\put(1606,-2986){\makebox(0,0)[lb]{\smash{\SetFigFont{11}{13.2}{\rmdefault}{\mddefault}{\updefault}\special{ps: 
gsave 0 0 0 setrgbcolor}{\small $x$}\special{ps: grestore}}}}
\put(1606,-3246){\makebox(0,0)[lb]{\smash{\SetFigFont{11}{13.2}{\rmdefault}{\mddefault}{\updefault}\special{ps: 
gsave 0 0 0 setrgbcolor}{\small $y$}\special{ps: grestore}}}}
\put(2056,-2986){\makebox(0,0)[lb]{\smash{\SetFigFont{11}{13.2}{\rmdefault}{\mddefault}{\updefault}\special{ps: 
gsave 0 0 0 setrgbcolor}{\small $z$}\special{ps: grestore}}}}
\put(2026,-3246){\makebox(0,0)[lb]{\smash{\SetFigFont{11}{13.2}{\rmdefault}{\mddefault}{\updefault}\special{ps:
gsave 0 0 0 setrgbcolor}{\small $xy$}\special{ps: grestore}}}}
\put(2046,-3506){\makebox(0,0)[lb]{\smash{\SetFigFont{11}{13.2}{\rmdefault}{\mddefault}{\updefault}\special{ps:
gsave 0 0 0 setrgbcolor}{\small $y^2$}\special{ps: grestore}}}}
\put(2441,-2986){\makebox(0,0)[lb]{\smash{\SetFigFont{11}{13.2}{\rmdefault}{\mddefault}{\updefault}\special{ps: 
gsave 0 0 0 setrgbcolor}{\small $xz$}\special{ps: grestore}}}}
\put(2441,-3246){\makebox(0,0)[lb]{\smash{\SetFigFont{11}{13.2}{\rmdefault}{\mddefault}{\updefault}\special{ps:
gsave 0 0 0 setrgbcolor}{\small $zy$}\special{ps: grestore}}}}
\put(2441,-3506){\makebox(0,0)[lb]{\smash{\SetFigFont{11}{13.2}{\rmdefault}{\mddefault}{\updefault}\special{ps:
gsave 0 0 0 setrgbcolor}{\small $xy^2$}\special{ps: grestore}}}}
\put(2451,-3766){\makebox(0,0)[lb]{\smash{\SetFigFont{11}{13.2}{\rmdefault}{\mddefault}{\updefault}\special{ps:
gsave 0 0 0 setrgbcolor}{\small $y^3$}\special{ps: grestore}}}}
\put(2946,-2986){\makebox(0,0)[lb]{\smash{\SetFigFont{11}{13.2}{\rmdefault}{\mddefault}{\updefault}\special{ps: 
gsave 0 0 0 setrgbcolor}{\small $z^2$}\special{ps: grestore}}}}
\put(2896,-3246){\makebox(0,0)[lb]{\smash{\SetFigFont{11}{13.2}{\rmdefault}{\mddefault}{\updefault}\special{ps:
gsave 0 0 0 setrgbcolor}{\small $xzy$}\special{ps: grestore}}}}
\put(2896,-3506){\makebox(0,0)[lb]{\smash{\SetFigFont{11}{13.2}{\rmdefault}{\mddefault}{\updefault}\special{ps:
gsave 0 0 0 setrgbcolor}{\small $zy^2$}\special{ps: grestore}}}}
\put(2896,-3766){\makebox(0,0)[lb]{\smash{\SetFigFont{11}{13.2}{\rmdefault}{\mddefault}{\updefault}\special{ps:
gsave 0 0 0 setrgbcolor}{\small $xy^3$}\special{ps: grestore}}}}
\put(2946,-4026){\makebox(0,0)[lb]{\smash{\SetFigFont{11}{13.2}{\rmdefault}{\mddefault}{\updefault}\special{ps:
gsave 0 0 0 setrgbcolor}{\small $y^4$}\special{ps: grestore}}}}
\put(1276,-451){\makebox(0,0)[lb]{\smash{\SetFigFont{11}{13.2}{\rmdefault}{\mddefault}{\updefault}\special{ps: 
gsave 0 0 0 setrgbcolor}{\small $H^*({\mathbb Z}_{2^r};{\mathbb F}_2)$}\special{ps: grestore}}}}
\put(3440,-2686){\makebox(0,0)[lb]{\smash{\SetFigFont{11}{13.2}{\rmdefault}{\mddefault}{\updefault}\special{ps: 
gsave 0 0 0 setrgbcolor}{\small $H^*({\mathbb Z}_{2^{r-2}}\oplus{\mathbb Z}_2;{\mathbb F}_2)$}\special{ps: grestore}}}}
\put(2551,-2236){\makebox(0,0)[lb]{\smash{\SetFigFont{11}{13.2}{\rmdefault}{\mddefault}{\updefault}\special{ps: 
gsave 0 0 0 setrgbcolor}{\small $d_2$}\special{ps: grestore}}}}
\end{picture}}

\vspace{.5cm}

\centerline{\setlength{\unitlength}{3400sp}%
\begingroup\makeatletter\ifx\SetFigFont\undefined%
\gdef\SetFigFont#1#2#3#4#5{%
  \reset@font\fontsize{#1}{#2pt}%
  \fontfamily{#3}\fontseries{#4}\fontshape{#5}%
  \selectfont}%
\fi\endgroup%
\begin{picture}(6000,3905)(4000,-4236)
\thinlines
\special{ps: gsave 0 0 0 setrgbcolor}\put(4801,-2761){\vector( 0, 1){2400}}
\special{ps: grestore}\special{ps: gsave 0 0 0 setrgbcolor}\put(4801,-2761){\vector( 1, 0){2400}}
\special{ps: grestore}\special{ps: gsave 0 0 0 setrgbcolor}\put(5251,-2761){\line( 0,-1){ 75}}
\special{ps: grestore}\special{ps: gsave 0 0 0 setrgbcolor}\put(5701,-2761){\line( 0,-1){ 75}}
\special{ps: grestore}\special{ps: gsave 0 0 0 setrgbcolor}\put(6151,-2761){\line( 0,-1){ 75}}
\special{ps: grestore}\special{ps: gsave 0 0 0 setrgbcolor}\put(6601,-2761){\line( 0,-1){ 75}}
\special{ps: grestore}\special{ps: gsave 0 0 0 setrgbcolor}\put(4801,-2461){\line(-1, 0){ 75}}
\special{ps: grestore}\special{ps: gsave 0 0 0 setrgbcolor}\put(4801,-2161){\line(-1, 0){ 75}}
\special{ps: grestore}\special{ps: gsave 0 0 0 setrgbcolor}\put(4801,-1861){\line(-1, 0){ 75}}
\special{ps: grestore}\special{ps: gsave 0 0 0 setrgbcolor}\put(4801,-1561){\line(-1, 0){ 75}}
\special{ps: grestore}\special{ps: gsave 0 0 0 setrgbcolor}\put(4801,-1261){\line(-1, 0){ 75}}
\special{ps: grestore}\special{ps: gsave 0 0 0 setrgbcolor}\put(4801,-961){\line(-1, 0){ 75}}
\special{ps: grestore}\special{ps: gsave 0 0 0 setrgbcolor}\put(4801,-661){\line(-1, 0){ 75}}
\special{ps: grestore}\special{ps: gsave 0 0 0 setrgbcolor}\put(4801,-2161){\vector( 3,-1){900}}
\special{ps: grestore}\special{ps: gsave 0 0 0 setrgbcolor}\put(5251,-2161){\vector( 3,-1){900}}
\special{ps: grestore}\special{ps: gsave 0 0 0 setrgbcolor}\put(5701,-2161){\vector( 3,-1){900}}
\special{ps: grestore}\special{ps: gsave 0 0 0 setrgbcolor}\put(4801,-961){\vector( 3,-1){900}}
\special{ps: grestore}\special{ps: gsave 0 0 0 setrgbcolor}\put(5251,-961){\vector( 3,-1){900}}
\special{ps: grestore}\special{ps: gsave 0 0 0 setrgbcolor}\put(5701,-961){\vector( 3,-1){900}}
\special{ps: grestore}
\put(9551,-2276){\makebox(0,0)[lb]{\smash{\SetFigFont{11}{13.2}{\rmdefault}{\mddefault}{\updefault}\special{ps: 
gsave 0 0 0 setrgbcolor}$r=3$\special{ps: grestore}}}}
\put(4696,-716){\makebox(0,0)[rb]{\smash{\SetFigFont{11}{13.2}{\rmdefault}{\mddefault}{\updefault}\special{ps: 
gsave 0 0 0 setrgbcolor}{\small $uv^3$}\special{ps: grestore}}}}
\put(4696,-1316){\makebox(0,0)[rb]{\smash{\SetFigFont{11}{13.2}{\rmdefault}{\mddefault}{\updefault}\special{ps: 
gsave 0 0 0 setrgbcolor}{\small $uv^2$}\special{ps: grestore}}}}
\put(4696,-1016){\makebox(0,0)[rb]{\smash{\SetFigFont{11}{13.2}{\rmdefault}{\mddefault}{\updefault}\special{ps: 
gsave 0 0 0 setrgbcolor}{\small $v^3$}\special{ps: grestore}}}}
\put(4696,-1616){\makebox(0,0)[rb]{\smash{\SetFigFont{11}{13.2}{\rmdefault}{\mddefault}{\updefault}\special{ps: 
gsave 0 0 0 setrgbcolor}{\small $c=v^2$}\special{ps: grestore}}}}
\put(4696,-1916){\makebox(0,0)[rb]{\smash{\SetFigFont{11}{13.2}{\rmdefault}{\mddefault}{\updefault}\special{ps: 
gsave 0 0 0 setrgbcolor}{\small $b=uv$}\special{ps: grestore}}}}
\put(4696,-2216){\makebox(0,0)[rb]{\smash{\SetFigFont{11}{13.2}{\rmdefault}{\mddefault}{\updefault}\special{ps: 
gsave 0 0 0 setrgbcolor}{\small $v$}\special{ps: grestore}}}}
\put(4696,-2516){\makebox(0,0)[rb]{\smash{\SetFigFont{11}{13.2}{\rmdefault}{\mddefault}{\updefault}\special{ps: 
gsave 0 0 0 setrgbcolor}{\small $a=u$}\special{ps: grestore}}}}
\put(5206,-2986){\makebox(0,0)[lb]{\smash{\SetFigFont{11}{13.2}{\rmdefault}{\mddefault}{\updefault}\special{ps: 
gsave 0 0 0 setrgbcolor}{\small $x$}\special{ps: grestore}}}}
\put(5206,-3246){\makebox(0,0)[lb]{\smash{\SetFigFont{11}{13.2}{\rmdefault}{\mddefault}{\updefault}\special{ps:
gsave 0 0 0 setrgbcolor}{\small $y$}\special{ps: grestore}}}} 
\put(5646,-2986){\makebox(0,0)[lb]{\smash{\SetFigFont{11}{13.2}{\rmdefault}{\mddefault}{\updefault}\special{ps: 
gsave 0 0 0 setrgbcolor}{\small $x^2$}\special{ps: grestore}}}}
\put(5626,-3246){\makebox(0,0)[lb]{\smash{\SetFigFont{11}{13.2}{\rmdefault}{\mddefault}{\updefault}\special{ps:
gsave 0 0 0 setrgbcolor}{\small $xy$}\special{ps: grestore}}}}
\put(5646,-3506){\makebox(0,0)[lb]{\smash{\SetFigFont{11}{13.2}{\rmdefault}{\mddefault}{\updefault}\special{ps:
gsave 0 0 0 setrgbcolor}{\small $y^2$}\special{ps: grestore}}}}
\put(6096,-2986){\makebox(0,0)[lb]{\smash{\SetFigFont{11}{13.2}{\rmdefault}{\mddefault}{\updefault}\special{ps: 
gsave 0 0 0 setrgbcolor}{\small $x^3$}\special{ps: grestore}}}}
\put(6056,-3246){\makebox(0,0)[lb]{\smash{\SetFigFont{11}{13.2}{\rmdefault}{\mddefault}{\updefault}\special{ps:
gsave 0 0 0 setrgbcolor}{\small $x^2y$}\special{ps: grestore}}}}
\put(6056,-3506){\makebox(0,0)[lb]{\smash{\SetFigFont{11}{13.2}{\rmdefault}{\mddefault}{\updefault}\special{ps:
gsave 0 0 0 setrgbcolor}{\small $xy^2$}\special{ps: grestore}}}}
\put(6096,-3766){\makebox(0,0)[lb]{\smash{\SetFigFont{11}{13.2}{\rmdefault}{\mddefault}{\updefault}\special{ps:
gsave 0 0 0 setrgbcolor}{\small $y^3$}\special{ps: grestore}}}}
\put(6546,-2986){\makebox(0,0)[lb]{\smash{\SetFigFont{11}{13.2}{\rmdefault}{\mddefault}{\updefault}\special{ps: 
gsave 0 0 0 setrgbcolor}{\small $x^4$}\special{ps: grestore}}}}
\put(6496,-3246){\makebox(0,0)[lb]{\smash{\SetFigFont{11}{13.2}{\rmdefault}{\mddefault}{\updefault}\special{ps:
gsave 0 0 0 setrgbcolor}{\small $x^3y$}\special{ps: grestore}}}}
\put(6496,-3506){\makebox(0,0)[lb]{\smash{\SetFigFont{11}{13.2}{\rmdefault}{\mddefault}{\updefault}\special{ps:
gsave 0 0 0 setrgbcolor}{\small $x^2y^2$}\special{ps: grestore}}}}
\put(6496,-3766){\makebox(0,0)[lb]{\smash{\SetFigFont{11}{13.2}{\rmdefault}{\mddefault}{\updefault}\special{ps:
gsave 0 0 0 setrgbcolor}{\small $xy^3$}\special{ps: grestore}}}}
\put(6546,-4026){\makebox(0,0)[lb]{\smash{\SetFigFont{11}{13.2}{\rmdefault}{\mddefault}{\updefault}\special{ps:
gsave 0 0 0 setrgbcolor}{\small $y^4$}\special{ps: grestore}}}}
\put(4876,-451){\makebox(0,0)[lb]{\smash{\SetFigFont{11}{13.2}{\rmdefault}{\mddefault}{\updefault}\special{ps: 
gsave 0 0 0 setrgbcolor}{\small $H^*({\mathbb Z}_8;{\mathbb F}_2)$}\special{ps: grestore}}}}
\put(7040,-2686){\makebox(0,0)[lb]{\smash{\SetFigFont{11}{13.2}{\rmdefault}{\mddefault}{\updefault}\special{ps: 
gsave 0 0 0 setrgbcolor}{\small $H^*({\mathbb Z}_2\oplus{\mathbb Z}_2;{\mathbb F}_2)$}\special{ps: grestore}}}}
\put(6151,-2236){\makebox(0,0)[lb]{\smash{\SetFigFont{11}{13.2}{\rmdefault}{\mddefault}{\updefault}\special{ps: 
gsave 0 0 0 setrgbcolor}{\small $d_2$}\special{ps: grestore}}}}
\end{picture}}
\caption{$E_2$ stage of the LHS spectral sequence for $H^*(\mbox{Hol}({\mathbb Z}_{2^r}); {\mathbb F}_2)$, $r \ge 3$}
\label{Lyndon_Hol_Z2r}
\end{figure}

\begin{note} The possible terms are linear combinations of $uy^2$ and $uxy.$
\end{note}

The mod 2 cohomology of $G_2$ (which also can by computed using Wall's method) is given by
$$H^q(G_2; {\mathbb F}_2) \cong \moplus_{q+1} {{\mathbb F}_2}. $$
It follows from the ranks of $H^q(G_2; {\mathbb F}_2)$'s that all the differentials in the LHS spectral sequence for
$H^*(G_2; {\mathbb F}_2)$ are zero.

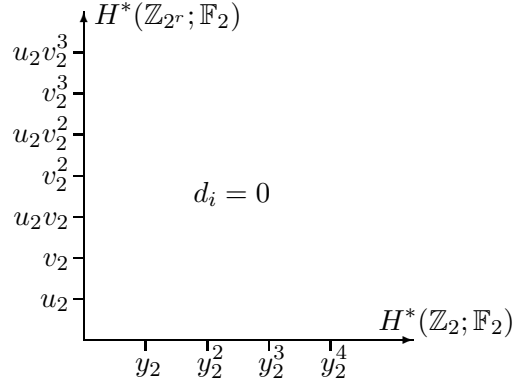
\begin{figure}
\centerline{\setlength{\unitlength}{3400sp}%
\begingroup\makeatletter\ifx\SetFigFont\undefined%
\gdef\SetFigFont#1#2#3#4#5{%
  \reset@font\fontsize{#1}{#2pt}%
  \fontfamily{#3}\fontseries{#4}\fontshape{#5}%
  \selectfont}%
\fi\endgroup%
\begin{picture}(3500,3105)(730,-2986)
\thinlines
\special{ps: gsave 0 0 0 setrgbcolor}\put(1201,-2761){\vector( 1, 0){2400}}
\special{ps: grestore}\special{ps: gsave 0 0 0 setrgbcolor}\put(1201,-2761){\vector( 0, 1){2400}}
\special{ps: grestore}\special{ps: gsave 0 0 0 setrgbcolor}\put(1126,-2461){\line( 1, 0){ 75}}
\special{ps: grestore}\special{ps: gsave 0 0 0 setrgbcolor}\put(1126,-2161){\line( 1, 0){ 75}}
\special{ps: grestore}\special{ps: gsave 0 0 0 setrgbcolor}\put(1126,-1861){\line( 1, 0){ 75}}
\special{ps: grestore}\special{ps: gsave 0 0 0 setrgbcolor}\put(1126,-1561){\line( 1, 0){ 75}}
\special{ps: grestore}\special{ps: gsave 0 0 0 setrgbcolor}\put(1126,-1261){\line( 1, 0){ 75}}
\special{ps: grestore}\special{ps: gsave 0 0 0 setrgbcolor}\put(1126,-961){\line( 1, 0){ 75}}
\special{ps: grestore}\special{ps: gsave 0 0 0 setrgbcolor}\put(1651,-2761){\line( 0,-1){ 75}}
\special{ps: grestore}\special{ps: gsave 0 0 0 setrgbcolor}\put(2101,-2761){\line( 0,-1){ 75}}
\special{ps: grestore}\special{ps: gsave 0 0 0 setrgbcolor}\put(2551,-2761){\line( 0,-1){ 75}}
\special{ps: grestore}\special{ps: gsave 0 0 0 setrgbcolor}\put(3001,-2761){\line( 0,-1){ 75}}
\special{ps: grestore}\special{ps: gsave 0 0 0 setrgbcolor}\put(1201,-661){\line(-1, 0){ 75}}
\put(1096,-2516){\makebox(0,0)[rb]{\smash{\SetFigFont{11}{13.2}{\rmdefault}{\mddefault}{\updefault}\special{ps: 
gsave 0 0 0 setrgbcolor}{\small $u_2$}\special{ps: grestore}}}}
\put(1096,-2216){\makebox(0,0)[rb]{\smash{\SetFigFont{11}{13.2}{\rmdefault}{\mddefault}{\updefault}\special{ps: 
gsave 0 0 0 setrgbcolor}{\small $v_2$}\special{ps: grestore}}}}
\put(1096,-1916){\makebox(0,0)[rb]{\smash{\SetFigFont{11}{13.2}{\rmdefault}{\mddefault}{\updefault}\special{ps: 
gsave 0 0 0 setrgbcolor}{\small $u_2v_2$}\special{ps: grestore}}}}
\put(1096,-1616){\makebox(0,0)[rb]{\smash{\SetFigFont{11}{13.2}{\rmdefault}{\mddefault}{\updefault}\special{ps: 
gsave 0 0 0 setrgbcolor}{\small $v_2^2$}\special{ps: grestore}}}}
\put(1096,-1316){\makebox(0,0)[rb]{\smash{\SetFigFont{11}{13.2}{\rmdefault}{\mddefault}{\updefault}\special{ps: 
gsave 0 0 0 setrgbcolor}{\small $u_2v_2^2$}\special{ps: grestore}}}}
\put(1096,-1016){\makebox(0,0)[rb]{\smash{\SetFigFont{11}{13.2}{\rmdefault}{\mddefault}{\updefault}\special{ps: 
gsave 0 0 0 setrgbcolor}{\small $v_2^3$}\special{ps: grestore}}}}
\put(1096,-716){\makebox(0,0)[rb]{\smash{\SetFigFont{11}{13.2}{\rmdefault}{\mddefault}{\updefault}\special{ps: 
gsave 0 0 0 setrgbcolor}{\small $u_2v_2^3$}\special{ps: grestore}}}}
\put(1576,-2986){\makebox(0,0)[lb]{\smash{\SetFigFont{11}{13.2}{\rmdefault}{\mddefault}{\updefault}\special{ps: 
gsave 0 0 0 setrgbcolor}{\small $y_2$}\special{ps: grestore}}}}
\put(2026,-2986){\makebox(0,0)[lb]{\smash{\SetFigFont{11}{13.2}{\rmdefault}{\mddefault}{\updefault}\special{ps: 
gsave 0 0 0 setrgbcolor}{\small $y_2^2$}\special{ps: grestore}}}}
\put(2471,-2986){\makebox(0,0)[lb]{\smash{\SetFigFont{11}{13.2}{\rmdefault}{\mddefault}{\updefault}\special{ps: 
gsave 0 0 0 setrgbcolor}{\small $y_2^3$}\special{ps: grestore}}}}
\put(2926,-2986){\makebox(0,0)[lb]{\smash{\SetFigFont{11}{13.2}{\rmdefault}{\mddefault}{\updefault}\special{ps: 
gsave 0 0 0 setrgbcolor}{\small $y_2^4$}\special{ps: grestore}}}}
\put(1276,-451){\makebox(0,0)[lb]{\smash{\SetFigFont{11}{13.2}{\rmdefault}{\mddefault}{\updefault}\special{ps: 
gsave 0 0 0 setrgbcolor}{\small $H^*({\mathbb Z}_{2^r};{\mathbb F}_2)$}\special{ps: grestore}}}}
\put(3350,-2686){\makebox(0,0)[lb]{\smash{\SetFigFont{11}{13.2}{\rmdefault}{\mddefault}{\updefault}\special{ps: 
gsave 0 0 0 setrgbcolor}{\small $H^*({\mathbb Z}_{2};{\mathbb F}_2)$}\special{ps: grestore}}}}
\put(2000,-1750){\makebox(0,0)[lb]{\smash{\SetFigFont{11}{13.2}{\rmdefault}{\mddefault}{\updefault}\special{ps:
gsave 0 0 0 setrgbcolor}{\small $d_i=0$}\special{ps: grestore}}}}
\end{picture}}
\caption{$E_2$ stage of the LHS spectral sequence for $H^*(G_2; {\mathbb F}_2)$}
\end{figure}

The spectral sequence for $H^*(\mbox{Hol}({\mathbb Z}_{2^r}); {\mathbb F}_2)$ maps into it; $x$ and $z$ are mapped to
0, $y$ to $y_2,$ $u$ to $u_2,$ and $v$ to $v_2.$ Therefore the term $uy^2$ is not present in the expression for $d_2(v)$.

Finally, considering the map from the spectral sequence for $H^*(\Hol(\Z_{2^r}); \F_2)$ into the spectral 
sequence for $H^*(G_3; \F_2)$, we conclude that the term $uxy$ is not present in the expression for $d_2(v)$ either. 
Therefore
$$d_2(v)=uz \mbox{ \ if } r>3,$$
$$d_2(v)=ux^2 \mbox{ \ if } r=3.$$
Then
$$d_2(v^{2n+1} x^i y^j z^k)= uv^{2n} x^i y^j z^{k+1} \mbox{ \ if } r>3,$$
$$d_2(v^{2n+1} x^i y^j)= uv^{2n} x^{i+2} y^j \mbox{ \ if } r=3,$$
and all other differentials are zero.

Let $a=u$, $b=uv$, and $c=v^2$. Then, additively, the $E_\infty$ stage in the spectral sequence for $\Hol(\Z_{2^r})$ is

$$E^{*,*}_\infty \cong \Lambda (a, b, x) \otimes
{\mathbb F_2}[c, \ z, \ y] / ab=az=0 \mbox{ \ \ if } r>3,$$
and
$$E^{*,*}_\infty \cong \Lambda (a, b) \otimes {\mathbb F}_2[c, \ x, \ y] /
ab=ax^2=0 \mbox{ \ \ if } r=3$$
where $|a|=|x|=|y|=1, \ |z|=2, \ |b|=3, \ |c|=4.$

However, the ring structure of the mod 2 cohomology is different.

First we compute the mod 2 cohomology rings of subgroups $G_1,$ $G_2,$ and $G_3$. This was done by the author using Bockstein 
homomorphisms and Steenrod square operations, however, it was learned later that these cohomology rings were worked out in \cite{H}. 
Thus we skip the details of our calculations. 

$G_1$ is a special case of theorem C, part (2), in \cite{H}: 
\begin{itemize}
\item If $r>3$, then $H^*(G_1;\F_2) \cong $ \\ 
$\cong \Lambda(x_1) \otimes \F_2[a_1, b_1, c_1, z_1] / a_1^2=a_1x_1, a_1z_1=0, a_1b_1=0, b_1^2=a_1x_1c_1 + b_1x_1z_1$ \\
where $|x_1| = |a_1| = 1,$ $|z_1|=2,$ $b_1|=3,$ and $|c_1|=4$.   
\item If $r=3$, then $H^*(G_1;\F_2) \cong $ \\
$ \cong \F_2[a_1, b_1, c_1, x_1] / a_1^2=a_1x_1, a_1x_1^2=0, a_1b_1=0, b_1^2=a_1x_1c_1 + b_1x_1^3 + c_1x_1^2$ \\
where $|x_1| = |a_1| = 1,$ $b_1|=3,$ and $|c_1|=4$.
\end{itemize}

$G_2$ is covered by theorem B, part (2): 
$$H^*(G_2; \F_2) \cong \F_2[u_2,v_2,y_2] / u_2 = u_2y_2,$$
where $|u_2|=|y_2|=1,$ and $|v_2|=2$. 

Finally, $G_3$ is covered also by theorem C, part (2). 

Next we compute the ring structure of the mod 2 cohomology of $\Hol(\Z_{2^r})$. First notice that the split sequence 
\eqref{split-short-seq-hol2r} gives the maps of cohomology rings
$$ H^*(\Z_{2^r}) \longleftarrow H^*(\Hol(\Z_{2^r});\F_2) \buildrel \curvearrowright \over \longleftarrow 
H^*(\Z_{2^{r-2}}\oplus\Z_2;\F_2).$$ 
Therefore all the products of $x$, $y$, and $z$ are as in $H^*(\Z_{2^{r-2}}\oplus\Z_2;\F_2)$, namely the only relation is that $x^2=0$ 
in case $r>3$. The products $ac$, $bc$, and $c^2$ are non-zero in $H^*(\Z_{2^r};\F_2)$, and hence are not linear combinations of  other 
products in $H^*(\Hol(\Z_{2^r});\F_2)$. The products $a^2$, $ab$, and $b^2$ all map to 0 in $H^*(\Z_{2^r};\F_2)$, thus these are 
linear combinations of other products each of which necessarily contains $x$, $y$, or $z$. On the other hand, there are choices of $a,$ 
$b,$ and $c$ that are mapped to 0 in $H^*(\Z_{2^{r-2}}\oplus\Z_2;\F_2)$, thus the expressions for $a^2$, $ab$, and $b^2$ in 
$H^*(\Hol(\Z_{2^r});\F_2)$ do not contain products of only $x$, $y$, and $z$. 
Thus $a^2$ is a linear combination of $ax$ and $ay$; The product $ab$ is a linear combination of $bx$, $by$, $axy^2$, $ay^2$, and 
either $z^2$ or $x^4$ depending on whether $p>3$ or $p=3$. Finally, $b^2$ is a linear combination of $acx$, $acy$, $cz/cx^2$, $cxy$, 
$cy^2$, $bxz/bx^3$, $byz/bx^2y$, $bxy^2$, $by^3$, $axy^4$, and $ay^4$.  

Now consider maps from $H^*(\Hol(\Z_{2^r}))$ into $H^*(G_i)$ for $1\le i\le 3$. 

As mentioned before, the map $H^*(\Hol(\Z_{2^r})) \rightarrow H^*(G_1)$ is given by $x\mapsto x_1$, $y\mapsto 0$, $z\mapsto z_1$, 
$a\mapsto a_1$, $b\mapsto b_1$, and $c\mapsto c_1$. Since $a_1^2=a_1x_1$, the term $ax$ is present in $a^2$. Since $a_1b_1=0$, neither 
of the terms not containing $y$ is present in $ab$. The formulas for $b_1^2$ above imply that if $r>3$ then $acx$ and $bxz$ are present 
in $b^2$ but $cz$ is not, while if $r=3$ then $acx$, $cx^2$, and $bx^3$ are all present. 

The map $H^*(\Hol(\Z_{2^r})) \rightarrow H^*(G_2)$ is given by $x\mapsto 0$, $y\mapsto y_2$, $z\mapsto 0$, $a\mapsto u_2$, 
$b\mapsto u_2v_2$, and $c\mapsto v_2^2$. The equation $u_2^2=u_2y_2$ implies that the term $ay$ is also present in $a^2$, so  
$a^2=ax+ay$. The product 
$ab$ is mapped into $u_2^2v_2 = u_2v_2y_2$ which is the image of $by$, therefore $by$ is present in $ab$, but $ay^3$ is not present. 
Further, $b^2$ is mapped into $(u_2v_2)^2 = u_2^2v_2^2 = u_2v_2^2y_2$, the image of $acy$, therefore $acy$ is present in $b^2$ but 
$cy^2$, $by^3$, and $ay^4$, are not present. 

Finally, considering the map $H^*(\Hol(\Z_{2^r})) \rightarrow H^*(G_3)$ we see that there are choices of $b$ and $c$ such that neither 
of the other terms is present in $ab$, and $b^2$, therefore the ring structure of the cohomology of $\Hol(\Z_{2^r})$ is as in the 
following theorem. 

\begin{theorem}\label{coh-ring-p2} There are choices of generators $a,$ $b,$ $c,$ $x,$ $y,$ and $z$ such that 
\begin{enumerate} 
\item if $r>3$, then  
$H^*(\Hol(\Z_{2^r};\F_2) \cong $\\
$\cong \Lambda(x) \otimes \F_2[a,b,c,y,z] / a^2=ax+ay, az=0, ab=by, b^2=acx+bxz+acy,$ \\
where $|x|=|y|=|a|=1,$ $|z|=2$, $|b|=3$, and $|c|=4$, 
\item  $H^*(\Hol(\Z_8;\F_2) \cong $\\
$\cong \F_2[a,b,c,x, y] / a^2=ax+ay, ax^2=0, ab=by, b^2=acx+bx^3+acy+cx^2,$ \\
where $|x|=|y|=|a|=1,$ $|b|=3$, and $|c|=4$.
\end{enumerate}
\end{theorem}


\section{Description and integer homology of $\Hol(\Z_{p^r})$ for $p$ odd}

\label{description-HolZpr}

Let $p$ be an odd prime. It is well-known that Aut$({\mathbb Z}_{p^r}) \cong {\mathbb Z}_{(p-1)p^{r-1}}$
(e.g., see \cite[p.98]{R}).
Let multiplication by $s$ be a generator of Aut$({\mathbb Z}_{p^r})$. Then the split short exact sequence
$$1 \rightarrow \Z_{p^r} \rightarrow \Hol(\Z_{p^r}) \buildrel \curvearrowleft \over \rightarrow \mbox{Aut}({\mathbb Z}_{p^r}) \rightarrow 1 $$
implies that
$$\mbox{Hol}({\mathbb Z}_{p^r}) \cong \{x, \ y | \ x^{(p-1)p^{r-1}}=y^{p^r}=1, \ yx=xy^s \}$$

Let A be Wall's free resolution for $G=\Hol(\Z_{p^r}):$ $A_{n,m}$ is free on one generator $a_{n,m}$ $(n\ge 0, m\ge 0)$ and 
the 
differentials as given in \cite{W}.
\begin{note}
We switched the indices so that the bidegree of the differential $d_k$ in the spectral sequence below is $(-k,k-1)$.
\end{note}
We take the tensor product with ${\mathbb Z}$ over ${\mathbb Z}(G).$ Let $g_{n,m}$ be the generator of
$A_{n,m}\otimes_{{\mathbb Z}(G)}{\mathbb Z} \cong {\mathbb Z}$. The differentials are given by

\begin{equation} \label{differentials-p} \begin{array}{@{}rcl@{}}
d_0 g_{n,2m+1} &=& 0, \\
d_0 g_{n,2m} &=& p^r g_{n,2m-1}, \\
d_1 g_{2n+1,2m} &=& (s^m-1) g_{2n,2m}, \\
d_1 g_{2n+1,2m-1} &=& -(s^m-1) g_{2n,2m-1}, \\
d_1 g_{2n,2m} &=& \displaystyle \left( \sum_{j=0}^{(p-1)p^{r-1}-1}s^{mj} \right) g_{2n-1,2m}, \\
\end{array} \end{equation}

$$\begin{array}{@{}rcl@{}}
d_1 g_{2n,2m-1} &=& \displaystyle - \left( \sum_{j=0}^{(p-1)p^{r-1}-1}s^{mj} \right) g_{2n-1,2m-1}, \\
d_2 g_{n,2m} &=& 0, \\
d_2 g_{n,2m-1} &=& -\frac{1}{p^r} \left( s^{ m(p-1)p^{r-1} } -1 \right) g_{n-2,2m},
\end{array} $$ 
and $d_k=0$ for $k>2$.

As Wall pointed out, the submodules
$$\begin{array}{rcl}
A^0 & = & \displaystyle \sum_{n\ge 0}{A_{n,0} \otimes_G {\mathbb Z}}, \\
A^m & = & \displaystyle \sum_{n\ge 0}{(A_{2m-1,n}+A_{2m,n}) \otimes_G {\mathbb Z}}
\end{array}$$
are invariant under $d$, and $A \otimes_G \Z$ is their direct sum. Hence it is sufficient to compute the homology of $d$ on
them separately. The homology of $A^0$ is that of ${\mathbb Z}_{(p-1)p^{r-1}}$. For $A^m$, $m>0$, consider the spectral
sequence of the filtration
$\displaystyle F_k(A^m) = \sum_{s\le k}A_{s,*}.$ $d_0$ is induced by $N_y$, thus is multiplication by
$p^r.$ Hence the $E_1$ stage collapses to the \mbox{$(2m-1)$-row.}

By Lemma \ref{obvious}, it suffices to know the highest powers of $p$ that divide the coefficients in
\eqref{differentials-p} in order to compute $H_* (A^m)$.

\begin{lemma} \label{surj_induces_surj}  
If $1 \le q \le r-1$, the natural surjection $\Z_{p^r} \buildrel j \over \rightarrow \Z_{p^q}$ induces a surjection of the
automorphism groups $\Aut (\Z_{p^r}) \buildrel j' \over \rightarrow \Aut (\Z_{p^q})$.
\end{lemma}

\begin{proof}
See example \ref{characteristic-example-Z}. 
\end{proof}

\begin{lemma} \label{divisible-p}
If $1 \le q \le r-1$ and $(k,p)=1,$ then $\nu_p(s^{(p-1)p^{q-1}k}-1)=q$,  
and $s^m-1$ is not divisible by $p$ if $(m,p-1)=1.$
\end{lemma}

\begin{proof}
Lemma \ref{surj_induces_surj} implies that multiplication by $s$ generates $\Aut(\Z_{p^q})$ for
$1 \le q \le r-1$. Thus $s^{(p-1)p^{q-1}}$ is the smallest power of $s$ that is congruent to 1 mod $p^q$.
\end{proof}  

Therefore 
\begin{itemize}
\item
if $(m,p-1)=1,$ then
$\nu_p(s^m-1)=0$ and $\displaystyle \nu_p \left( \sum_{j=0}^{(p-1)p^{r-1}-1}s^{mj} \right) = r$.
\item
if $m=(p-1)p^{q-1}k$ where $(k,p)=1$, then $\nu_p(s^m -1)=q$ and \\
$\displaystyle \nu_p \left( \sum_{j=0}^{(p-1)p^{r-1}-1}s^{mj} \right) = \nu_p \left( \frac{s^{m(p-1)p^{r-1}}-1}{s^m-1} 
\right) = r-1$;
\item
if $m=(p-1)p^{r-1}k$, then $s^m \equiv 1(\mbox{mod } p^r)$ and
$\displaystyle \nu_p \left( \sum_{j=0}^{(p-1)p^{r-1}-1}s^{mj}\right) = r-1$.
\end{itemize}

It follows that if $r=1$ then

$$H_i(A^m;{\mathbb Z}) = \left\{ \begin{array}{ll}
{\mathbb Z}_p & \mbox{if } m=(p-1)k \mbox{ and } i=2m-1 \\
0 & \mbox{otherwise}
\end{array} \right.$$

and if $r\ge 2$ then
$$H_i(A^m;{\mathbb Z})= \left\{ \begin{array}{ll}
{\mathbb Z}_{p^q} & \mbox{if } m=(p-1)p^{q-1}k, \  1 \le q \le r-1, (k,p)=1, \ i=2m-1 \\
{\mathbb Z}_{p^{q-1}} & \mbox{if } m=(p-1)p^{q-1}k, \  2 \le q \le r-1, (k,p)=1, \ i \ge 2m \\
{\mathbb Z}_{p^r} & \mbox{if } m=(p-1)p^{r-1}k  \mbox{ and } i=2m-1 \\
{\mathbb Z}_{p^{r-1}} & \mbox{if } m=(p-1)p^{r-1}k \mbox{ and } i \ge 2m \\
0 & \mbox{otherwise }
\end{array} \right.$$

Now we just add up the homologies of $A^m$, and we get:

$$H_q(\mbox{Hol}({\mathbb Z}_p);{\mathbb Z}) \cong \left(\moplus_{n_1} {\mathbb Z}_p \right) \moplus
\left(\moplus_m {\mathbb Z}_{p-1}\right)$$

where 
\bigskip

$ \begin{array}{lll}
n_1 &=& \left\{ \begin{array}{ll}
                1 & \mbox{ if } q=2(p-1)k-1 \\
                0 & \mbox{ otherwise}
                \end{array} \right. \\
 & & \\
m &=& \left\{ \begin{array}{ll}
                1 & \mbox{ if } q \mbox{ is odd } \\
                0 & \mbox{ otherwise}
                \end{array} \right.
\end{array}
$

$$H_q(\mbox{Hol}({\mathbb Z}_{p^2});{\mathbb Z}) \cong \moplus_{i=1}^2\left(\moplus_{n_i}{\mathbb
Z}_{p^i}\right)\moplus\left(\moplus_m{\mathbb Z}_{p-1}\right)$$

where 

\bigskip

$ \begin{array}{@{}lll@{}}
n_1 &=& \displaystyle \left[ \frac{q}{2(p-1)p} \right] +
                \left\{ \begin{array}{ll}
                1 & \mbox{if } q=2(p-1)k-1, \  (k,p)=1 \\
                0 & \mbox{otherwise}
                \end{array} \right. +
                        \left\{ \begin{array}{ll}
                        1 & \mbox{if } q \mbox{ is odd } \\
                        0 & \mbox{otherwise}
                        \end{array} \right. \\
\end{array}$

$ \begin{array}{@{}lll@{}} 
n_2 &=& \left\{ \begin{array}{ll}
                1 & \mbox{ if } q=2(p-1)pk-1 \\
                0 & \mbox{ otherwise} 
                \end{array} \right. \\
 & & \\ 
m &=& \left\{ \begin{array}{ll}
                1 & \mbox{ if } q \mbox{ is odd } \\
                0 & \mbox{ otherwise}
                \end{array} \right.
\end{array}
$

\bigskip

If $r\ge 3$,
$$H_q(\mbox{Hol}({\mathbb Z}_{p^r});{\mathbb Z}) \cong \moplus_{i=1}^r\left(\moplus_{n_i}{\mathbb
Z}_{p^i}\right)\moplus\left(\moplus_m{\mathbb Z}_{p-1}\right)$$

where 

\bigskip

$ \begin{array}{@{}lll@{}}
n_i &=& \displaystyle \left[ \frac{q}{2(p-1)p^i} \right] - \left[ \frac{q}{2(p-1)p^{i+1}} \right] + \\
 & & \\
 & & +  \left\{ \begin{array}{ll}
        1 & \mbox{if } q=2(p-1)p^{i-1}k-1, \  (k,p)=1 \\
	0 & \mbox{otherwise}
        \end{array}  \right. \hfill (1 \le i \le r-2) \\
 & & \\ 
n_{r-1} &=& \displaystyle \left[ \frac{q}{2(p-1)p^{r-1}} \right] +
                \left\{ \begin{array}{ll}
                1 & \mbox{if } q=2(p-1)p^{r-2}k-1, \ (k,p)=1 \\
		0 & \mbox{otherwise}  
                \end{array}  \right. + \\
 & & \\
 & & +                  \left\{ \begin{array}{ll}   
                        1 & \mbox{if } q \mbox{ is odd } \\
                        0 & \mbox{otherwise}
                        \end{array} \right. \\
 & & \\
n_r &=& \left\{ \begin{array}{ll}
        1 & \mbox{if } q=2(p-1)p^{r-1}k-1 \\
	0 & \mbox{otherwise}
        \end{array}  \right. \\
 & & \\ 
m &=& \left\{ \begin{array}{ll}
        1 & \mbox{if } q \mbox{ is odd } \\
        0 & \mbox{otherwise}
        \end{array} \right.
\end{array}
$
\bigskip

In the future sections we will need integer cohomology rather than homology, but it is easy to compute using the universal
coefficients theorem.
 
\section{$H^*(\Hol(\Z_{p^r});\F_p)$ for $p$ odd prime}
                
The mod $p$ cohomology of $\Hol(\Z_p)$ and of $\Hol(\Z_{p^2})$ are well-known (let us just mention that the
LHS spectral sequence collapses for $r=1$ but does not collapse for $r=2$), so we only give the result for $r\ge 3$. 
First we use the universal coefficients theorem to compute $H^*(\mbox{Hol}({\mathbb Z}_{p^r};{\mathbb F}_p)$ additively:

$$H^q(\Hol(\Z_{p^r});\F_p) \cong \moplus_n \F_p $$
where

$$ n = \left\{ \begin{array}{lll}   
2k+2 & \mbox{if} & q=2(p-1)pk+2(p-1)l-1, \ 1 \le l \le p \\
     &           & \mbox{or } q=2(p-1)pk+2(p-1)l, \ 1 \le l \le p-1 \\
2k+1 & \mbox{if} & q=2(p-1)pk+m, \ 0 \le m \le 2(p-1)p-2, \\
     &           & \mbox{but } q \mbox{ is not of either of the above forms}
\end{array} \right.
$$
  
To find the ring structure of $H^*(\mbox{Hol}({\mathbb Z}_{p^r}); {\mathbb F}_p)$ we use the LHS
spectral sequence. First consider the top row in the diagram

\begin{equation} \label{p-diagram} \begin{CD}
  @.   1        @.   1              @. 1                   @.      \\
@.     @VVV          @VVV                @VVV                   @. \\
1 @>>> \Z_{p^r} @>>> G              @>>> \Z_{p^{r-1}}      @>>> 1  \\
@.     @|            @VVV                @VVV                   @. \\  
1 @>>> \Z_{p^r} @>>> \Hol(\Z_{p^r}) @>>> \Z_{(p-1)p^{r-1}} @>>> 1  \\
@.     @VVV          @VVV                @VVV                   @. \\
1 @>>> 1        @>>> Z_{p-1}        @=   \Z_{p-1}          @>>> 1  \\
@.     @VVV          @VVV                @VVV                   @. \\              
  @.   1        @.   1              @.   1                 @.
\end{CD} \end{equation}

(Let us choose $g_1=s^{p-1}$ as a generator for ${\mathbb Z}_{p^{r-1}}$, and $g_2=s^{p^{r-1}}$ as a generator for ${\mathbb
Z}_{p-1}.$)

$G$ is a split extension of cyclic by cyclic groups, so we can use Wall's method to compute
its integer homology, and then use the universal coefficients theorem to compute its mod $p$ cohomology:

$$H^q(G; {\mathbb F}_p) \cong \moplus_n {\mathbb F}_p$$
where
$$n=\left\{ \begin{array}{rl}
2k+1 & \mbox{if } q=2pk \\   
2k+2 & \mbox{if } 2pk+1 \le q \le 2p(k+1) -1
\end{array} \right.$$

Now consider the LHS spectral sequence for $H^*(G;{\mathbb F}_p).$  By lemma \ref{divisible-p}, \\
\mbox{$g_1=s^{p-1}\equiv 1 (\mbox{mod }p),$} therefore
$E_2^{p,q} \cong H^p({\mathbb Z}_{p^{r-1}}; H^q({\mathbb Z}_{p^r}; {\mathbb F}_p))$ has trivial local
coefficients, so
$E_2^{p,q} \cong H^p({\mathbb Z}_{p^{r-1}}; {\mathbb F}_p) \otimes  H^q({\mathbb Z}_{p^r}; {\mathbb F}_p)$.
  
\begin{figure}
\centerline{\setlength{\unitlength}{3400sp}%
\begingroup\makeatletter\ifx\SetFigFont\undefined%
\gdef\SetFigFont#1#2#3#4#5{%
  \reset@font\fontsize{#1}{#2pt}%
  \fontfamily{#3}\fontseries{#4}\fontshape{#5}%
  \selectfont}%
\fi\endgroup%
\begin{picture}(4350,5005)(40,-3236)
\thinlines
\special{ps: gsave 0 0 0 setrgbcolor}\put(1201,-2761){\vector( 1, 0){2400}}
\special{ps: grestore}\special{ps: gsave 0 0 0 setrgbcolor}\put(1201,-2761){\vector( 0, 1){4000}}
\special{ps: grestore}\special{ps: gsave 0 0 0 setrgbcolor}\put(1126,-2461){\line( 1, 0){ 75}}
\special{ps: grestore}\special{ps: gsave 0 0 0 setrgbcolor}\put(1126,-2161){\line( 1, 0){ 75}}
\special{ps: grestore}\special{ps: gsave 0 0 0 setrgbcolor}\put(1126,-1861){\line( 1, 0){ 75}}
\special{ps: grestore}\special{ps: gsave 0 0 0 setrgbcolor}\put(1126,-1561){\line( 1, 0){ 75}}
\special{ps: grestore}\special{ps: gsave 0 0 0 setrgbcolor}\put(1126,-661){\line( 1, 0){ 75}}
\special{ps: grestore}\special{ps: gsave 0 0 0 setrgbcolor}\put(1126,-361){\line( 1, 0){ 75}}
\special{ps: grestore}\special{ps: gsave 0 0 0 setrgbcolor}\put(1126,-61){\line( 1, 0){ 75}}
\special{ps: grestore}\special{ps: gsave 0 0 0 setrgbcolor}\put(1126,239){\line( 1, 0){ 75}}
\special{ps: grestore}\special{ps: gsave 0 0 0 setrgbcolor}\put(1126,539){\line( 1, 0){ 75}}
\special{ps: grestore}\special{ps: gsave 0 0 0 setrgbcolor}\put(1126,839){\line( 1, 0){ 75}}
\special{ps: grestore}\special{ps: gsave 0 0 0 setrgbcolor}\put(1651,-2761){\line( 0,-1){ 75}}
\special{ps: grestore}\special{ps: gsave 0 0 0 setrgbcolor}\put(2101,-2761){\line( 0,-1){ 75}}
\special{ps: grestore}\special{ps: gsave 0 0 0 setrgbcolor}\put(2551,-2761){\line( 0,-1){ 75}}
\special{ps: grestore}\special{ps: gsave 0 0 0 setrgbcolor}\put(3001,-2761){\line( 0,-1){ 75}}
\special{ps: grestore}\special{ps: gsave 0 0 0 setrgbcolor}\put(1201,-2161){\vector( 3,-1){900}}
\special{ps: grestore}\special{ps: gsave 0 0 0 setrgbcolor}\put(1651,-2161){\vector( 3,-1){900}}
\special{ps: grestore}\special{ps: gsave 0 0 0 setrgbcolor}\put(2101,-2161){\vector( 3,-1){900}}
\special{ps: grestore}\special{ps: gsave 0 0 0 setrgbcolor}\put(1201,-1561){\vector( 3,-1){900}}
\special{ps: grestore}\special{ps: gsave 0 0 0 setrgbcolor}\put(1651,-1561){\vector( 3,-1){900}}
\special{ps: grestore}\special{ps: gsave 0 0 0 setrgbcolor}\put(2101,-1561){\vector( 3,-1){900}}
\special{ps: grestore}\special{ps: gsave 0 0 0 setrgbcolor}\put(1201,-361){\vector( 3,-1){900}}
\special{ps: grestore}\special{ps: gsave 0 0 0 setrgbcolor}\put(1651,-361){\vector( 3,-1){900}}
\special{ps: grestore}\special{ps: gsave 0 0 0 setrgbcolor}\put(2101,-361){\vector( 3,-1){900}}
\special{ps: grestore}\special{ps: gsave 0 0 0 setrgbcolor}\put(1201,839){\vector( 3,-1){900}}
\special{ps: grestore}\special{ps: gsave 0 0 0 setrgbcolor}\put(1651,839){\vector( 3,-1){900}}
\special{ps: grestore}\special{ps: gsave 0 0 0 setrgbcolor}\put(2101,839){\vector( 3,-1){900}}
\special{ps: grestore}\special{ps: gsave 0 0 0 setrgbcolor}\put(1201,-661){\line(-1, 0){ 75}}
\put(1096,-2516){\makebox(0,0)[rb]{\smash{\SetFigFont{11}{13.2}{\rmdefault}{\mddefault}{\updefault}\special{ps: 
gsave 0 0 0 setrgbcolor}{\small $a_0=u$}\special{ps: grestore}}}}
\put(1096,-2216){\makebox(0,0)[rb]{\smash{\SetFigFont{11}{13.2}{\rmdefault}{\mddefault}{\updefault}\special{ps: 
gsave 0 0 0 setrgbcolor}{\small $v$}\special{ps: grestore}}}}
\put(1096,-1916){\makebox(0,0)[rb]{\smash{\SetFigFont{11}{13.2}{\rmdefault}{\mddefault}{\updefault}\special{ps: 
gsave 0 0 0 setrgbcolor}{\small $a_1=uv$}\special{ps: grestore}}}}
\put(1096,-1616){\makebox(0,0)[rb]{\smash{\SetFigFont{11}{13.2}{\rmdefault}{\mddefault}{\updefault}\special{ps: 
gsave 0 0 0 setrgbcolor}{\small $v^2$}\special{ps: grestore}}}}
\put(1096,-716){\makebox(0,0)[rb]{\smash{\SetFigFont{11}{13.2}{\rmdefault}{\mddefault}{\updefault}\special{ps: 
gsave 0 0 0 setrgbcolor}{\small $a_{p-2}=uv^{p-2}$}\special{ps: grestore}}}}
\put(1096,-416){\makebox(0,0)[rb]{\smash{\SetFigFont{11}{13.2}{\rmdefault}{\mddefault}{\updefault}\special{ps:
gsave 0 0 0 setrgbcolor}{\small $v^{p-1}$}\special{ps: grestore}}}}
\put(1096,-116){\makebox(0,0)[rb]{\smash{\SetFigFont{11}{13.2}{\rmdefault}{\mddefault}{\updefault}\special{ps:
gsave 0 0 0 setrgbcolor}{\small $b=uv^{p-1}$}\special{ps: grestore}}}}
\put(1096,184){\makebox(0,0)[rb]{\smash{\SetFigFont{11}{13.2}{\rmdefault}{\mddefault}{\updefault}\special{ps:
gsave 0 0 0 setrgbcolor}{\small $c=v^p$}\special{ps: grestore}}}}
\put(1096,484){\makebox(0,0)[rb]{\smash{\SetFigFont{11}{13.2}{\rmdefault}{\mddefault}{\updefault}\special{ps:
gsave 0 0 0 setrgbcolor}{\small $uv^p$}\special{ps: grestore}}}}
\put(1096,784){\makebox(0,0)[rb]{\smash{\SetFigFont{11}{13.2}{\rmdefault}{\mddefault}{\updefault}\special{ps:
gsave 0 0 0 setrgbcolor}{\small $v^{p+1}$}\special{ps: grestore}}}}
\put(1576,-2986){\makebox(0,0)[lb]{\smash{\SetFigFont{11}{13.2}{\rmdefault}{\mddefault}{\updefault}\special{ps: 
gsave 0 0 0 setrgbcolor}{\small $x$}\special{ps: grestore}}}}
\put(2026,-2986){\makebox(0,0)[lb]{\smash{\SetFigFont{11}{13.2}{\rmdefault}{\mddefault}{\updefault}\special{ps: 
gsave 0 0 0 setrgbcolor}{\small $z$}\special{ps: grestore}}}}
\put(2401,-2986){\makebox(0,0)[lb]{\smash{\SetFigFont{11}{13.2}{\rmdefault}{\mddefault}{\updefault}\special{ps: 
gsave 0 0 0 setrgbcolor}{\small $xz$}\special{ps: grestore}}}}
\put(2926,-2986){\makebox(0,0)[lb]{\smash{\SetFigFont{11}{13.2}{\rmdefault}{\mddefault}{\updefault}\special{ps: 
gsave 0 0 0 setrgbcolor}{\small $z^2$}\special{ps: grestore}}}}
\put(1276,1249){\makebox(0,0)[lb]{\smash{\SetFigFont{11}{13.2}{\rmdefault}{\mddefault}{\updefault}\special{ps: 
gsave 0 0 0 setrgbcolor}{\small $H^*({\mathbb Z}_{p^r};{\mathbb F}_p)$}\special{ps: grestore}}}}
\put(3390,-2686){\makebox(0,0)[lb]{\smash{\SetFigFont{11}{13.2}{\rmdefault}{\mddefault}{\updefault}\special{ps: 
gsave 0 0 0 setrgbcolor}{\small $H^*({\mathbb Z}_{p^{r-1}};{\mathbb F}_p)$}\special{ps: grestore}}}}
\put(2551,-2236){\makebox(0,0)[lb]{\smash{\SetFigFont{11}{13.2}{\rmdefault}{\mddefault}{\updefault}\special{ps: 
gsave 0 0 0 setrgbcolor}{\small $d_2$}\special{ps: grestore}}}}
\put(1700,-1100){\makebox(0,0)[lb]{\smash{\SetFigFont{11}{13.2}{\rmdefault}{\mddefault}{\updefault}\special{ps:
gsave 0 0 0 setrgbcolor}{\small \ldots}\special{ps: grestore}}}}
\end{picture}}
\caption{$E_2$ stage of the LHS spectral sequence for $H^*(G;{\mathbb F}_p)$}
\label{Lyndon_Hol_Zpr}
\end{figure}
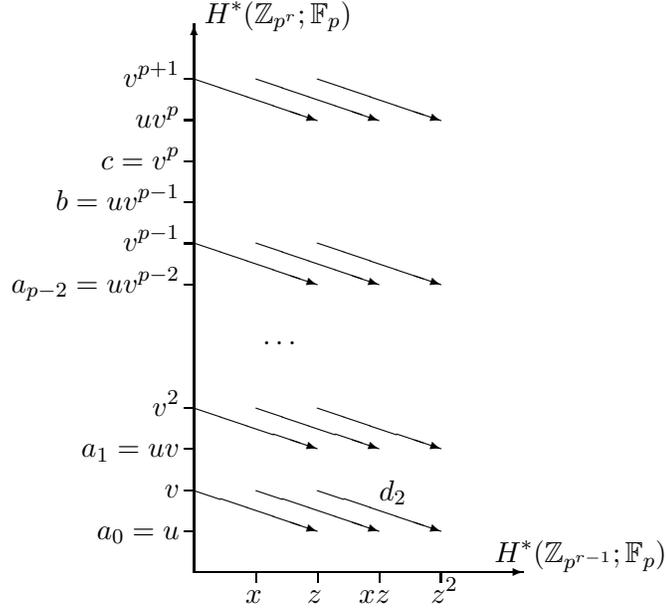

Since the rows of \eqref{p-diagram} are split, all the differentials landing on the horizontal axis are zero. In particular,
$d_2(u)=0.$

There is a possible non-zero differential $$d_2(v)=(\mbox{scalar})uz.$$
The fact that $\displaystyle H^2(G; {\mathbb F}_p) \cong \moplus_2{\mathbb F}_p$ implies that $d_2(v)\ne0$, i.e.
$$d_2(v)=muz$$
where $(m,p)=1$. The coefficient $m$ depends on the choice of the generator $v$ for $H^2({\mathbb Z}_{p^r};{\mathbb F}_p).$
It can be assumed that $m\equiv 1(\mbox{mod } p)$ for some choice of $v.$

Then
$$d_2(v^q x^i z^j) = q u v^{q-1} x^i z^{j+1},$$
$$d_2(uv^q x^i z^j) = 0.$$
Notice that $d_2(v^q x^i z^j) \ne 0$ if and only if $(q,p)=1.$   
Now we see that the rank of $\displaystyle \moplus_{i+j=q} E_3^{i,j}$ equals the rank of $H^q(G; {\mathbb F}_p)$ for
each
$q\ge 1,$ therefore all higher differentials are zero.
  
The mod $p$ cohomology ring of $G$ was computed by the author using Bockstein homomorphisms (since we already know the 
integer cohomology), but the calculation is long and is omitted here since this cohomology ring is worked out in \cite{H}. 
More precisely, this is a special case of theorem C, part (1). 

Let $a_i=uv^i (0 \le i \le p-2),$ $b=uv^{p-1},$ $c=v^p.$ Then
$$H^*(G; {\mathbb F}_p) \cong \Lambda(a_0, a_1, \ldots, a_{p-2}, b, x) \otimes {\mathbb F}_p[c,z]/ a_i a_j = a_i b = a_i z = 
0$$
where $|a_i|=2i+1,$ $|b|=2p-1,$ $|c|=2p,$ $|x|=1,$ and $|z|=2.$

Now consider the LHS spectral sequence for the middle column in \eqref{p-diagram}:
$$E_2^{p,q} \cong H^p({\mathbb Z}_{p-1}; H^q(G,{\mathbb F}_p))$$
Therefore $E_2^{p,q}=0$ for $p>0$, and $E_2^{0,q}\cong \left( H^q(G; {\mathbb F}_p) \right) ^{{\mathbb Z}_{p-1}}.$
The action of ${\mathbb Z}_{p-1}$ on $H^q(G,{\mathbb F}_p)$ is given by
$$g_2(x)=x, \ \ g_2(z)=z, \ \ g_2(a_i)=s^{(i+1)p^{r-1}}a_i, \ \ g_2(b)=s^{p^r}b, \ \ g_2(c)=s^{p^r}c.$$
It is easy to write out all the invariants under the action of $g_2$: the ring of invariants is generated by 
$x,$ $z,$ $a_{p-i}c^{i-2}$ ($2 \le i \le p$), $c^{p-1},$ and $bc^{p-2}.$

\begin{figure}
\centerline{\setlength{\unitlength}{3400sp}%
\begingroup\makeatletter\ifx\SetFigFont\undefined%
\gdef\SetFigFont#1#2#3#4#5{%
  \reset@font\fontsize{#1}{#2pt}%
  \fontfamily{#3}\fontseries{#4}\fontshape{#5}%
  \selectfont}%
\fi\endgroup%
\begin{picture}(4400,5005)(-60,-3236)
\thinlines
\special{ps: gsave 0 0 0 setrgbcolor}\put(1201,-2761){\vector( 1, 0){2400}}
\special{ps: grestore}\special{ps: gsave 0 0 0 setrgbcolor}\put(1201,-2761){\vector( 0, 1){4000}}
\special{ps: grestore}\special{ps: gsave 0 0 0 setrgbcolor}\put(1126,-1861){\line( 1, 0){ 75}}
\special{ps: grestore}\special{ps: gsave 0 0 0 setrgbcolor}\put(1126,-1561){\line( 1, 0){ 75}}
\special{ps: grestore}\special{ps: gsave 0 0 0 setrgbcolor}\put(1126,-661){\line( 1, 0){ 75}}
\special{ps: grestore}\special{ps: gsave 0 0 0 setrgbcolor}\put(1126,-361){\line( 1, 0){ 75}}
\special{ps: grestore}\special{ps: gsave 0 0 0 setrgbcolor}\put(1126,539){\line( 1, 0){ 75}}
\special{ps: grestore}\special{ps: gsave 0 0 0 setrgbcolor}\put(1126,839){\line( 1, 0){ 75}}
\special{ps: grestore}\special{ps: gsave 0 0 0 setrgbcolor}\put(1651,-2761){\line( 0,-1){ 75}}
\special{ps: grestore}\special{ps: gsave 0 0 0 setrgbcolor}\put(2101,-2761){\line( 0,-1){ 75}}
\special{ps: grestore}\special{ps: gsave 0 0 0 setrgbcolor}\put(2551,-2761){\line( 0,-1){ 75}}
\special{ps: grestore}\special{ps: gsave 0 0 0 setrgbcolor}\put(3001,-2761){\line( 0,-1){ 75}}
\special{ps: grestore}\special{ps: gsave 0 0 0 setrgbcolor}\put(1201,-1561){\vector( 3,-1){900}}
\special{ps: grestore}\special{ps: gsave 0 0 0 setrgbcolor}\put(1651,-1561){\vector( 3,-1){900}}
\special{ps: grestore}\special{ps: gsave 0 0 0 setrgbcolor}\put(2101,-1561){\vector( 3,-1){900}}
\special{ps: grestore}\special{ps: gsave 0 0 0 setrgbcolor}\put(1201,-361){\vector( 3,-1){900}}
\special{ps: grestore}\special{ps: gsave 0 0 0 setrgbcolor}\put(1651,-361){\vector( 3,-1){900}}
\special{ps: grestore}\special{ps: gsave 0 0 0 setrgbcolor}\put(2101,-361){\vector( 3,-1){900}}
\special{ps: grestore}\special{ps: gsave 0 0 0 setrgbcolor}\put(1201,-661){\line(-1, 0){ 75}}
\put(1096,-1916){\makebox(0,0)[rb]{\smash{\SetFigFont{11}{13.2}{\rmdefault}{\mddefault}{\updefault}\special{ps: 
gsave 0 0 0 setrgbcolor}{\small $d_0=uv^{p-2}$}\special{ps: grestore}}}}
\put(1096,-1616){\makebox(0,0)[rb]{\smash{\SetFigFont{11}{13.2}{\rmdefault}{\mddefault}{\updefault}\special{ps: 
gsave 0 0 0 setrgbcolor}{\small $v^{p-1}$}\special{ps: grestore}}}}
\put(1096,-716){\makebox(0,0)[rb]{\smash{\SetFigFont{11}{13.2}{\rmdefault}{\mddefault}{\updefault}\special{ps: 
gsave 0 0 0 setrgbcolor}{\small $d_1=uv^{2p-3}$}\special{ps: grestore}}}}
\put(1096,-416){\makebox(0,0)[rb]{\smash{\SetFigFont{11}{13.2}{\rmdefault}{\mddefault}{\updefault}\special{ps:
gsave 0 0 0 setrgbcolor}{\small $v^{2p-2}$}\special{ps: grestore}}}}
\put(1096,484){\makebox(0,0)[rb]{\smash{\SetFigFont{11}{13.2}{\rmdefault}{\mddefault}{\updefault}\special{ps:
gsave 0 0 0 setrgbcolor}{\small $e=uv^{p^2-p-1}$}\special{ps: grestore}}}}
\put(1096,784){\makebox(0,0)[rb]{\smash{\SetFigFont{11}{13.2}{\rmdefault}{\mddefault}{\updefault}\special{ps:
gsave 0 0 0 setrgbcolor}{\small $f=v^{p^2-p}$}\special{ps: grestore}}}}
\put(1576,-2986){\makebox(0,0)[lb]{\smash{\SetFigFont{11}{13.2}{\rmdefault}{\mddefault}{\updefault}\special{ps: 
gsave 0 0 0 setrgbcolor}{\small $x$}\special{ps: grestore}}}}
\put(2026,-2986){\makebox(0,0)[lb]{\smash{\SetFigFont{11}{13.2}{\rmdefault}{\mddefault}{\updefault}\special{ps: 
gsave 0 0 0 setrgbcolor}{\small $z$}\special{ps: grestore}}}}
\put(2401,-2986){\makebox(0,0)[lb]{\smash{\SetFigFont{11}{13.2}{\rmdefault}{\mddefault}{\updefault}\special{ps: 
gsave 0 0 0 setrgbcolor}{\small $xz$}\special{ps: grestore}}}}
\put(2926,-2986){\makebox(0,0)[lb]{\smash{\SetFigFont{11}{13.2}{\rmdefault}{\mddefault}{\updefault}\special{ps: 
gsave 0 0 0 setrgbcolor}{\small $z^2$}\special{ps: grestore}}}}
\put(1276,1249){\makebox(0,0)[lb]{\smash{\SetFigFont{11}{13.2}{\rmdefault}{\mddefault}{\updefault}\special{ps: 
gsave 0 0 0 setrgbcolor}{\small $H^*({\mathbb Z}_{p^r};{\mathbb F}_p)$}\special{ps: grestore}}}}
\put(3390,-2686){\makebox(0,0)[lb]{\smash{\SetFigFont{11}{13.2}{\rmdefault}{\mddefault}{\updefault}\special{ps: 
gsave 0 0 0 setrgbcolor}{\small $H^*({\mathbb Z}_{(p-1)p^{r-1}};{\mathbb F}_p)$}\special{ps: grestore}}}}
\put(2551,-1636){\makebox(0,0)[lb]{\smash{\SetFigFont{11}{13.2}{\rmdefault}{\mddefault}{\updefault}\special{ps: 
gsave 0 0 0 setrgbcolor}{\small $d_2$}\special{ps: grestore}}}}
\put(1500,600){\makebox(0,0)[lb]{\smash{\SetFigFont{11}{13.2}{\rmdefault}{\mddefault}{\updefault}\special{ps:
gsave 0 0 0 setrgbcolor}{\small $d_2=0$}\special{ps: grestore}}}}
\put(1700,100){\makebox(0,0)[lb]{\smash{\SetFigFont{11}{13.2}{\rmdefault}{\mddefault}{\updefault}\special{ps:
gsave 0 0 0 setrgbcolor}{\small \ldots}\special{ps: grestore}}}}
\end{picture}}
\caption{$E_2$ stage of the LHS spectral sequence for $H^*(\Hol(\Z_{p^r});{\mathbb F}_p)$, $r\ge 2$}
\label{Lyndon_Hol_Zpr_inv}
\end{figure}
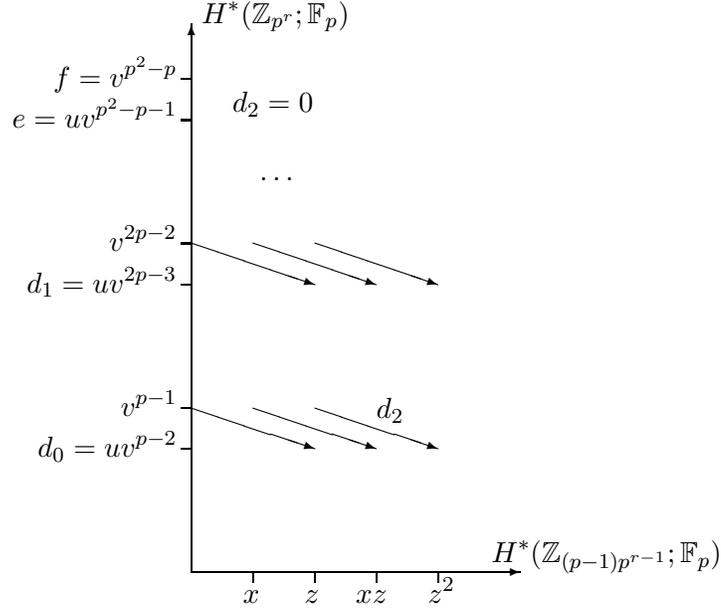

Let $d_i=a_{p-2-i}c^i$ ($0 \le i \le p-2$), $e=bc^{p-2},$ and $f=c^{p-1}.$ Then we have

\begin{theorem}\label{coh-ring-podd}
$$H^*(\mbox{Hol}({\mathbb Z}_{p^r}); {\mathbb F}_p) \cong \Lambda (d_0, d_1, \ldots, d_{p-2}, e, x) \otimes {\mathbb 
F}_p[f,z]/
d_i d_j = d_i e = d_i z = 0$$
where $|d_i|=2(p-1)(i+1)-1,$ $|e|=2p(p-1)-1,$ $|f|=2p(p-1),$ $|x|=1,$ $|z|=2.$
\end{theorem}

This also follows from theorem A, part (2), of \cite{H}. 

The LHS spectral sequence for $H^*(\Hol(\Z_{p^r});{\mathbb F}_p)$ is shown in figure \ref{Lyndon_Hol_Zpr_inv} and will be used in 
section \ref{wilson_section}. 

\section{Bockstein homomorphisms}

In this section we use the integer cohomology of $\Hol({\mathbb Z}_{p^r})$ (which can be obtained from the last formula of
section 2 using the universal coefficients formula) to compute the Bockstein homomorphisms in the Bockstein spectral sequence   
for $H^*(\mbox{Hol}({\mathbb Z}_{p^r});{\mathbb F}_p)$.

Cases $r=1$ and $r=2$ are special, so let $r \ge 3$.

Since $H^2(\Hol(\Z_{p^r});{\mathbb Z}) \cong {\mathbb Z}_{p^{r-1}} \oplus {\mathbb Z}_{p-1},$ there exists a cycle in
$E^2_{r-1}$ which is in the
image of $\beta _{r-1}.$ It is easy to see that $H^2(\mbox{Hol}({\mathbb Z}_{p^r});{\mathbb F}_p)$ is generated by $z$ and
$H^1(\Hol(\Z_{p^r});{\mathbb F}_p)$ is generated by $x$, so $$\beta_{r-1}(x)=c_x z, \ \ \ c_x\ne 0 \ (\mbox{mod } p).$$
Then $\beta_{r-1}(z)=0$ because $\beta_{r-1}^2=0$, $\beta_{r-1}(xz^m)=c_x z^{m+1}$, and $\beta_n(x)=\beta_n(z)=0$ for 
$n<r-1.$

Since $\Z_p$ is a direct summand of $H^{2(p-1)}(\mbox{Hol}({\mathbb Z}_{p^r});{\mathbb Z})$, and $\{ d_0x, \ z^{p-1} \}$ is a
basis for $H^{2(p-1)}\Hol(\Z_{p^r});{\mathbb F}_p)$, some non-zero linear combination of $d_0 x$ and $z^{p-1}$ is in the
image of $\beta_1$.
$\{ d_0, \ xz^{p-2} \} $ is a basis for $H^{2(p-1)-1}(\mbox{Hol}({\mathbb Z}_{p^r});{\mathbb F}_p)$, and we know that
$\beta_1(x)=\beta_1(z)=0$, so $\beta_1(d_0)\ne 0$. Let $\beta_1(d_0) = c_{d_0} d_0 x + c_{d_0}' z^{p-1}$. Then
$\beta_1(d_0 x) = c_{d_0}' z^{p-1}x,$ but $H^{2(p-1)+1}(\Hol(\Z_{p^r});{\mathbb Z})$ does not have a ${\mathbb Z}_p$
summand, therefore $c_{d_0}'=0.$ Thus $$\beta_1(d_0) = c_{d_0} d_0 x, \ \ \ c_{d_0}\ne 0 \ (\mbox{mod } p).$$

Similarly, $$\beta_1(d_i) = c_{d_i} d_i x, \ \ \ c_{d_i}\ne 0 \ (\mbox{mod } p), \ \ \ 0\le i\le p-2.$$

Further, $\{e, xz^{p(p-1)-1} \}$ is a basis for $H^{2p(p-1)-1}(\mbox{Hol}({\mathbb Z}_{p^r});{\mathbb F}_p)$.
If $r>3$, then ${\mathbb Z}_{p^2}$ is a direct summand of $H^{2p(p-1)}(\mbox{Hol}({\mathbb Z}_{p^r});{\mathbb Z})$,
and $\beta_2(xz^{p(p-1)-1})=0$, therefore $\beta_2(e)\ne 0$.
If $r=3$, then ${\mathbb Z}_{p^2} \oplus {\mathbb Z}_{p^2}$ is a direct summand of
$H^{2p(p-1)}(\mbox{Hol}({\mathbb Z}_{p^r});{\mathbb Z})$, and $\beta_2(xz^{p(p-1)-1}) = c_x z^{p(p-1)}$, and again we have 
that
$\beta_2(e)$ must be non-zero.

Since ${\mathbb Z}_p$ is a direct summand of $H^{2p(p-1)+1}(\mbox{Hol}({\mathbb Z}_{p^r});{\mathbb Z})$,
$\{f, ex, z^{p(p-1)}\}$ is a basis for $H^{2p(p-1)}(\mbox{Hol}({\mathbb Z}_{p^r});{\mathbb F}_p)$, and
$\beta_1(e)=\beta_1(x)=\beta_1(z)=0$, we conclude that $\beta_1(f)\ne 0$.
Let $\beta_1(f)=c_f ez + c_f' fx + c_f'' xz^{p(p-1)}$. Also, ${\mathbb Z}_p$ is a direct
summand of $H^{2p(p-1)+2}(\mbox{Hol}({\mathbb Z}_{p^r});{\mathbb Z})$, and $\{ fx, ez, xz^{p(p-1)} \}$ is a basis for
$H^{2p(p-1)+1}(\mbox{Hol}({\mathbb Z}_{p^r});{\mathbb F}_p)$, therefore $\beta_1(fx)=c_f ezx \ne 0$, so $c_f \ne 0$.

Now back to $\beta_2(e)$. Since $\beta_1(f)\ne 0$, $f$ doesn't survive to the second stage of the Bockstein spectral 
sequence,
therefore $\beta_2(e)=c_e ex + c_e' z^{p(p-1)}$.

If $r>3$, then $\beta_2(x)=0$ and thus $\beta_2(ex)=c_e' z^{p(p-1)} x =0 $ because ${\mathbb Z}_{p^2}$ is not a direct
summand of $H^{2p(p-1)+1}(\mbox{Hol}({\mathbb Z}_{p^r});{\mathbb Z})$. Thus $c_e'=0$, and we have

\begin{equation} \label{beta2e} \beta_2(e)=c_e ex, \ \ \ c_e\ne 0  \ (\mbox{mod } p). \end{equation}

$\beta_1(fe)=\beta_1(f)\cdot e = c_f' fxe + c_f''xz^{p(p-1)}e$. ${\mathbb Z}_p$ but not ${\mathbb Z}_p \oplus {\mathbb
Z}_p$ is a direct summand of $H^{4p(p-1)}(\mbox{Hol}({\mathbb Z}_{p^r});{\mathbb Z})$, and
$\beta_1(fxz^{p(p-1)-1})=c_f exz^{p(p-1)}$ where $c_f\ne 0$. Therefore $c_f'=0$. So $\beta_1(f)=c_f ez + c_f'' xz^{p(p-1)}$.
Let $e_{new}=e+\frac{c_f''}{c_f}xz^{p(p-1)-1}$. Then we'll have
$\beta_2(e_{new})= c_e ex = c_e e_{new}x$, so we'll replace $e$ by $e_{new}$ and will right just $e$ again.
This gives \begin{equation} \label{beta1f} \beta_1(f)=c_f ez, \ \ \ c_f\ne 0 \ (\mbox{mod } p). \end{equation}

If $r=3$, then $\beta_2(ex)=-c_x ez + c_e' z^{p(p-1)}x$ (where $c_x\ne 0$) must be equal to 0 module the image of $\beta_1$,
and the only way this may happen is if $\beta_1(f) = m(-c_x ez + c_e' z^{p(p-1)}x) =  c_f ez + c_f'' z^{p(p-1)}x$. Again, let
$e_{new}=e+\frac{c_f''}{c_f}xz^{p(p-1)-1}$. Then $\beta_2(e_{new})= c_e e_{new}x$ and $\beta_1(f)=c_f e_{new}z$, so after
replacing $e$ by $e_{new}$ we'll have \eqref{beta2e} and \eqref{beta1f}.

Notice that \eqref{beta1f} implies that $\beta_1(ef)=0$. Moreover, it is easy to compute $\beta_1$ of any product of powers
of $d_i,$ $e,$ $x,$ $f,$ and $z,$ and see that only elements of the form $d_i x f^n$, $ef^n$, $x$, $z$, $f^p$ (where $0\le n<
p$) and their products are mapped to 0 under $\beta_1$. All higher bocksteins of $z$ and $d_i x$ are zero because there are
no ${\mathbb Z}_{p^*}$ summands in the integer cohomology of $\Hol({\mathbb Z}_{p^r})$ of degrees one more than the degrees
of $z$ and $d_i x$. All Bockstein homomorphisms of  $x$ we already know.

If $r>3$, then ${\mathbb Z}_{p^2}$ is a direct summand of $H^{4p(p-1)}(\mbox{Hol}({\mathbb Z}_{p^r});{\mathbb Z})$ and
$\{ ef, $ $ez^{p(p-1)},$ $xz^{2p(p-1)-1} \}$ is a basis for the kernel of $\beta_1$ in degree $4p(p-1)-1$; but
$\beta_2(ez^{p(p-1)}) = \beta_2(xz^{2p(p-1)-1}) = 0$, so $\beta_2(ef) \ne 0$.
Let $\beta_2(ef) = c_{ef} efx + c_{ef}'z^{2p(p-1)}$. (The terms with $f^2$, $fz^{p(p-1)}$, and $exz^{p(p-1)}$ are not present
in this linear combination because $\beta_1$ of the first two elements is non-zero, and the third one is in the image of
$\beta_1$.)
Then $\beta_2(efx)=c_{ef}'z^{2p(p-1)}x=0$ (because ${\mathbb Z}_{p^2}$ is not a direct summand of
$H^{4p(p-1)+1}(\Hol(\Z_{p^r});{\mathbb Z})$), and $z^{2p(p-1)}x$ is not in the image of $\beta_1$, therefore
$c_{ef}'=0$ and we have
\begin{equation} \label{beta2ef} \beta_2(ef)=c_{ef} efx, \ \ \ c_{ef}\ne 0 \ (\mbox{mod } p). \end{equation}

If $r=3$, then $\Z_{p^2}\oplus\Z_{p^2}$ is a direct summand of $H^{4p(p-1)}(\mbox{Hol}({\mathbb Z}_{p^r});{\mathbb Z})$.
Now $\beta_2(ez^{p(p-1)}) = 0$ and $\beta_2(xz^{2p(p-1)-1}) = c_x z^{2p(p-1)}$, so we still have $\beta_2(ef) \ne 0$. Let
$\beta_2(ef) = c_{ef} efx + c_{ef}'z^{2p(p-1)}$. Then $\beta_2(efx)=c_{ef}'z^{2p(p-1)}x - c_x efz$, but $efz$ is in the image
of $\beta_1$, since $\beta_1(f^2)=2f c_f ez$. $z^{2p(p-1)}x$ is not in the image of $\beta_1$, therefore we still have
$c_{ef}'=0$ and hence \eqref{beta2ef}.

Similarly, $$\beta_2(ef^n)=c_{ef^n} ef^n x, \ \ \ c_{ef^n}\ne 0 \ (\mbox{mod } p), \ 2 \le n \le p-2.$$

If $r\ne 4$, then $\Z_{p^3}$ is a direct summand of $H^{2p^2(p-1)}(\mbox{Hol}(\Z_{p^r});{\mathbb Z})$, so there must be 
an element in $E^{2p^2(p-1)-1}_3$ whose $\beta_3$ is non-zero. If $r=4$, then $\Z_{p^3} \oplus {\mathbb Z}_{p^3}$ is a direct 
summand of $H^{2p^2(p-1)}(\mbox{Hol}({\mathbb Z}_{p^r});\Z)$; in this case $\beta_3(xz^{p^2(p-1)-1})\ne 0$, but there must be 
another such element. Here is a basis for $H^{2p^2(p-1)-1}(\Hol({\mathbb Z}_{p^r});{\mathbb F}_p)$: $ef^{p-1},$
$ef^{p-1-i}z^{ip(p-1)}$ ($1\le i\le p-1$), $f^{p-i}xz^{ip(p-1)-1}$ ($1\le i \le p-1$), $xz^{p^2(p-1)-1}$. However,
$\beta_3(ef^{p-1-i}z^{ip(p-1)} =0$ (it is easy to see that $\beta_3(ef^nz)=0$ and, as mentioned above, $\beta_3(z)=0$),
$\beta_1(f^{p-i}xz^{ip(p-1)-1)}=c_f(p-i)f^{p-i-1}exz^{ip(p-1)}\ne 0$, and in the case $r>4$, $\beta_3(xz^{p^2(p-1)-1})=0$
while in the case $r=3$ $\beta_3(xz^{p^2(p-1)-1})$ is not even defined. So $\beta_3(ef^{p-1})\ne 0$.

Since $\Z_{p^2}$ is a direct summand of $H^{2p^2(p-1)+1}(\mbox{Hol}({\mathbb Z}_{p^r});{\mathbb Z})$, there must be a cycle 
in $E^{2p^2(p-1)-1}_3$ which is in the image of $\beta_2$. An analysis of a basis for
$H^{2p^2(p-1)}(\Hol(\Z_{p^r});{\mathbb F}_p)$ similar to the preceding paragraph shows that the only possibility is that 
$\beta_2(f^p)\ne 0$. An analysis of a basis for $H^{2p^2(p-1)+1}(\mbox{Hol}({\mathbb Z}_{p^r});{\mathbb F}_p)$ gives that
$\beta_2(f^p) = c_{f^p} ef^{p-1}z + c_{f^p}' f^px + c_{f^p}'' xz^{p^2(p-1)}$. Also, if $r>3$, then 
$\beta_2(f^px)= c_{f^p} ef^{p-1}zx \ne 0$ (I bet nobody will ever read this mess, but just in case I will say that this is because 
$\Z_{p^2}$ is a direct summand of $H^{2p^2(p-1)+2}(\mbox{Hol}({\mathbb Z}_{p^r});\Z)$), so $c_{f^p}\ne 0$. 
If $r=3$, then ${\mathbb Z}_{p^2} \oplus {\mathbb Z}_{p^2}$ is a summand of
$H^{2p^2(p-1)+2}(\mbox{Hol}({\mathbb Z}_{p^r});\Z)$, so in addition to $\beta_2(xz^{p^2(p-1)})=z^{p^2(p-1)+1}$ we must have 
$\beta_2(f^px)=c_{f^p} ef^{p-1}zx + c_x f^pz \ne 0$. But $\beta_2(\beta_2(xz^{p^2(p-1)}))=0$ which implies that $c_{f^p}\ne 
0$,
$c_{f^p}'=0$, and $c_{f^p}''=0$. Thus if $r=3$ then

\begin{equation} \label{beta2fp} \beta_2(f^p) =  c_{f^p} ef^{p-1}z, \ \ \ c_{f^p} \ne 0 \ (\mbox{mod } p). \end{equation}

Now look at $\beta_3(ef^{p-1})$ again. Since $\beta_2(f^p)\ne 0$, $f^p$ doesn't survive to the third stage of the Bockstein
spectral sequence, therefore $\beta_3(ef^{p-1}) = c_{ef^{p-1}} ef^{p-1}x +  c_{ef^{p-1}}' z^{p^2(p-1)}$. If $r>4$, then
$\beta_3(ef^{p-1}x)=c_{ef^{p-1}}' z^{p^2(p-1)}x =0$ therefore $c_{ef^{p-1}}'=0$, and we have

\begin{equation} \label{beta3efn} 
\beta_3(ef^{p-1}) = c_{ef^{p-1}} ef^{p-1}x, \ \ \ c_{ef^{p-1}} \ne 0 \ (\mbox{mod } p). 
\end{equation}

If $r=4$, then $\beta_3(ef^{p-1}x)=c_{ef^{p-1}}' z^{p^2(p-1)}x - c_x ef^{p-1}z $ (where $c_x \ne 0$) must be equal to 0
module the image of $\beta_2$, and the only way this can happen is if this is a multiple of $\beta_2(f^p) = c_{f^p}
ef^{p-1}z$. Therefore, we have $c_{ef^{p-1}}'=0$ and thus \eqref{beta3efn}.

If $r=3$, then $z^{p^2(p-1)}$ is in the image of $\beta_2$, so we have \eqref{beta3efn} again.

Now we wish to show that \eqref{beta2fp} holds even if $r>3$.
$\beta_2(f^p\cdot ef^{p-1}) = (c_{f^p}' exf^{2p-1} + c_{f^p}''exz^{p^2(p-1)} f^{p-1}$. ${\mathbb Z}_{p^2}$ but not
${\mathbb Z}_{p^2} \oplus {\mathbb Z}_{p^2}$ is a summand of \\ $H^{4p^2(p-1)}(\mbox{Hol}({\mathbb Z}_{p^r});{\mathbb Z})$, and
$\beta_2(xz^{p^2(p-1)-1}f^p = c_{f^p}exz^{p^2(p-1)} f^{p-1}$, therefore $c_{f^p}'=0$.
Similarly, it can be shown that $c_{f^p}''=0$, and
$$ \beta_k(ef^{p^{k-2}-1}) = c_{ef^{p^{k-2}-1}} ef^{p^{k-2}-1}x, \ \ \ c_{ef^{p^{k-2}-1}} \ne 0 \ (\mbox{mod } p) $$
and
$$ \beta_{k-1}(f^{p^{k-2}}) =  c_{f^{p^{k-2}}} ef^{p^{k-2}-1}z, \ \ \ c_{f^{p^{k-2}}} \ne 0 \ (\mbox{mod } p). $$

\section{Continuous cohomology of $\Hol(\Z_{p^r})$}

Recall (see \cite{S}) that if we have a group $\displaystyle G=\lim_{\leftarrow}G_i,$ and coefficients 
$\displaystyle A=\lim_{\rightarrow}A_j$ where the limits are taken over surjections $p_i:G_{i+1}\rightarrow G_i$ and injections 
$i_j:A_j\rightarrow A_{j+1}$, and the actions of $G_n$'s on $A_n$'s are compatible, that is, $(p_n g_{n+1})(a_n) = g_{n+1}(i_n(a_n))$ 
for all $a_n\in A_n$ and $g_{n+1}\in G_{n+1}$, then continuous cohomology of $G$ with coefficients in $A$ is defined by 
$$H_{cont}^*(G;A) = \lim_{\rightarrow} (H^*(G_n;A_n)).$$

It is well-known that for $G_n=\Z_{p^n}$ with natural surjections and $A_n=\Z_p$ with $i_n$ isomorphism for each $n$, 
$H_{cont}^* (G;A) = \Lambda(x)$, $|x|=1$. 

Let $G_n=\Hol(\Z_{p^n})$ with surjections from example \ref{characteristic-example-Z}, and $A_n=\Z_p$ with $i_n$ 
isomorphism for each $n$. 
Then the LHS spectral sequence for $\Hol(\Z_{p^n})$ maps into that for $\Hol(\Z_{p^{n+1}})$. For $p=2$ (see figure 
\ref{Lyndon_Hol_Z2r}), $x_n \mapsto x_{n+1}$, $y_n \mapsto y_{n+1}$, $z_n \mapsto 0$, $u_n \mapsto u_{n+1}$, and 
$v_n \mapsto 0$. For $p$ odd (see figure \ref{Lyndon_Hol_Zpr_inv}), $x_n \mapsto x_{n+1}$, $z_n \mapsto 0$, 
$u_n \mapsto u_{n+1}$, and $v_n \mapsto 0$. Therefore we have 

\begin{theorem}\label{my-continuous-theorem}
$H_{cont}^* (\lim_{\leftarrow} ( \Hol(\Z_{p^n}) ); \F_p) = $

$ \hfill = \left\{ \begin{array}{@{}ll@{\,}}
\Lambda(x)\otimes \F_2[u,y]/u^2=ux+uy, & |x|=|u|=|y|=1, \ \hfill p=2, \\
\Lambda(x),                  & |x|=1,         \hfill p \mbox{ is odd.}  
\end{array}\right.$ \\
The Poincar\'e series of $H_{cont}^* (\lim_{\leftarrow} ( \Hol(\Z_{2^n}) ); \F_2)$ is $\frac{(1+t)^2}{1-t}$. 
\end{theorem}

\chapter{Holomorph of $\displaystyle \moplus_n \Z_{p^r} $}

\label{matrices-chapter}

\thispagestyle{kisa}

\section {Structure of $\displaystyle \Hol \left( \moplus_n R\right) $}

\label{matrices-section}

Let $R$ be a commutative ring with a unit. Then $\displaystyle \Aut\left(\moplus_n R\right)$, the group of automorphisms of 
$\displaystyle \moplus_n R$ which is regarded as a module over $R$, is isomorphic to $GL(n,R)$, the group of all invertible 
$n \times n$ matrices with entries in $R$. 
\begin{lemma} \label{hol-matrix-lemma}
\begin{enumerate}
\item The group $\displaystyle \Hol \left( \moplus_n R \right)$ is isomorphic to the subgroup of $GL(n+1, R)$
consisting of all the matrices of the form
\begin{equation} \renewcommand{\baselinestretch}{1} \normalsize \label{hol-matrix}
\left[
\begin{array}{c|c}
1      & * \hfill \ldots \hfill * \\
\hline
0      & n \times n        \\
\vdots & \mbox{invertible} \\
0      & \mbox{matrix}     \\
\end{array}
\right] , 
\end{equation} 
the group of automorphisms of $\displaystyle \moplus_{n+1} R$ fixing 
$$ \renewcommand{\baselinestretch}{1} \normalsize 
e_1 = \left[ \begin{array}{c} 1 \\ 0 \\ \vdots \\ 0 \end{array} \right] .$$
\item Alternatively, $\displaystyle \Hol \left( \moplus_n R \right)$ is isomorphic to the subgroup of $GL(n+1, R)$
consisting of all the matrices of the form

$$ \left[ \renewcommand{\baselinestretch}{1} \normalsize 
\begin{array}{c|c}
1      & 0 \hfill \ldots \hfill 0 \\
\hline  
*      & n \times n        \\
\vdots & \mbox{invertible} \\
*      & \mbox{matrix}     \\
\end{array}
\right] , $$ 

\noindent the stabilizer of a rank $n$ submodule of $\displaystyle \moplus_{n+1} R$.  
\end{enumerate}
\end{lemma}

\begin{proof}
\begin{enumerate}
\item Recall that $\displaystyle \Hol \left( \moplus_n R \right)$ is defined by the split extension
$$ 1 \rightarrow \moplus_n R \rightarrow \Hol \left( \moplus_n R \right)
\buildrel \curvearrowleft \over \rightarrow \Aut \left( \moplus_n R \right) \rightarrow 1,$$
where multiplication is defined by $(f,x)(g,y)=(fg,g^{-1}(x)+y)$ for
$\displaystyle f,g \in \Aut \left( \moplus_n R \right)$ and
$\displaystyle x,y \in \moplus_n R$.

Think of $\displaystyle x \in \moplus_n R$ as a vertical $n-$vector.
For $\displaystyle f \in \Aut \left( \moplus_n R \right)$, let $M_f$ be the $n\times n$ matrix such that
$f(x)=M_fx$ for every $\displaystyle x \in \moplus_n R$. There is a one-to-one correspondence between
elements $(f,x)$ of $\displaystyle \Hol \left( \moplus_n R \right)$ and matrices of the form
\eqref{hol-matrix}, namely 

$$ (f,x) \leftrightarrow
\left[
\begin{array}{c|c}     1 & x^T    \\    \hline    0 & (M_f^{-1})^T    \end{array}
\right] .$$

This correspondence is a group homomorphism: multiplication of elements in $\displaystyle \Hol \left( \moplus_n R \right)$ 
regarded as a subgroup of $GL(n+1, R)$ is as follows:

$ (f,x)(g,y) = 
 \left[
\begin{array}{c|c}     1 & x^T    \\    \hline    0 & (M_f^{-1})^T    \end{array}
\right]
\left[
\begin{array}{c|c}     1 & y^T    \\    \hline    0 & (M_g^{-1})^T    \end{array}
\right]
=
\left[
\begin{array}{c|c}     1 & y^T + x^T(M_g^{-1})^T     \\   \hline    0 & (M_f^{-1})^T (M_g^{-1})^T    \end{array}
\right]
=$

$= \left[
\begin{array}{c|c}     1 & (M_g^{-1}x+y)^T     \\   \hline    0 & ((M_fM_g)^{-1})^T     \end{array}
\right]
=(fg, g^{-1}(x)+y).$
\item An isomorphism is given by 
$$ (f,x) \leftrightarrow
\left[
\begin{array}{c|c}     1 & 0   \\    \hline    -M_fx & M_f    \end{array}
\right] .$$
\end{enumerate}
\end{proof}

\begin{remark}
This subgroup of $GL(n+1, R)$  isomorphic to $\displaystyle \Hol \left( \moplus_n R \right)$ 
is not normal.
\end{remark}

It is a classical fact that the stabilizer of a submodule $\displaystyle \moplus_n \Z_2 \subset \moplus_{n+1} \Z_2$ is a maximal 
subgroup of $GL(n+1, \Z_2)$ (see \cite{AFG}). We give a proof of this fact for completeness. 

For the rest of this section we will identify $\displaystyle \Hol \left( \moplus_n R \right)$ with the subgroup of 
$GL(n+1, R)$ consisting of matrices of the form \eqref{hol-matrix}. 

\begin{corollary} $\displaystyle \Hol \left( \moplus_n \Z_2 \right)$ regarded as a subgroup of $GL(n+1, \Z_2)$ as 
above is maximal. 
\end{corollary}
\begin{proof}

First we show that $\displaystyle \Hol \left( \moplus_n \Z_2 \right)$ and the permutation matrix 
$$\renewcommand{\baselinestretch}{1}\normalsize P_{1,2} = \left[ \begin{array}{c|c|ccc}  
		& 1 & &   & \\ 
\hline
	      1 &   & &   & \\
\hline
                &   & &   & \\
                &   & & I & \\
                &   & &   &  		
\end{array}\right] $$ 
(where missing cells are 0) generate all of $GL(n+1, \Z_2)$. Let $A\in GL(n+1, \Z_2).$ If the first column of $A$ is 
$e_1$, then 
$A\in \Hol \left( \moplus_n \Z_2 \right)$. If the first column of $A$ has the form 
$\renewcommand{\baselinestretch}{1}\normalsize 
A_1 = \left[ \begin{array}{c} a \\ 1 \\ 0 \\ \vdots \\ 0 \end{array} \right]$, 
then 
$$\renewcommand{\baselinestretch}{1}\normalsize A = \left[ \begin{array}{c|ccc} 
	 	a & & * & \\
\hline
		1 & & * & \\
\hline
                  & &   & \\
                  & & * & \\
                  & &   &         \end{array}\right]  =
\left[ \begin{array}{c|c|ccc}
		1 & a & &   & \\
\hline
		  & 1 & &   & \\
\hline
 		  &   & &   & \\
	 	  &   & & I & \\
		  &   & &   &      \end{array}\right] 
        \left[ \begin{array}{c|c|ccc}
                  & 1 & &   & \\
\hline
                1 &   & &   & \\
\hline
                  &   & &   & \\
                  &   & & I & \\
                  &   & &   &      \end{array}\right] 				
\left(
        \left[ \begin{array}{c|c|ccc}
                  & 1 & &   & \\
\hline
                1 & a & &   & \\
\hline
                  &   & &   & \\
                  &   & & I & \\
                  &   & &   &      \end{array}\right] 
A \right) .$$ 
\noindent The first matrix in this decomposition is in 
$\displaystyle \Hol \left( \moplus_n \Z_2 \right)$, the second one is our permutation matrix $P_{1,2}$, and it is easy to 
check that the product in parenthesis is again in $\displaystyle \Hol \left( \moplus_n \Z_2 \right)$ as it has the 
first column $e_1$. 

Here is a geometric interpretation of this decomposition. Given a basis in $\moplus_{n+1} \Z_2$ with the first 
vector of the form $A_1$, there is a sequence of linear 
transformations each of which either is given by a matrix in 
$\displaystyle \Hol \left( \moplus_n \Z_2 \right) \hookrightarrow GL(n+1, \Z_2)$ or 
permutes $e_1$ and $e_2$ and fixes all the other basis vectors, such that the composite of the transformations in 
this sequence transforms the standard basis into the given one. This sequence is as follows. First transform $e_2$ 
into $A_1$ and leave all the other basis vectors 
fixed. Then  interchange the first two basis vectors. Finally, leaving the first vector fixed, transform all the 
other vectors as needed to obtain the given basis. 

If the $i$-th entry of $A$ is 1 for some $i>2$, then 

$\renewcommand{\baselinestretch}{1}\normalsize
A = \left[ \begin{array}{c|ccccccc}
                a & & & & * & & & \\
\hline
                b & & & & * & & & \\
\hline
                  & & & &   & & & \\
                x & & & & * & & & \\
                  & & & &   & & & \\
\hline  
                1 & & & & * & & & \\
\hline      
                  & & & &   & & & \\
                y & & & & * & & & \\
                  & & & &   & & &       \end{array}\right] =
\left[ \begin{array}{c|c|ccc|c|ccc}
		1 &   & &   & & a & &   & \\
\hline
		  & 1 & &   & & b & &   & \\
\hline
		  &   & &   & &   & &   & \\
		  &   & & I & & x & &   & \\
		  &   & &   & &   & &   & \\
\hline
		  &   & &   & & 1 & &   & \\
\hline
		  &   & &   & &   & &   & \\
		  &   & &   & & y & & I & \\
                  &   & &   & &   & &   & \end{array}\right] \times $

$$\renewcommand{\baselinestretch}{1}\normalsize \left[ \begin{array}{c|c|ccc|c|ccc@{}}
                1 &   & &   & &   & &   & \\
\hline
                  &   & &   & & 1 & &   & \\
\hline
                  &   & &   & &   & &   & \\
                  &   & & I & &   & &   & \\
                  &   & &   & &   & &   & \\
\hline
                  & 1 & &   & &   & &   & \\
\hline
                  &   & &   & &   & &   & \\
                  &   & &   & &   & & I & \\
                  &   & &   & &   & &   & \end{array}\right]    
\left[ \begin{array}{c|c|ccc@{\hspace{1mm}}c@{\hspace{1mm}}ccc}
                  & 1 & &   & &   & &   & \\
\hline
                1 &   & &   & &   & &   & \\
\hline
                  &   & &   & &   & &   & \\
                  &   & &   & &   & &   & \\
                  &   & &   & &   & &   & \\
                  &   & &   & & I & &   & \\
                  &   & &   & &   & &   & \\
                  &   & &   & &   & &   & \\
                  &   & &   & &   & &   & \end{array}\right]
\left(
\left[ \begin{array}{c|c|ccc|c|ccc@{}}
                  &   & &   & & 1 & &   & \\
\hline
                1 &   & &   & & a & &   & \\
\hline
                  &   & &   & &   & &   & \\
                  &   & & I & & x & &   & \\
                  &   & &   & &   & &   & \\
\hline
                  & 1 & &   & & b & &   & \\
\hline
                  &   & &   & &   & &   & \\
                  &   & &   & & y & & I & \\
                  &   & &   & &   & &   & \end{array}\right]
A \right)$$ 

The first two matrices (and therefore their product, as will be needed below) are in 
$\displaystyle \Hol \left( \moplus_n \Z_2 \right)$, the third one is $P_{1,2}$, and, again, it is easy to check that the product in 
parenthesis is in $\displaystyle \Hol \left( \moplus_n \Z_2 \right)$. A geometric interpretation is similar (but now 
we added one step, namely, interchange the second and $i$-th vectors). 

Thus we proved that for any matrix $A\in GL(n+1, \Z_2)$ there exist matrices $B$ and $C$ in 
$\displaystyle \Hol \left( \moplus_n \Z_2 \right)$ such that $A=BP_{1,2}C$. Now let $M$ also be any matrix in $GL(n+1, \Z_2)$. Then 
there exist matrices $K$ and $L$ in $\displaystyle \Hol \left( \moplus_n \Z_2 \right)$ such that $M=KP_{1,2}L$, therefore 
$P_{1,2} = K^{-1} M L^{-1}$, and 
$$A = BP_{1,2}C = B K^{-1} M L^{-1} C.$$
This shows that $A$ is in the subgroup of $GL(n+1, \Z_2)$ generated by $\displaystyle \Hol \left( \moplus_n \Z_2 \right)$ and $M$. 
Thus $\displaystyle \Hol \left( \moplus_n \Z_2 \right)$ is maximal in $GL(n+1, \Z_2)$.  

\end{proof}

\begin{remark} The above corollary does not hold for $R \not\cong \Z_2$ because if $x\in R$ such that $x\ne 0$ or $1$, then all 
entries except the first one of the first column of any product of matrices in $\displaystyle \Hol \left( \moplus_n \Z_2 
\right)$ and the diagonal matrix with diagonal entries $x,$ $1,$ \ldots $1,$ are 0. 
\end{remark}

Let $p$ be prime.
For a group $G$, let $\Syl_p(G)$ denote the $p$-Sylow subgroup of $G$.

\begin{lemma} $\displaystyle \Syl_p \left(\Hol \left( \moplus_n \Z_p \right)\right)$ is isomorphic to 
$\Syl_p \left( GL(n+1, \Z_p)\right) $.
\end{lemma}
\begin{proof}
It is well-known, and is easy to check by computing the order, that the upper-triangular matrices with 1's on the diagonal form a 
$p$-Sylow subgroup of $GL(n+1, \Z_p)$. Namely, the order of $GL(n+1, \Z_p)$ is
$$\displaystyle (p^{n+1}-1)(p^{n+1}-p) \ldots (p^{n+1}-p^n) = \prod_{i=1}^{n}p^i \prod_{i=1}^{n+1}(p^i-1),$$
and the order of the subgroup consisting of upper-triangular matrices with 1's on the diagonal is
$\displaystyle \prod_{i=1}^{n}p^i$.

Since $\Syl_p \left(GL(n+1, \Z_p)\right)$ is contained in the subgroup of $GL(n+1, \Z_p)$ isomorphic to
$\displaystyle \Hol \left( \moplus_n \Z_p \right)$,
$\displaystyle \Syl_p \left(\Hol \left( \moplus_n \Z_p \right)\right) \cong \ \Syl_p \left(GL(n+1, \Z_p)\right)$.
\end{proof}

\begin{corollary}
The inclusion of $\displaystyle \Hol \left( \moplus_n \Z_p \right)$ into $GL(n+1, \Z_p)$ induces an epimorphism in mod $p$ homology and 
a monomorphism in mod $p$ cohomology. 
\end{corollary}

Let $\Syl_p(GL(n, \Z_p))_i$ denote the subgroup of the $\Syl_p(GL(n, \Z_p))$ consisting of all the upper-triangular
matrices with 1's on the diagonal and $a_{kl}=0$ if $k\ne l, k>i$ (and anything to the right from the diagonal in
lines 1 through $i$).

\begin{lemma}
If $i<j$, $\Syl_p (GL(n, \Z_p))_i$ is normal in $\Syl_p(GL(n, \Z_p))_j$.
\end{lemma}

\begin{proof}
Let $\left[ \begin{array}{c|c} A & B \\ \hline 0 & I \end{array} \right] $ where $A$ is an $i\times i$ upper-triangular matrix
with 1's on the diagonal be any element of $\Syl_p(GL(n, \Z_p))_i$, and let
$\left[ \begin{array}{c|c} C & D \\ \hline 0 & E \end{array} \right] $ where $C$ and
$E$ are $i\times i$ and $(n-i)\times(n-i)$ respectively upper-triangular matrices with 1's on the diagonal be any element of
$_p GL(n, \Z_p)$. Then the inverse of $\left[ \begin{array}{c|c} C & D \\ \hline 0 & E \end{array}\right] $ is
$\left[ \begin{array}{c|c} C^{-1} & -C^{-1}DE^{-1} \\ \hline 0 & E^{-1} \end{array} \right] $, and
$ \left[ \begin{array}{c|c} C & D \\ \hline 0 & E \end{array} \right]
\left[ \begin{array}{c|c} A & B \\ \hline 0 & I \end{array} \right]
\left[ \begin{array}{c|c} C^{-1} & -C^{-1}DE^{-1} \\ \hline 0 & E^{-1} \end{array} \right] = $ \\
$ = \left[ \begin{array}{c|c} CAC^{-1} & -CAC^{-1}DE^{-1} +CBE^{-1} + DE^{-1} \\ \hline 0 & I \end{array} \right] $
where $CAC^{-1}$ is an upper-triangular matrix with 1's on the diagonal because $C$ and $A$ are such.
\end{proof}

However, I had no luck with calculating the mod $p$ homology of \\ $\Syl_p(GL(n, \Z_p))$ by running the LHS spectral sequence for
$$1 \rightarrow \ \Syl_p(GL(n, \Z_p))_{i-1} \rightarrow \ \Syl_p(GL(n, \Z_p))_i \rightarrow \moplus_{n-i}\Z_p \rightarrow 1$$
for each $i$ from 2 to $n-1$ $\left( E^2_{qr} \cong H_q\left(\moplus_{n-i}\Z_p ; H_r( \Syl_p(GL(n, \Z_p))_{i-1} ; \F_p )\right) 
\right)$.

If $r\ge 2$, 
$\displaystyle \Syl_p\left( \Hol \left( \moplus_n \Z_{p^r} \right)\right) \not\cong \Syl_p \left(GL(n+1, \Z_{p^r})\right)$ 
because
$$|GL(n, \Z_{p^r})| = (p^{rn}-p^{(r-1)n})(p^{rn}-p^{(r-1)n+1}) \ldots (p^{rn}-p^{(r-1)n+n-1}),$$ thus 
$$| _p GL(n, \Z_{p^r})| = p^{rn^2-\frac{n(n+1)}{2}},$$ and 
$$\displaystyle \left| _p \Hol \left( \moplus_n \Z_{p^r} \right) \right| =
\left| _p \moplus_n \Z_{p^r} \right| \cdot \left| _p GL(n, \Z_{p^r}) \right| =
p^{rn+rn^2-\frac{n(n+1)}{2}} \ne \left| _p GL(n+1, \Z_{p^r}) \right|.$$

\begin{remark}
The subgroup of $GL(n, \Z_{p^r})$ consisting of matrices of the form
$$\renewcommand{\baselinestretch}{1}\normalsize \left[
\begin{array}{ccccc}
a_1 &        &   &   &     \\
    & .      &   & * &     \\
    &        & . &   &     \\
    & b_{ij} &   & . &     \\
    &        &   &   & a_n
\end{array}
\right]$$
where $a_i \equiv 1 (\mbox{mod } p)$ and $b_{ij} \equiv 0 (\mbox{mod } p)$ is a $p-$Sylow subgroup of $GL(n, \Z_{p^r})$.
\end{remark}

\section{Lyndon-Hochschild-Serre Spectral Sequence for $\displaystyle H^* \left( \Hol\left( \moplus_n \Z_{p^r} \right) \right) $}

\label{wilson_section} 

\begin{lemma} \label{productFp} (Wilson's theorem)
$\displaystyle \prod_{c\in\F_p^*}c = -1$.
\end{lemma}
\begin{proof}
$\F_p^*$ is a cyclic group of order $p-1$. Let $s$ be a generator of $\F_p^*$. Then for $p$ odd we have
$$\prod_{c\in\F_p^*} c = \prod_{i=0}^{p-2}s^i = s^{\displaystyle \left( \sum_{i=0}^{p-2}i\right)} =
s^{\displaystyle \frac{(p-1)(p-2)}{2}} = \left(s^{\displaystyle \frac{p-1}{2}}\right)^{p-2} = (-1)^{p-2} = -1.$$
\end{proof}

\begin{theorem} \label{not-collapse-theorem}
The LHS spectral sequence for the extension
$$1 \rightarrow \moplus_n \Z_{p^r} \rightarrow \Hol\left( \moplus_n \Z_{p^r} \right) \rightarrow \Aut\left( \moplus_n \Z_{p^r} \right)
\rightarrow 1$$
does not collapse at the $E_2$ stage if $p=2$ and $r\ge 3$, or $p$ is odd and $r\ge 2$.
\end{theorem}

\begin{proof}
Consider the map of extensions
$$\begin{CD}
1  @>>> \displaystyle \moplus_n \Z_{p^r} @>>> (\Hol(\Z_{p^r}))^n @>>> (\Aut(\Z_{p^r}))^n @>>> 1  \\           
@.      @|                                    @VVV                    @VVV                    @. \\   
1  @>>> \displaystyle \moplus_n \Z_{p^r} @>>> \displaystyle \Hol\left( \moplus_n \Z_{p^r} \right) @>>>
\displaystyle \Aut\left( \moplus_n \Z_{p^r} \right) @>>> 1
\end{CD}$$

The LHS spectral sequence for $\displaystyle H^* \left(\Hol\left( \moplus_n \Z_{p^r} \right); \F_p \right)$ whose $E_2$ stage is given
by
$$E_2^{pq} = H^p \left(GL(n,\Z_{p^r}); H^q \left(\moplus_n \Z_{p^r} ; \F_p \right) \right)$$
maps into that for
$H^* \left( \left(\Hol\left(\Z_{p^r}\right)\right)^n; \F_p\right)$. (The LHS spectral sequence for \\
$H^* \left(\Hol\left(\Z_{p^r}\right) ; \F_p \right)$ is shown in figure \ref{Lyndon_Hol_Z2r} for $p=2$, and in
figure \ref{Lyndon_Hol_Zpr_inv} for $p$ odd.) Recall that
$$\displaystyle H^*\left( \moplus_n \Z_{p^r} ; \F_p \right) \cong \Lambda [u_1, \ldots, u_n] \otimes \F_p[v_1, \ldots, v_n].$$
The vertical axis of our spectral sequence is given by
$$\displaystyle E^{0q}_2 \cong H^0 \left(GL(n,\Z_{p^r}); H^q \left(\moplus_n \Z_{p^r} ; \F_p \right)\right) \cong
\left( H^q\left(\moplus_n\Z_{p^r} ; \F_p \right)\right)^{GL(n,\Z_{p^r})},$$
the ring of invariants in $\displaystyle H^q \left( \moplus_n \Z_{p^r} ; \F_p \right)$ under the action of $GL(n,\Z_{p^r})$.
This contains an analogue of Dickson algebra $\F_p[d_1, \ldots , d_n]$ (defined in \cite{D}; see also \cite{AM}).
The action of $GL(n,\Z_{p^r})$ on $\displaystyle \F_p[v_1, \ldots, v_n] \subset H^q \left( \moplus_n \Z_{p^r} ; \F_p \right)$ is given
by $M(x)=M^Tx$ for $M\in GL(n,\Z_{p^r})$ and
$x=\left[ \begin{array}{c}x_1 \\ \vdots \\ x_n \end{array} \right] \in (\F_p)^n$,
where we identify $v_i$ with the vector with $i^{th}$ coordinate equal to 1, and other coordinates 0.
Every invertible $n\times n$ matrix permutes the $p^n$ elements of $(\F_p)^n$, thus it permutes the set of linear 
combinations of the $v_i$'s, $\{ c_1v_1 + \ldots c_nv_n \ | \ c_i\in \F_p\}$.
It follows that all the coefficients of the polynomial 
$\displaystyle \prod_{c_1,\ldots,c_n\in\F_p} (t+c_1v_{1}+\ldots+c_nv_{n})$ are invariant under the action of $GL(n,\Z_{p^r})$.
Consider the coefficient of $t$. It is a homogeneous polynomial in $v_i$'s of degree $p^n-1$. In fact, it is the last generator 
$d_n$ of the Dickson algebra because the degrees of all other generators are divisible by $p$. 
It contains all the terms of the form 
\vspace{-.2cm}
$$\displaystyle (-1)^n v_{\sigma(1)}^{p^n-p^{n-1}} v_{\sigma(2)}^{p^{n-1}-p^{n-2}} \ldots v_{\sigma(n-1)}^{p^2-p} v_{\sigma(n)}^{p-1}$$
\vspace{-1cm}

\noindent where $\sigma$ is an element of the symmetric group $S_n$, because there is only one way to get this term:
exactly $\frac{1}{p}$th of the multiples in the product $\displaystyle \prod_{c_1,\ldots,c_n\in\F_p} (t+c_1v_{1}+\ldots+c_nv_{n})$
(i.e. $p^{n-1}$ of them) does not contain $v_{\sigma(1)}$ (when $c_{\sigma(1)}=0$),
thus $p^n-p^{n-1}$ of the multiples contain $v_{\sigma(1)}$;
exactly $p^{n-1}-p^{n-2}$ of the remaining multiples contain $v_{\sigma(2)}$, and so on.

\begin{note} The term we actually get this way is
\vspace{-.2cm}
$$\displaystyle \left( \prod_{c\in\F_p^*}c^{p^{n-1}} \right) v_{\sigma(1)}^{p^n-p^{n-1}} \ldots
\left( \prod_{c\in\F_p^*} c \right) v_{\sigma(n)}^{p-1}$$
\vspace{-.5cm}

\noindent but by lemma \ref{productFp} each product $\displaystyle \prod_{c\in\F_p^*}c^{p^i} \equiv -1$ (mod $p$).
\end{note}

In the spectral sequence for $ H^* \left(\left(\Hol\left(\Z_{p^r}\right)\right)^n\right)$ we have
$$\begin{array}{l}
\displaystyle d_2 \left( \sum_{\sigma\in \Sigma_n} v_{\sigma(1)}^{p^n-p^{n-1}} v_{\sigma(2)}^{p^{n-1}-p^{n-2}} \ldots
v_{\sigma(n-1)}^{p^2-p} v_{\sigma(n)}^{p-1} \right) \\
= \displaystyle \sum_{\sigma\in \Sigma_n} \left(
d_2 \left(v_{\sigma(1)}^{p^n-p^{n-1}} \right) v_{\sigma(2)}^{p^{n-1}-p^{n-2}}\ldots v_{\sigma(n)}^{p-1} +
\ldots + v_{\sigma(1)}^{p^n-p^{n-1}} \ldots v_{\sigma(n-1)}^{p^2-p} d_2 \left( v_{\sigma(n)}^{p-1} \right)
\right) \\
= \displaystyle \sum_{\sigma\in \Sigma_n} \left(
v_{\sigma(1)}^{p^n-p^{n-1}} \ldots v_{\sigma(n-1)}^{p^2-p} d_2 \left( v_{\sigma(n)}^{p-1} \right) \right) \\
= \left\{ \begin{array}{l@{}l@{}l}
\displaystyle \sum_{\sigma\in \Sigma_n} \left(
v_{\sigma(1)}^{p^n-p^{n-1}} \ldots v_{\sigma(n-1)}^{p^2-p} m v_{\sigma(n)}^{p-2}u_{\sigma(n)}z_{\sigma(n)}
\right) & \mbox{ if } & p\mbox{ is odd, or } p=2 \mbox{ and } r>3\\
\displaystyle \sum_{\sigma\in \Sigma_n} \left(
v_{\sigma(1)}^{2^n-2^{n-1}} \ldots v_{\sigma(n-1)}^{2^2-2} u_{\sigma(n)}x_{\sigma(n)}^2 \right) & \mbox{ if } & 
p=2 \mbox{ and } r=3
\end{array}\right. \\
\ne 0
\end{array}$$
since all the terms in the sum are linearly independent.
The other terms in the coefficient of $t$ are of the form $v_1^{e_1} v_2^{e_2} \ldots v_n^{e_n}$ where the set of exponents
$\{ e_1, \ldots, e_n\}$ is not a permutation of $\{ p^n-p^{n-1}, p^{n-1}-p^{n-2}, \ldots, p^2-p, p-1\}$, so
$d_2(v_1^{e_1} v_2^{e_2} \ldots v_n^{e_n})$ in the spectral sequence for
$ H^* \left(\left(\Hol\left(\Z_{p^r}\right)\right)^n; \F_p \right)$ is either 0 or linearly independent from
$\{ d_2(v_{\sigma(1)}^{p^n-p^{n-1}} v_{\sigma(2)}^{p^{n-1}-p^{n-2}} \ldots v_{\sigma(n-1)}^{p^2-p} v_{\sigma(n)}^{p-1}) | 
\sigma\in\Sigma_n \}$. 
Therefore the image of the coefficient of $t$ under $d_2$ in the spectral sequence for
$\displaystyle H^* \left( \Hol\left( \moplus_n \Z_{p^r} \right) ; \F_p \right)$ is non-zero.
\end{proof}

\begin{corollary} \label{nonzero-cohomology-corollary} 
Let $n\ge 2$. Then
$$\displaystyle H^2\left( GL(n,\Z_{p^r}); H^{2p^n-3}\left( \moplus_n \Z_{p^r} ; \F_p \right)\right) \ne 0$$
if $p=2$ and $r\ge 3$, or if $p$ is odd and $r\ge 2$, where the action of $GL(n,\Z_{p^r})$ on 
$H^{2p^n-3}\left( \moplus_n \Z_{p^r} ; \F_p \right)$ is induced by the natural action of $GL(n,\Z_{p^r})$ on 
$\displaystyle \Z_{p^r}$.
\end{corollary}

\section{Congruence subgroups}

\label{congruence_subgroups}

We have seen in example \ref{characteristic-example-Z} that for $k>l$ there is a reduction homomorphism 
$\Hol(\Z_{p^k})\rightarrow\Hol(\Z_{p^l})$. Similarly, there is a reduction homomorphism
$\displaystyle \Hol\left(\moplus_n \Z_{p^k}\right) \rightarrow \Hol\left(\moplus_n \Z_{p^l}\right)$.
The goal of this section is to compute the cohomology of the kernel of
\begin{equation}\label{holomorph-reduction-map}
\displaystyle \Hol\left(\moplus_n \Z_{p^k}\right) \rightarrow \Hol\left(\moplus_n \Z_p\right).
\end{equation}
for odd $p$.

First recall a few definitions and a theorem from \cite{P} (everything we need from \cite{P} is also contained in \cite{BP}). Let $p$ be 
an odd prime. 

For a finite $p$-group $G$, let $\Omega_1(G)$ be the subgroup generated by $g\in G$ such that $g^p=1$.
The group $G$ is called $p$-central if $\Omega_1(G)$ is central in $G$.

A tower of groups is a sequence of $p$-central groups $\{G_1, G_2, \ldots, G_n\}$ with surjective homomorphisms
$\pi_i:G_i\rightarrow G_{i-1}$ with $\Ker(\pi_i)\cong\Omega_1(G_i)$.

The central extension 
\begin{equation}\label{central-extension}
1 \rightarrow \Omega_1(G_i) \rightarrow G_i \rightarrow G_{i-1} \rightarrow 1
\end{equation}
pulls back to
$$1 \rightarrow \Omega_1(G_i) \rightarrow E_i \rightarrow \Omega_1(G_{i-1}) \rightarrow 1.$$
The $p$-power map $\phi:\Omega_1(G_{i-1}) \rightarrow \Omega_1(G_i)$ is defined by $\phi(x)=\bar{x}^p\in \Omega_1(G_i)$ where
$\bar{x}\in E_i$ is a preimage of $x\in \Omega_1(G_{i-1})$.

A tower of groups is called uniform if all the $p$-power maps of the tower are isomorphisms.

\begin{theorem}\label{jonathan-theorem} (W. Browder, J. Pakianathan, T.Weigel) Let $\{G_1, G_2, \ldots, G_N\}$ be a uniform tower 
with dim$(G_1)=n.$
Then for $1 \le k < N$ we have
$$H^*(G_k;\F_p)=\Lambda(x_1,\ldots,x_n)\otimes\F_p[s_1,\ldots,s_n]$$
where $|x_i|=1$ (and for $k>1$, $x_i$ is induced from $G_{k-1}$) and $|s_i|=2$ (and for $k>1$, $s_i$ ``comes'' from $\Omega(G_k))$, and
$H^*(G_N;\F_p)$ is given by the same formula if and only if the tower extends to a uniform tower $\{G_1,\ldots,G_N,G_{N+1}\}$.
\end{theorem}

It follows from the LHS spectral sequence for the extension \eqref{central-extension} (see \cite{P} for details) that the map 
$H^*(G_{k-1}) \rightarrow H^*(G_k)$ is the identity map on the exterior part and trivial on the polynomial part. Let 
$G= \lim_{\leftarrow} G_k$. Then we have 

\begin{corollary}
$H_{cont}^*(G;\F_p) = \Lambda(x_1,\ldots,x_n),$ \  $|x_i|=1.$
\end{corollary} 

We will show that the kernels of the reduction map \eqref{holomorph-reduction-map} fit into an infinite uniform tower, and hence theorem
\ref{jonathan-theorem} is applicable.

\ As in \cite{P}, let $\Gamma_{n,k}$ be the kernel of the reduction map \mbox{$GL(n,\Z_{p^{k+1}}) \rightarrow GL(n,\Z_p)$}. \
Let \ $\Gamma_{n,k}^{\scriptHol}$ \ be \ the \ kernel \ of \ the \ reduction \ map \
\mbox{$\displaystyle \Hol\left(\moplus_n \Z_{p^{k+1}}\right) \rightarrow \Hol\left(\moplus_n \Z_{p}\right)$}.

We then have the following maps of extensions:
\begin{equation}\label{reduction-map-gamma}\begin{CD}
1  @>>> \Gamma_{n,k}   @>>> GL(n,\Z_{p^{k+1}}) @>>> GL(n,\Z_p) @>>> 1  \\
@.      @VVV                @VVV                    @|              @. \\ 
1  @>>> \Gamma_{n,k-1} @>>> GL(n,\Z_{p^{k}})   @>>> GL(n,\Z_p) @>>> 1
\end{CD}\end{equation}

and similarly

\begin{equation}\label{reduction-map-gamma-hol}\begin{CD}
1 @>>> \Gamma_{n,k}^{\scriptHol}   @>>> \displaystyle \Hol\left(\moplus_n\Z_{p^{k+1}}\right) @>>>
\displaystyle \Hol\left(\moplus_n\Z_p\right)  @>>> 1 \\
@.     @VVV                             @VVV                                                       @|       @. \\ 
1 @>>> \Gamma_{n,k-1}^{\scriptHol} @>>> \displaystyle \Hol\left(\moplus_n\Z_{p^{k}}\right)   @>>> 
\displaystyle \Hol\left(\moplus_n\Z_p\right)  @>>> 1
\end{CD}\end{equation}

Also, $\Gamma_{n,k}^{\scriptHol}$ and $\Gamma_{n,k}$ fit into a split extension as follows:

\begin{equation}\label{split-ext-gamma-hol-gamma}\begin{CD}
   @.  1    @. 1    @. 1    @.     \\
@.     @VVV    @VVV    @VVV    @.  \\
1 @>>> \displaystyle \moplus_n p\Z_{p^{k+1}} @>>> \Gamma_{n,k}^{\scriptHol} @>{\curvearrowleft}>> \Gamma_{n,k} @>>> 1 \\
@.     @VVV    @VVV    @VVV    @.   \\
1 @>>> \displaystyle \moplus_n \Z_{p^{k+1}}  @>>> \displaystyle \Hol\left(\moplus_n \Z_{p^{k+1}}\right) @>{\curvearrowleft}>>
GL(n,\Z_{p^{k+1}}) @>>> 1 \\
@.     @VVV    @VVV    @VVV    @.   \\
1 @>>> \displaystyle \moplus_n \Z_p          @>>> \displaystyle \Hol\left(\moplus_n \Z_p\right)         @>{\curvearrowleft}>> 
GL(n,\Z_p)         @>>> 1 \\
@.     @VVV    @VVV    @VVV    @.   \\
   @.  1    @. 1    @. 1    @.     
\end{CD}\end{equation}

Combining the first row of \eqref{split-ext-gamma-hol-gamma} and the reduction maps in \eqref{reduction-map-gamma} and
\eqref{reduction-map-gamma-hol}, we get

\begin{equation}\label{main-gamma-diagram}\begin{CD}
   @.  1    @. 1    @. 1    @.     \\
@.     @VVV    @VVV    @VVV    @.   \\
1 @>>> \displaystyle\moplus_n p^k\Z_{p^{k+1}} @>>> \Ker_{n,k}^{\scriptHol}     @>{\curvearrowleft}>> \Ker_{n,k}   @>>> 1\\
@.     @VVV    @VVV    @VVV    @.   \\
1 @>>> \displaystyle \moplus_n p\Z_{p^{k+1}}  @>>> \Gamma_{n,k}^{\scriptHol}   @>{\curvearrowleft}>> \Gamma_{n,k} @>>> 1\\
@.     @VVV    @VVV    @VVV    @.   \\
1  @>>> \displaystyle \moplus_n p \Z_{p^k  }  @>>> \Gamma_{n,k-1}^{\scriptHol} @>{\curvearrowleft}>> \Gamma_{n,k-1} @>>> 1\\
@.     @VVV    @VVV    @VVV    @.   \\
   @.  1    @. 1    @. 1    @.     
\end{CD}\end{equation}

We think of $\Gamma_{n,k}$ as the group consisting of matrices of the form $1+pA$ (mod $p^{k+1}$),
$\Gamma_{n,k-1}$ consisting of matrices of the form $1+pA$ (mod $p^k$),
and $\Ker_{n,k}$ then consists of matrices of the form $1+p^kA$ (mod $p^{k+1}$).

Let $\pi_{n,k}$ be the surjection in the middle column of \eqref{main-gamma-diagram}.

\begin{theorem}
For each $n\in\N$, the groups $\Gamma_{n,1}^{\scriptHol},$ $\Gamma_{n,2}^{\scriptHol},$ $ \ldots$ together with the surjections
$\pi_{n,k}$ form an infinite
uniform tower.
\end{theorem}
The following four lemmas give a proof of this theorem.

\begin{lemma}
$\Gamma_{n,k}^{\scriptHol}$ is a $p$-group.
\end{lemma}
\begin{proof}
$\displaystyle \left| \Gamma_{n,k}^{\scriptHol} \right| = \left| \moplus_n p \Z_{p^{k+1}} \right| \cdot \left| \Gamma_{n,k} \right| =
(p^k)^n \frac{|GL(n,\Z_{p^{k+1}})|}{|GL(n,\Z_p)|} = $ \\
$\displaystyle = p^{kn}
\frac{(p^{(k+1)n}-p^{kn})(p^{(k+1)n}-p^{kn+1})\ldots(p^{(k+1)n}-p^{kn+n-1})}{(p^n-1)(p^n-p)\ldots(p^n-p^{n-1})} =
p^{k(n^2+n)}.$
\end{proof}

\begin{lemma}
$\Ker \left( \pi_{n,k}:\Gamma_{n,k}^{\scriptHol} \rightarrow \Gamma_{n,k-1}^{\scriptHol} \right) \cong
\Omega_1 \left( \Gamma_{n,k}^{\scriptHol} \right)$.
\end{lemma}
\begin{proof}
Let $g=(1+pA,px)\in \Gamma_{n,k}^{\scriptHol}$. Then it can be easily shown by induction that
$g^i = \left(1+ipA+\frac{i(i-1)}{2}p^2A^2+p^3(\ldots), ipx-\frac{i(i-1)}{2}p^2Ax+p^3(\ldots)\right)$.
For $i=p$ (recall $p$ is odd) we have $g^p=(1+p^2A+p^3(\ldots),p^2x+p^3(\ldots))$, so $g^p=1$ \ $\Leftrightarrow$ \
$p^2A=0$ and $p^2x=0$ (mod $p^{k+1}$) $\Leftrightarrow$ $A=p^{k-1}B$, $x=p^{k-1}y$ $\Leftrightarrow$ $g=(1+p^kB, p^ky)$
$\Leftrightarrow$ $g\in \Ker(\pi_{n,k})$.
\end{proof}

\begin{lemma}
$\Omega_1 \left( \Gamma_{n,k}^{\scriptHol} \right)$ is central in $\Gamma_{n,k}^{\scriptHol}$ (thus $\Gamma_{n,k}^{\scriptHol}$ is
$p$-central).
\end{lemma}
\begin{proof}
Every element of $\Omega_1 \left( \Gamma_{n,k}^{\scriptHol} \right)$ commutes with every element of $\Gamma_{n,k}^{\scriptHol}$: \\
$(1+p^kB, p^ky)(1+pA, px) \equiv (1+p^kB+pA, p^ky+px) \equiv $\\
$\equiv (1+pA, px)(1+p^kB, p^ky)$ (mod $p^{k+1}$).
\end{proof}

Thus we have an infinite tower $\{\Gamma_{n,1}^{\scriptHol}, \Gamma_{n,2}^{\scriptHol}, \ldots \}$
with $\Gamma_{n,1}^{\scriptHol} \cong \Omega_1 \left( \Gamma_{n,1}^{\scriptHol} \right)$ an elementary abelian
$p$-group of rank $n^2+n$. It is useful to say more about $\Gamma_{n,1}^{\scriptHol}$. Its elements have the form  $(1+pA,px)$
where $1+pA\in GL(n,\Z_{p^2})$ and $\displaystyle px\in \moplus_n \Z_{p^2}$. Let $H$ be the group whose elements are pairs
$(A,x)$ with $A\in Mat_{n\times n}(\Z_p)$ and $\displaystyle x\in \moplus_n \Z_p$, and the multiplication is given by
$(A,x)(B,y)=(A+B,x+y)$. Thus $H$ is an elementary abelian group of rank $n^2+n$. It is easy to check that the map
\begin{equation}\label{gamma_n1-iso}
\Gamma_{n,1}^{\scriptHol} \buildrel h \over \rightarrow H
\end{equation}
defined by $h((1+pA,px)) = (A,x)$ (mod $p$) is an isomorphism.

Now consider the $p$-power map $\phi$ for the central extension
$$1 \rightarrow \Omega_1\left( \Gamma_{n,k}^{\scriptHol} \right)
\rightarrow E_k \rightarrow \Omega_1\left( \Gamma_{n,k-1}^{\scriptHol} \right) \rightarrow 1$$
$\phi(g)=\bar{g}^p$.

\begin{lemma}
The $p$-power map $\phi:\Omega_1\left( \Gamma_{n,k-1}^{\scriptHol} \right) \rightarrow \Omega_1\left( \Gamma_{n,k}^{\scriptHol} \right)$
is an isomorphism (so the tower is uniform).
\end{lemma}
\begin{proof}
Let $g\in \Omega_1\left( \Gamma_{n,k-1}^{\scriptHol} \right)$, then $g=(1+p^{k-1}A,p^{k-1}x)$ (mod $p^k$). Its lift to
$\Gamma_{n,k}^{\scriptHol}$ is $\bar{g}=(1+p^{k-1}A,p^{k-1}x)$ (mod $p^{k+1}$). Then
$\bar{g}^p=(1+p^kA,p^kx)$ (mod $p^{k+1}$). Since $p^{k-1}\Z_{p^k}\cong p^k\Z_{p^{k+1}}$, we have an isomorphism.
\end{proof}

\begin{corollary}
If $p$ is odd,
$$H^* \left( \Gamma_{n,k}^{\scriptHol}; \F_p \right) \cong \Lambda(x_{ij}, x_l) \otimes \F_p[s_{ij}, s_l], 
\hspace{1cm} 1 \le i,j,l \le n $$
where $|x_{ij}|=|x_l|=1$ and $|s_{ij}|=|s_l|=2$. For $k>1$, $x_{ij}$'s and $x_l$'s are induced from 
$\Gamma_{n,k-1}^{\scriptHol}$, and $s_{ij}$'s and $s_l$'s ``come'' from $\Omega_1\left(\Gamma_{n,k}^{\scriptHol}\right)$.  
(The reason for using two indices for some generators and one index for others will be explained later.)
\end{corollary}

Recall also (see \cite{P}) that for a central extension
$$1 \rightarrow V \rightarrow G \rightarrow W \rightarrow 1$$
where $V$ and $W$ are elementary abelian $p$-groups, and the $p$-power map $\phi:W\rightarrow V$ is an isomorphism, one defines a
bracket on $W$ by
$$[x,y]=\phi^{-1}(\bar{x}\bar{y}\bar{x}^{-1}\bar{y}^{-1})$$
where $\bar{x}, \bar{y}\in G$ are preimages of $x, y\in W$.

Let $\{G_1,G_2\}$ be a uniform tower, i.e. $1\rightarrow \Omega_1(G_2) \rightarrow G_2 \rightarrow G_1 \rightarrow 1$ is a central
extension with $\phi:G_1 \rightarrow \Omega_1(G_2)$ an isomorphism of elementary abelian $p$-groups. Then we have a bracket on $G_1$.
J. Pakianathan proves that the tower $\{G_1,G_2\}$ extends to a uniform tower of length 3 if and only if this bracket is a
Lie algebra bracket, i.e. it satisfies the Jacobi identity.

Now let $\{G_1, G_2, \ldots, G_N\}$ be a uniform tower, and let $c_{ij}^k$'s be the structure constants of the Lie algebra associated to
$G_2$ as above and denoted by $Log(G_2)$. Then for $p\ge 5$ (the restriction $p\ge 5$ is needed because this 
approach uses Weigel's functors {\bf exp} and {\bf log} for powerful $p$-central Lie algebras and groups, see 
\cite{We}), he obtains the following formulas for the first Bockstein homomorphisms $\beta: H^*(G_k;\F_p)\rightarrow 
H^{*+1}(G_k;\F_p)$ for $1 < k < N$:
\begin{equation}\label{jonathan-bocksteins}
\begin{array}{rcl}
\beta(x_t) & = & \displaystyle -\sum_{i<j} c_{ij}^t x_i x_j \\
\beta(s_t) & = & \displaystyle \sum_{i,j=1}^n c_{ij}^t s_i x_j
\end{array}
\end{equation}

We will use these formulas to calculate the first Bocksteins for $\Gamma_{n,k}^{\scriptHol}$.
First we derive an explicit formula for the bracket on $\Gamma_{n,1}^{\scriptHol}$.
Recall that an element of $\Gamma_{n,1}^{\scriptHol}$ has the form $(1+pA,px)$ where $1+pA\in GL(n,\Z_{p^2})$ and
$\displaystyle px\in \moplus_n \Z_{p^2}$. \\
$$[(1+pA,px),(1+pB,py)] = \phi^{-1}((1+pA,px)(1+pB,py)(1+pA,px)^{-1}(1+pB,py)^{-1})$$
where we denote lifts of elements by the same symbols, only now the matrices and vectors are considered mod $p^3$ rather than mod $p^2$.\\
$(1+pA,px)(1+pB,py)(1+pA,px)^{-1}(1+pB,py)^{-1} \equiv$ \\
$\equiv (1+pA,px)(1+pB,py)((1+pA)^{-1},(1+pA)(-px))((1+pB)^{-1},(1+pB)(-py)) \equiv$ \\
$\equiv (1+pA,px)(1+pB,py)(1-pA+p^2A^2,-px-p^2Ax)(1-pB+p^2B^2,-py-p^2By) \equiv $ \\
$\equiv (1+p^2(AB-BA), p^2(Ay-Bx))$ (mod $p^3$).

The inverse of the $p$-power isomorphism ``reduces'' the exponent of $p$ by 1, thus we get
$$[(1+pA,px),(1+pB,py)] = (1+p(AB-BA), p(Ay-Bx)) \ \ (\mbox{mod } p^2)$$

For simplicity, we use the isomorphism \eqref{gamma_n1-iso} to rewrite the bracket in terms of elements of the form $(A,x)$ where 
$A\in Mat_{n\times n}(\Z_p)$ and $\displaystyle x\in \moplus_n \Z_p$. 

\begin{theorem} The Lie algebra bracket in $\Gamma_{n,k}^{Hol}$ is given by 
\begin{equation} \label{bracket-formula}
[(A,x),(B,y)] = (AB-BA,Ay-Bx).
\end{equation}
\end{theorem}

Let $E_{ij}$ be the $n \times n$-matrix with the $ij$-th entry 1 and the other entries 0. Let $E_l$ be the $n$-vector with $l$-th entry
1 and the other entries 0. Let $e_{ij}=(E_{ij},0)$ and $e_l=(0,E_l)$. Then $\{e_{11}, \ldots, e_{nn}, e_1, \ldots, e_n\}$ is a basis
for
$H$ isomorphic to $\Gamma_{n,1}^{\scriptHol}$ by \eqref{gamma_n1-iso}.
In this basis, the structure constants of the associated Lie algebra are as follows.

$$\begin{array}{rll}
c_{ij \ lm}^{tu} & = & \delta_{it} \delta_{um} \delta_{jl} - \delta_{tl} \delta_{uj} \delta_{im} \\
c_{ij \ lm}^{t}  & = & 0 \\
c_{ij \ l}^{tu}  & = & 0 \\
\end{array}$$

$$\begin{array}{rll}
c_{ij \ l}^{t}   & = & \delta_{ti} \delta_{jl} \\
c_{i \ lm}^{tu}  & = & 0 \\
c_{i \ lm}^{t}   & = & - \delta_{tl} \delta_{im} \\
c_{i \ l}^{tu}   & = & 0 \\
c_{i \ l}^{t}    & = & 0
\end{array}$$

Using \eqref{jonathan-bocksteins}, we have

$$\begin{array}{lll}
\beta(x_{tu}) & = & \displaystyle -\sum_{ij<lm} c_{ij\ lm}^{tu} x_{ij} x_{lm}
                    -\sum_{ij,l} c_{ij\ l}^{tu} x_{ij} x_{l} -\sum_{i<l} c_{i\ l}^{tu} x_i x_l \\
              & = & \displaystyle -\sum_{ij<lm}^{ } \delta_{it} \delta_{um} \delta_{jl} x_{ij} x_{lm}
                    + \sum_{ij<lm} \delta_{tl} \delta_{uj} \delta_{im} x_{ij} x_{lm} \\
              & = & \displaystyle -\sum_{tl<lu}^{} x_{tl} x_{lu} + \sum_{iu<ti} x_{iu} x_{ti} \\
              & = & - \left\{ \begin{array}{lll}
                        \displaystyle \sum_{l=t}^n x_{tl} x_{lu} & \mbox{if} & t<u \\
                        \displaystyle \sum_{l=t+1}^n x_{tl} x_{lu} & \mbox{if} & t\ge u
 \end{array}\right. + \left\{ \begin{array}{lll}
                        \displaystyle \sum_{i=1}^{t} x_{iu} x_{ti} & \mbox{if} & u<t \\
                        \displaystyle \sum_{i=1}^{t-1} x_{iu} x_{ti} & \mbox{if} & u\ge t
\end{array}\right. \\
\beta(x_t) & = & \displaystyle -\sum_{ij<lm} c_{ij\ lm}^{t}  x_{ij} x_{lm}
                 -\sum_{ij,l} c_{ij\ l}^{t} x_{ij} x_{l} -\sum_{i<l} c_{i\ l}^{t} x_i x_l \\
           & = & \displaystyle -\sum_{ij,l}^{ } \delta_{ti} \delta_{jl} x_{ij} x_{l} \\
           & = & \displaystyle -\sum_{l=1}^{n} x_{tl} x_l \\ 
\beta(s_{tu}) & = & \displaystyle \sum_{ij,lm} c_{ij\ lm}^{tu} s_{ij} x_{lm} + \sum_{ij,l} c_{ij\ l}^{tu} s_{ij} x_{l}
                    + \sum_{i,lm}c_{i\ lm}^{tu} s_{i} x_{lm} + \sum_{i,l} c_{i\ l}^{tu} s_i x_l \\
              & = & \displaystyle \sum_{ij,lm}^{ } \delta_{it} \delta_{um} \delta_{jl} s_{ij} x_{lm}
                    - \sum_{ij,lm} \delta_{tl} \delta_{uj} \delta_{im} s_{ij} x_{lm} \\
              & = & \displaystyle \sum_{tl,lu}^{ } s_{tl} x_{lu} - \sum_{iu,ti} s_{iu} x_{ti} \\
              & = & \displaystyle \sum_{l=1}^n s_{tl} x_{lu} - \sum_{i=1}^n s_{iu} x_{ti}
\end{array}$$

$$\begin{array}{lll}
\beta(s_t) & = & \displaystyle \sum_{ij,lm} c_{ij\ lm}^{t} s_{ij} x_{lm} + \sum_{ij,l} c_{ij\ l}^{t} s_{ij} x_{l}
                    + \sum_{i,lm}c_{i\ lm}^{t} s_{i} x_{lm} + \sum_{i,l} c_{i\ l}^{t} s_i x_l \\
           & = & \displaystyle \sum_{ij,l}^{ } \delta_{ti} \delta_{jl} s_{ij} x_{l}
                 - \sum_{i,lm} \delta_{tl} \delta_{im} s_{i} x_{lm} \\
           & = & \displaystyle \sum_{tl,l}^{ } s_{tl} x_l - \sum_{i,ti} s_i x_{ti} \\
           & = & \displaystyle \sum_{l=1}^n s_{tl} x_l - \sum_{i=1}^n s_i x_{ti}
\end{array}$$

Thus we have 
\begin{theorem}\label{congruence-ring-Bock}
If $p\ge 3$, then
$$H^* \left( \Gamma_{n,k}^{\scriptHol}; \F_p \right) \cong \Lambda(x_{ij}, x_{l}) \otimes \F_p[s_{ij}, s_{l}], 
\hspace{1cm} 1 \le i,j,k, \le n$$
where $|x_{ij}|=|x_l|=1$ and $|s_{ij}|=|s_l|=2$. Also, if $p\ge 5$, the first Bocksteins are given by
$$\begin{array}{rll}
\beta(x_{tu}) & = & \left\{
        \begin{array}{lll}
        \displaystyle - \sum_{i=t}^n x_{ti}x_{iu} + \sum_{i=1}^{t-1} x_{iu}x_{ti} & \mbox{if} & t<u \\
        \displaystyle - \sum_{i=t+1}^n x_{ti}x_{iu} + \sum_{i=1}^{t-1} x_{iu}x_{ti} & \mbox{if} & t=u \\
        \displaystyle - \sum_{i=t+1}^n x_{ti}x_{iu} + \sum_{i=1}^{t} x_{iu}x_{ti} & \mbox{if} & t>u
        \end{array} \right. \\
\beta(x_t)    & = & \displaystyle - \sum_{i=1}^n x_{ti} x_i \\
\beta(s_{tu}) & = & \displaystyle \sum_{i=1}^n (s_{ti}x_{iu}-s_{iu}x_{ti}) \\
\beta(s_t)    & = & \displaystyle \sum_{i=1}^n (s_{ti} x_i - s_i x_{ti})
\end{array}$$
\end{theorem}

\begin{remark} 
We will see in the next section that these formulas hold for $p=2$ and $p=3$ as well. Notice that the formula for $\beta_1(x_{tu})$ 
simplifies a little bit for $p=2$. 
\end{remark}

\section{Minh and Symonds's approach: covers cases $p=2$ and $3$} 

\label{minh_simonds}

In this section we will show that Minh and Symonds's approach described in \cite{MS} applies to 
$\Gamma_{n,k}^{\scriptHol}$ for any prime $p$. 
Recall that if $G$ is a pro-$p$ group, then a closed, normal, finitely generated subgroup $N$ is almost powerfully 
embedded if 
\begin{itemize}
\item $[G,N] \subset N^p $ for $p$ odd 
\item $[G,N] \subset N^2$ and $[N,N] \subset N^4$ for $p=2$
\end{itemize}
where $N^p$ denotes the closure of the subgroup of $N$ generated by $p$-th powers of its elements. 

\begin{lemma} \label{mike-lemma}
\begin{enumerate} 
\item If $p$ is and odd prime, then any element of $\Z_{p^{k+1}}$ of the form $1+p^2a$ is the $p$-th power of some element of 
the form $1+pb$. 
\item Any element of $\Z_{2^{k+1}}$ of the form $1+8a$ is the square of some element of the form $1+4b$.
\item Any element of $\Z_{2^{k+1}}$ of the form $1+16a$ is the fourth power of some element of the form $1+4b$. 
\end{enumerate}
\end{lemma} 

\begin{proof}
The author had quite an ugly proof of this lemma. The idea of the proof given below belongs to Michael Knapp.  
\begin{enumerate}

\item First notice that $(1+pb)^p \equiv 1+p^2a \ (\mbox{mod } p^{k+1})$, and conversely,  
since $x^p \equiv x \ (\mbox{mod } p),$ $x^p\equiv 1+pa \ (\mbox{mod } p^{k+1})$ implies that $x\equiv 1 \ 
(\mbox{mod } p)$, therefore $x\in \Aut(\Z_{p^{k+1}})$. Let $s$ be a generator of $\Aut(\Z_{p^{k+1}})\cong \Z_{(p-1)p^k}.$ 
Then $1+p^2a=s^t$ for some $t$, mod $(p-1)p^k$. By theorem 5-1 of \cite{A}, given $t$, the equation 
$qp\equiv t \ (\mbox{mod } (p-1)p^k)$ has at most $p$ solutions. Therefore the equation 
$(s^q)^p\equiv s^t\ (\mbox{mod } p^{k+1})$ has at most $p$ solutions. Thus for each $a$, there are at 
most $p$ numbers of the form $1+pb$ such that $(1+pb)^p \equiv 1+p^2a \ (\mbox{mod } p^{k+1})$. 
There are $p^k$ numbers of the form $1+pb$ and $p^{k-1}$ numbers of the form $1+p^2a$. 
Hence every equation $(1+pb)^p \equiv 1+p^2a \ (\mbox{mod } p^{k+1})$ has a solution. 

\item As in part 1, $(1+4b)^4 \equiv 1+8a\ (\mbox{mod } 2^{k+1})$, and conversely,
since the square of an even number is even, $x^2\equiv 1+8a \ (\mbox{mod } 2^{k+1})$ implies $x$ is odd, i.e.
$x\in \Aut(\Z_{2^{k+1}}) \cong \Z_2 \oplus \Z_{2^{k-1}} = <-1> \oplus <3>$. In this basis, $1+8a$ is of the form
$(1,3^t)$ (it is easy to see that it can not be of the form $(-1,3^t)$). The equation then can be rewritten as
$(b_1, b_2)^2\equiv (1,3^t)$. Since $(b_1,b_2)^2 = (1,b_2^2) = (1,(3^q)^2),$ the equation becomes
$(3^q)^2\equiv 3^t \ (\mbox{mod } 2^{k+1})$ or $2q \equiv t \ (\mbox{mod } 2^{k-1})$. By 5-1 of \cite{A}, the latter has at
most 2 solutions. Since we are only interested in solutions of the form $1+4b$ (but not $3+4b$), $b_1$ is determined by
$b_2$. Since there are $2^{k-1}$ elements of the form $1+4b$ and $2^{k-2}$ elements of the form $1+8a$, every equation
$(1+4b)^2 \equiv 1+8a \ (\mbox{mod } 2^{k+1})$ has a solution.

\item The proof of this part is the proof of part 2 with minor modifications. 
First, $(1+4b)^4 \equiv 1+16a\ (\mbox{mod } 2^{k+1})$, and conversely,
since the fourth power of an even number is even, $x^4\equiv 1+16a \ (\mbox{mod } 2^{k+1})$ implies $x$ is odd, i.e. 
$x\in \Aut(\Z_{2^{k+1}}) \cong \Z_2 \oplus \Z_{2^{k-1}} = <-1> \oplus <3>$. In this basis, $1+16a$ is of the form 
$(1,3^t)$. The equation then can be rewritten as 
$(b_1, b_2)^4\equiv (1,3^t)$. Since $(b_1,b_2)^4 = (1,b_2^4) = (1,(3^q)^4),$ the equation becomes 
$(3^q)^4\equiv 3^t \ (\mbox{mod } 2^{k+1})$ or $4q \equiv t \ (\mbox{mod } 2^{k-1})$. By 5-1 of \cite{A}, the latter has at 
most 4 solutions. Since we are only interested in solutions of the form $1+4b$ (but not $3+4b$), $b_1$ is determined by 
$b_2$. Since there are $2^{k-1}$ elements of the form $1+4b$ and $2^{k-3}$ elements of the form $1+16a$, every equation 
$(1+4b)^4 \equiv 1+16a \ (\mbox{mod } 2^{k+1})$ has a solution. 

\end{enumerate}
\end{proof}

\begin{theorem} Let $k\ge 2$, and $G=\Gamma_{n,k}^{\scriptHol}=\{ (1+pA,px)(\mbox{mod } p^{k+1}) \}$.
\begin{enumerate}
\item If $p$ is odd, let $N=G$.
\item If $p=2$, let $N=\{ (1+4A,4x)(\mbox{mod } 2^{k+1}) \}$. 
\end{enumerate}
Then $N$ is almost powerful embedded in $G$. 
\end{theorem}

\begin{proof}
\begin{enumerate}
\item We have seen that \\
$[(1+pA,px),(1+pB,py)] \equiv (1+p^2(AB-BA), p^2(Ay-Bx)) \hfill (\mbox{mod } p^3),$ \\
thus $[G,N] \subset \{(1+p^2C,p^2z)(\mbox{mod } p^{k+1}) \}$. Now we will show that any element of the form
$(1+p^2C,p^2z)$ is generated by the $p$-th powers of some elements of the form $(1+pD,pw)$. 
$$(1+p^2C,p^2z) = (1+p^2C,0)(1,p^2z)=\prod_{(i,j)\in I}(1+p^2E_{ij},0)(1,p^2z)$$ 
where $E_{ij}$ is the matrix whose $ij$-th entry is 1 and all the other entries are 0. By part 1 of lemma \ref{mike-lemma}, there exists 
$x\equiv 0(\mbox{mod }p)$ such that $(1+x)^p\equiv 1+p^2(\mbox{mod }p^{k+1})$. Then $(1+p^2E_{ii},0) \equiv ((1+xE_{ii})^p,0)$. If $i\ne 
j$, then $(1+p^2E_{ij},0) \equiv ((1+pE_{ij})^p,0)$. Also, $(1,p^z) \equiv (1,pz)^p$. Thus $[G,N] \subset N^p$. 

\item Let $(1+2A,2x)\in G$ and $1+4B,4x)\in N$. Then \\
$[(1+2A,2x),(1+4B,4y)] \equiv (1+8(AB-BA), 8(Ay-Bx)) \hfill (\mbox{mod } 16),$ \\
thus $[G,N] \subset \{(1+8C,8z)(\mbox{mod } 2^{k+1}) \}$. The proof that any element of the form $(1+8C,8z)$ is generated by the squares 
of elements of the form $(1+4D,4w)$ is analogous to the proof of part 1 but is using \ref{mike-lemma} (2) instead of \ref{mike-lemma} 
(1). Thus $[G,N] \subset N^2$. 

Now let $(1+4A,4x)$ and $(1+4B,4y)$ be elements of $N$. Then \\
$[(1+4A,4x),(1+4B,4y)] \equiv (1+16(AB-BA), 16(Ay-Bx)) \hfill (\mbox{mod } 32),$ \\
thus $[N,N] \subset \{(1+16C,16z)(\mbox{mod } 2^{k+1}) \}$. The proof that any element of the form $(1+16C,16z)$ is generated by the 
squares of elements of the form $(1+4D,4w)$ is analogous to the proof of part 1 but is using \ref{mike-lemma} (3). Thus 
$[N,N] \subset N^4$.
\end{enumerate}
\end{proof}

As in \cite{MS}, let 
$$\Phi(N)=N^p[G,N]=[G,N]= \left\{  
\begin{array}{lll}
\{ (1+p^2A,p^2x) \ (\mbox{mod }p^{k+1})\} & \mbox{ if } & p \mbox{ is odd}, \\
\{ (1+8A,8x) \ (\mbox{mod }p^{k+1})\} & \mbox{ if } & p=2. \\
\end{array} \right.$$
For $i\ge 1$, define recursively a sequence of subgroups of $G$ by 
$$\begin{array}{rcl}
\Phi_G^1(N) & = & N \\
\Phi_G^{i+1}(N) & = & \Phi_G(\Phi_G^i(N)) = \left\{ 
\begin{array}{lll}
\{ (1+p^iA,p^ix) \ (\mbox{mod }p^{k+1})\} & \mbox{ if } & p \mbox{ is odd}, \\
\{ (1+2^{i+1}A,2^{i+1}x) \ (\mbox{mod }2^{k+1})\} & \mbox{ if } & p=2. \\
\end{array} \right. 
\end{array}$$

The $p$-power map induces a map 
$$\Phi_G^i(N) \buildrel p \over \longrightarrow \Phi_G^{i+1}.$$
For $i\ge 0$, let 
$$\begin{array}{l}
G_i = G/ \Phi_G^{i+1}(N) \cong \Gamma_{n,i}^{\scriptHol}, \\
A_i = \Phi_G^i(N)/\Phi_G^{i+1}(N) \cong \Ker_{n,i}^{\scriptHol}
\end{array}$$
where $\Gamma_{n,i}^{\scriptHol}$ and $\Ker_{n,i}^{\scriptHol}$ are as in section \ref{congruence_subgroups}. 
For $1\le i\le k-1$ in case $p$ is odd, and for $0\le i\le k-2$ in case $p=2$, $A_i$ is a vector space over $F_p$ of dimension $n^2+n$, 
thus $N$ is uniform, or almost uniformly embedded in $G$. Also, each $A_i$ is central in $G_i$, and we have successive central 
extensions 
$$1 \longrightarrow A_i \longrightarrow G_i \longrightarrow G_{i-1} \longrightarrow 1 $$
for $1\le i\le k-1$ in case $p$ is odd, and for $0\le i\le k-2$ in case $p=2$, 
i.e. this construction brings us to the uniform tower of groups $G_i\equiv \Gamma_{n,i}^{\scriptHol}$ studied in section 
\ref{congruence_subgroups}. Furthermore, since $\{ \Gamma_{n,k}^{\scriptHol}\}$ is an {\it infinite} uniform tower, $N$ is 
$\Omega$-extendable in $G$ for any $k$, thus all the results of section 3 of \cite{MS} apply to our 
tower. Note that the formula \ref{bracket-formula}, and thus our formulas for the Lie algebra structure constants and formulas for 
the Bockstein homomorphisms are valid for any $p$, so theorem \ref{congruence-ring-Bock} holds for any prime number $p$. 

\chapter{Wreath product}

\thispagestyle{kisa}

\label{wreath-chapter}

\section{Definition of wreath product}

\begin{definition}
The wreath product of a permutation group $P$ on $q$ letters (i.e. a subgroup of $S_q$) and a group $G$ is the 
semi-direct product 
$$P \wr G = P \ltimes G^q$$
where $P$ acts on the $G^q$ by permuting the factors. More precisely, elements of $ P\ltimes G^q$ are 
$(q+1)$-tuples $(\sigma, x_1, \ldots , x_q)$ where $\sigma \in P$ and $f_i \in G$, and the product is given by 
$$(\sigma, x_1, \ldots , x_q)(\tau, y_1, \ldots , y_q) = (\sigma\tau, x_{\tau^{-1}(1)}y_1, \ldots , x_{\tau^{-1}(q)}y_q).$$
\end{definition}

\begin{remark} There are different conventions for the multiplication of permutations, and hence different definitions of the 
multiplication in the wreath product. We use the convention $(\sigma \tau ) (i) = \sigma (\tau (i))$. If the reader prefers the 
multiplication of permutations to be given by $(\sigma \tau ) (i) = \tau (\sigma (i))$, then they should replace the permutations 
on the right hand sides of all our definitions by their inverses. 
\end{remark}

Define $P \buildrel \phi \over \hookrightarrow \Aut(G^q)$ by 
$$\phi(\sigma)(x_1, \ldots , x_q) = (x_{\sigma(1)}, \ldots , x_{\sigma(q)}).$$
Then the wreath product $P \wr G$ is the pullback $\xi$ in the following diagram:

$$\begin{CD}
\xi         @>>>  P            \\
@VVV              @VV{\phi}V   \\
\Hol(G^q)   @>>>  \Aut(G^q)
\end{CD}  
$$

\noindent i.e. we have a map of split extensions

$$\begin{CD}
1  @>>>  G^q  @>>>  P \wr G    @>>>  P          @>>>  1  \\
@.       @|         @VVV             @VVV             @. \\
1  @>>>  G^q  @>>>  \Hol(G^q)  @>>>  \Aut(G^q)  @>>>  1
\end{CD}
$$

We know from example \ref{product_example} that $(\Aut(G))^q$ and $(\Hol(G))^q$ can be regarded as subgroups 
of $\Aut(G^q)$ and $\Hol(G^q)$ respectively. We will now show that there are intermediate subgroups of 
$\Aut(G^q)$ and $\Hol(G^q)$ which are isomorphic to $S_q \wr \Aut(G)$ and $S_q \wr \Hol(G)$.

An element of $S_q \wr \Aut(G)$ is an $(n+1)$-tuple $(\sigma, f_1, \ldots , f_q)$ where
$\sigma\in S_q$ and $f_i\in \Aut(G)$, and an element of $S_q \wr \Hol(G)$ has the form
$(\sigma, (f_1,x_1), \ldots , (f_q,x_q))$ where $\sigma\in S_q$, $f_i\in \Aut(G)$, and $x_i\in G$.

Define $S_q \wr \Aut(G) \buildrel i \over \rightarrow \Aut(G^q)$ by
$$(i(\sigma, f_1, \ldots , f_q)) (x_1, \ldots x_q) = (f_{\sigma(1)}(x_{\sigma(1)}), \ldots , f_{\sigma(q)}(x_{\sigma(q)})).$$
  
Define now $S_q \wr \Hol(G) \buildrel j \over \rightarrow \Hol(G^q)$ by
$$j(\sigma, (f_1,x_1), \ldots, (f_q,x_q)) = (i(\sigma, f_1, \ldots , f_q), (x_1, \ldots, x_q)).$$

Here is the verification that $i$ and $j$ are group homomorphisms: 

\bigskip

\noindent $\begin{array}{l}
(i(\sigma, f_1,\ldots,f_q)) ((i(\tau,g_1, \ldots,g_q))(x_1,\ldots,x_q)) = \\
= (i(\sigma, f_1,\ldots,f_q)) ((g_{\tau(1)}(x_{\tau(1)}),\ldots,g_{\tau(q)}(x_{\tau(q)}))) = \\
= (f_{\sigma(1)}(g_{\sigma(\tau(1))}(x_{\sigma(\tau(1))})), \ldots, f_{\sigma(q)}(g_{\sigma(\tau(q))}(x_{\sigma(\tau(q))}))),
\end{array}$ \ 

\bigskip

\noindent $\begin{array}{l}
(i((\sigma, f_1,\ldots,f_q)(\tau, g_1,\ldots,g_q)))(x_1,\ldots,x_q) = \\
= (i(\sigma\tau, f_{\tau^{-1}(1)}g_1, \ldots, f_{\tau^{-1}(q)}g_q)) (x_1, \ldots, x_q) = \\
= (f_{(\sigma\tau)\tau^{-1}(1)}g_{(\sigma\tau)(1)}(x_{(\sigma\tau)(1)}), \ldots, 
f_{(\sigma\tau)\tau^{-1}(q)}g_{(\sigma\tau)(q)}(x_{(\sigma\tau)(q)})).
\end{array}$  \ 

\bigskip

\noindent $\begin{array}{l}
(j(\sigma, (f_1,x_1), \ldots, (f_q,x_q))) (j(\tau, (g_1,y_1), \ldots, (g_q,y_q))) = \\
= (i(\sigma, f_1,\ldots,f_q), (x_1,\ldots,x_q)) (i(\tau, g_1,\ldots,g_q), (y_1,\ldots,y_q)) = \\
= ((i(\sigma, f_1,\ldots,f_q))(i(\tau, g_1,\ldots,g_q)), ((i(\tau, g_1,\ldots,g_q))^{-1}(x_1,\ldots,x_q))(y_1,\ldots,y_q)) = \\
= (i((\sigma, f_1,\ldots,f_q)(\tau, g_1,\ldots,g_q)), ((i((\tau, g_1,\ldots,g_q)^{-1}))(x_1,\ldots,x_q))(y_1,\ldots,y_q)) = \\
= (i(\sigma\tau, f_{\tau^{-1}(1)}g_1, \ldots, f_{\tau^{-1}(q)}g_q), 
((i(\tau^{-1}, g_{\tau(1)}^{-1}, \ldots, g_{\tau(q)}^{-1}))(x_1,\ldots,x_q))(y_1,\ldots,y_q)) = \\
= (i(\sigma\tau, f_{\tau^{-1}(1)}g_1, \ldots, f_{\tau^{-1}(q)}g_q), 
(g_1^{-1}(x_{\tau^{-1}(1)}), \ldots, g_q^{-1}(x_{\tau^{-1}(q)}))(y_1, \ldots, y_q)) = \\
= (i(\sigma\tau, f_{\tau^{-1}(1)}g_1, \ldots, f_{\tau^{-1}(q)}g_q), 
(g_1^{-1}(x_{\tau^{-1}(1)})y_1, \ldots, g_q^{-1}(x_{\tau^{-1}(q)})y_q)), 
\end{array}$  \ 

\bigskip

\noindent $\begin{array}{l}
j((\sigma, (f_1,x_1), \ldots, (f_q,x_q)) (\tau, (g_1,y_1), \ldots, (g_q,y_q))) = \\
= j(\sigma\tau, (f_{\tau^{-1}(1)},x_{\tau^{-1}(1)})(g_1,y_1), \ldots, (f_{\tau^{-1}(q)},x_{\tau^{-1}(q)})(g_q,y_q)) = \\
= j(\sigma\tau, (f_{\tau^{-1}(1)}g_1, g_1^{-1}(x_{\tau^{-1}(1)})y_1), \ldots, (f_{\tau^{-1}(q)}g_q, g_q^{-1}(x_{\tau^{-1}(q)})y_q)) = \\
= (i(\sigma\tau, f_{\tau^{-1}(1)}g_1,\ldots,f_{\tau^{-1}(q)}g_q), (g_1^{-1}(x_{\tau^{-1}(1)})y_1,\ldots,g_q^{-1}(x_{\tau^{-1}(q)})y_q)).
\end{array}$ \

\bigskip

Therefore we have a map of split extensions
$$\begin{CD}
1  @>>>  G^q  @>>>  S_q \wr \Hol(G)  @>>>  S_q \wr \Aut(G)  @>>>  1   \\ 
@.       @|         @VVjV                  @VViV                  @.  \\
1  @>>>  G^q  @>>>  \Hol(G^q)        @>>>  \Aut(G^q)        @>>>  1
\end{CD}
$$

Also, the following diagrams commute:
$$\begin{CD}
S_q \wr S_n \wr \Aut(G)  @>>>  S_{qn} \wr \Aut(G)  \\
@VVV                           @VVV                \\
S_q \wr \Aut(G^n)        @>>>  \Aut(G^{qn})
\end{CD}$$

and

$$\begin{CD}
S_q \wr S_n \wr \Hol(G)  @>>>  S_{qn} \wr \Hol(G)  \\
@VVV                           @VVV                \\
S_q \wr \Hol(G^n)        @>>>  \Hol(G^{qn})
\end{CD}$$

\section{Cohomology of $\displaystyle \Hol \left ( \moplus_n \Z_p \right)$ at noncharacteristic primes}

Wreath products come up in the calculation of the homology of classical groups: $G(k) =$ $GL(k,\F_p)$, $O(k,\F_p)$, $Sp(k,\F_p)$, or 
$U(k,\F_{p^2}$) with coefficients in $\F_l$ where $(l,p)=1$ (see \cite{FP}). Namely, the inclusion 
$$S_q \wr G(c) \hookrightarrow G(qc)$$
(this is a special case of the inclusion $S_q \wr \Aut(G) \rightarrow \Aut(G^q)$ above) 
induces an epimorphism in mod $l$ homology, where $c=2$ if $l=2$ and $c$ is the smallest index for which $|G(c)|$ is divisible by $l$ 
if $l$ is odd. Using Quillen's lemma, this is improved as follows: the composite 
$$G(c)^q \hookrightarrow S_q \wr G(c) \hookrightarrow G(qc)$$ induces an epimorphism in mod $l$ homology. 

Since $$\left| \Hol \left( \moplus_k \Z_p \right) \right| = |\moplus_k \Z_p| \cdot |GL(k,\Z_p)| = p^k \cdot |GL(k,\Z_p)|,$$
the results of \cite[III, \S 5]{FP} apply to holomorphs. However, we are not getting any new information this way because 
$$H^*\left(\Hol \left( \moplus_k \Z_{p^r} \right) ; \F_l \right) \cong H^* \left( GL(n, \Z_{p^r}) ; \F_l \right) $$ 
for $(l,p)=1$.

\section{Permutative categories}

Recall (see e.g. \cite{FP}) that a permutative category $(\C, \square, 0, c)$ is a small category $\C$ with a bifunctor 
$\square: \C \times \C \rightarrow \C,$ an object $0,$ and a natural isomorphism $c:\square \rightarrow \square \tau$ (where 
$\tau: \C \times \C \rightarrow \C \times \C$ is a transformation), i.e. a function 
$c: \Ob \C \times \Ob \C \rightarrow \Mor \C$ such that for any $\alpha, \beta \in \Ob\C$, 
$\mbox{source}(c(\alpha,\beta)) = \alpha \square \beta$ and $\mbox{target}(c(\alpha,\beta)) = \beta \square \alpha$, and such 
that for any $A,B \in \Hom (\alpha, \beta)$, the diagram 
\begin{equation}\label{permutative-def-diagram}
\begin{CD}
\alpha \square \beta  @>{A \square B}>>  \gamma \square \delta \\ 
@VVc(\alpha,\beta)V                      @VVc(\gamma,\delta)V  \\
\beta \square \alpha  @>{B \square A}>>  \delta \square \gamma 
\end{CD}\end{equation} 
commutes, which satisfy the following properties: 
\begin{enumerate}
\item $\square$ is associative,
\item $0$ is a left and right unit for $\square$,
\item for all $\alpha \in \Ob \C$, $c(\alpha,0)=1_{\alpha}=c(0,\alpha)$,
\item for all $\alpha, \beta \in \Ob \C,$ $c(\beta,\alpha)=(c(\alpha,\beta))^{-1}$,
\item for all $\alpha, \beta, \gamma \in \Ob \C,$ 
$(c(\gamma,\alpha)\square 1_{\beta}) \circ c(\alpha\square\beta,\gamma) = 1_{\alpha} \square c(\beta,\gamma).$ 
\end{enumerate}

A permutative functor between permutative categories is a functor that commutes with $\square$ and $c$, and sends $0$ to $0$. 

As shown in \cite{FP}, each family of classical groups with direct sum homomorphism $G(m) \times G(n) \rightarrow G(m+n)$ forms a 
permutative category. Below we give some other examples of permutative categories. 
Let $R$ be a commutative ring with a unit. Recall the inclusion 
$\Hol \left( \moplus_m R \right) \times \Hol \left( \moplus_n R \right) \hookrightarrow \Hol \left( \moplus_{m+n} R \right)$ 
from example \ref{product_example}. The family $\left\{ \Hol\left(\moplus_n R\right), n \ge 0 \right\}$ 
forms a permutative category (denoted by ${\cal H}ol(R)$) as follows. 
Objects are nonnegative integers, and morphisms are given 
by $$\Hom(m,n)=\left\{ \begin{array}{ll} 
0,                 & \mbox{ if } m\ne n \\
\Hol\left(\moplus_n R\right), & \mbox{ if } m = n.   
\end{array} \right. $$ 
Consider $\Hol\left(\moplus_n R\right)$ as a subgroup of $GL(n+1,R)$ as described in lemma \ref{hol-matrix-lemma}. 
Define $\square$ on objects by 
$$m \square n = m+n$$
and on morphisms by 
$$\left[ \begin{array}{c|c} 
1 & x^T   \\    \hline    0 & (M^{-1})^T
\end{array} \right] 
\square 
\left[ \begin{array}{c|c}
1 & y^T   \\   \hline   0 & (N^{-1})^T
\end{array} \right] = 
\left[  \begin{array}{c|c|c} 
1 & x^T & y^T   \\   \hline   0 & (M^{-1})^T & 0   \\   \hline   0 & 0 & (N^{-1})^T 
\end{array} \right] $$
and define the function $c$ by 
$$c(m,n)=\left[  \begin{array}{c|c|c} 
1 & 0   & 0    \\   \hline   0 & 0   & I_n   \\    \hline   0 & I_m & 0 
\end{array} \right]$$
It is easy to check that these functions $\square$ and $c$ satisfy the required properties.

The surjective and splitting homomorphisms between $\Hol\left( \moplus_n R \right)$ and $GL(n,R)$ in the split extension for the 
holomorph form permutative functors between the permutative categories ${\cal H}ol(R)$ and ${\cal GL}(R)$. However, the inclusions 
$\Hol\left( \moplus_n R \right) \hookrightarrow GL(n+1,R)$ do not give a permutative functor. 

More generally, the inclusions 
$$\Aut(G^m) \times \Aut(G^n) \rightarrow \Aut(G^{m+n}),$$ 
$$\Hol(G^m) \times \Hol(G^n) \rightarrow \Hol(G^{m+n}),$$ 
and
$$\Aut\left(\free_m G\right) \times \Aut\left(\free_n G\right) \rightarrow \Aut\left(\free_{m+n}G\right)$$ 
give rise to permutative categories. 

\begin{theorem} \label{permutative-theorem} 
For any group $G$, the families 
$$\begin{array}{ll}
 (i) & \left\{ \Aut(G^n), n \ge 0 \right\}, \\ 
(ii) & \left\{ \Hol(G^n), n \ge 0 \right\}, \mbox{ and} \\ 
(iii) & \displaystyle \left\{ \Aut\left(\free_n G\right), n \ge 0 \right\}
\end{array}$$
form permutative categories, denoted by ${\cal A}ut(G^n)$ and ${\cal H}ol(G^n)$, with objects nonnegative integers, morphisms 
given by 
$$ \Hom(m,n)=\left\{ \begin{array}{ll}
0,                 & \mbox{ if } m\ne n \\
\begin{array}{@{}ll}
 (i)  & \Aut(G^n) \\
(ii)  & \Hol(G^n) \\
(iii) & \Aut\left(\free_n G\right)
\end{array},
& \mbox{ if } m = n,
\end{array} \right. $$
bifunctor $\square$ defined on objects by 
$$m \square n = m+n$$
and on morphisms by
$$\begin{array}{ll}
 (i)  & f_m \square g_n = f_m \times g_n, \\  
(ii)  & (f_m, (x_1, \ldots, x_m)) \square (g_n, (y_1, \ldots, y_n)) = (f_m \times g_n, (x_1,\ldots,x_m,y_1,\ldots,y_n)), \\
(iii) & f_m \square g_n = f_m \times g_n
\end{array}$$ 
where in (i) and (ii) $f_m\in\Aut(G^m),$ $g_n\in\Aut(G^n),$ and $f_m \times g_n\in\Aut(G^{m+n})$ is the product of $f_m$ and $g_n$, 
in (iii) $f_m\in\Aut\left(\free_m G\right),$ $g_n\in\Aut\left(\free_n G\right),$ and $f_m \times g_n\in\Aut\left(\free_{m+n} G \right)$ 
is the product of $f_m$ and $g_n$, and 
$$\begin{array}{ll}
 (i)  & c(m,n) = s_{m,n}, \\ 
(ii)  & c(m,n) = (s_{m,n}, 1_{G^{m+n}}), \\ 
(iii) & c(m,n)(x_i) = \left\{ 
\begin{array}{lll}
x_{i+n} & \mbox{if} & 1\le i\le m, \\
x_{i-m} & \mbox{if} & m+1\le i\le n, 
\end{array}
\right. 
\end{array}$$ 
where $s_{m,n}\in\Aut(G^{m+n})$ is defined by  
$$s_{m,n}(x_1,\ldots,x_m,x_{m+1},\ldots,x_{m+n}) = (x_{m+1},\ldots,x_{m+n},x_1,\ldots,x_m),$$
and $1_{G^{m+n}}$ is the identity element in $G^{m+n}$.    
There are permutative functors between ${\cal A}ut(G^n)$ and ${\cal H}ol(G^n)$ given by the identity maps on objects 
and the natural maps between $\Aut(G^m)$ and $\Hol(G^m)$ on morphisms.
\end{theorem} 

The proof of this theorem is a straightforward verification given next.

\begin{proof}
Commutativity of \ref{permutative-def-diagram}: 
\begin{enumerate}
\item[(i)] since 
$(x_1,\ldots,x_{k})=(x_1,1,\ldots,1)(1,x_2,1,\ldots,1) \ldots (1,\ldots,1,x_k)$, it suffices to check 
commutativity of the diagram on elements of the form \\ $(1,\ldots,1,x,1,\ldots,1)$. Let $x$ be in the 
$i$-th place. If $1\le i\le m$, we have 
$$\begin{CD}
(1,\ldots,1,x,1,\ldots,1;1,\ldots,1) @>{f\times g}>> (x_1,\ldots,x_m;1,\ldots,1) \\
@VVc(m,n)V                                           @VVc(m,n)V                  \\          
(1,\ldots,1;1,\ldots,1,x,1,\ldots,1) @>{g\times f}>> (1,\ldots,1;x_1\ldots,x_m). 
\end{CD}$$
The case $m+1\le i\le m+n$ is similar.
\item[(ii)] $c(m,n)((f, (x_1, \ldots, x_m))\square(g,(y_1, \ldots, y_n))) =$ \\ 
$=(s_{m,n},1_{G^{m+n}})(f\times g, (x_1, \ldots, x_m, y_1, \ldots, y_n)) =$ \\
$=(s_{m,n}\circ (f\times g), (x_1, \ldots, x_m, y_1, \ldots, y_n))$. \\

$(g, (y_1, \ldots, y_n))\square((f, (x_1, \ldots, x_m)))c(m,n) = $ \\
$=(g\times f, (y_1, \ldots, y_n, x_1, \ldots, x_m))(s_{m,n},1_{G^{m+n}}) = $ \\
$=((g\times f)\circ s_{m,n}, (s_{m,n})^{-1}(y_1, \ldots, y_n, x_1, \ldots, x_m)1_{G^{m+n}}) = $ \\
$=((g\times f)\circ s_{m,n}, (s_{n,m})(y_1, \ldots, y_n, x_1, \ldots, x_m)) = $ \\
$=((g\times f)\circ s_{m,n}, (x_1, \ldots, x_m, y_1, \ldots, y_n)) = $ \\
$=(s_{m,n}\circ (f\times g), (x_1, \ldots, x_m, y_1, \ldots, y_n))$ by part (i). 

\item[(iii)] As in part (i), we will check the commutativity of the diagram on the generators. If $1\le i\le m$, we have 
$$\begin{CD}
x_i         @>{f\times g}>> f(x_i)     \\
@VVc(m,n)V                  @VVc(m,n)V  \\
x_{i+n}     @>{g\times f}>> f(x_{i+n}) 
\end{CD}$$
The case $m+1\le i\le m+n$ is similar. 
\end{enumerate}
Axioms: 
\begin{enumerate}
\item $\square$ is associative because multiplication in $\Aut(G^n)$ is. 
\begin{enumerate}
\item[(i)] $(f_m \times g_n)\times h_p = f_m\times (g_n \times h_p),$
\item[(ii)] $(f_m \times g_n)\times h_p, (x_1, \ldots, x_m,y_1, \ldots, y_n, z_1, \ldots, z_p))=$ \\
$ = (f_m\times (g_n \times h_p), (x_1, \ldots, x_m,y_1, \ldots, y_n, z_1, \ldots, z_p)),$
\item[(iii)] same as (i).
\end{enumerate}

\item 0 is a unit in $\Z$.

\item 
\begin{enumerate}
\item[(i)] $c(m,0)=s_{m,0}=1_{\scriptAut(G^m)}=s_{0,m}=c(0,m)$,
\item[(ii)] $c(m,0)=(s_{m,0}, 1_{G^m}) = (s_{0,m}, 1_{G^m}) = c(0,m),$
\item[(iii)] $\displaystyle c(m,0)=1_{\scriptAut\left( \displaystyle \free_m G\right)}=c(0,m).$
\end{enumerate}

\item
\begin{enumerate}
\item[(i)] $c(m,n)c(n,m)(x_1,\ldots, x_{m+n}) = s_{m,n}s_{n,m}(x_1,...,x_{m+n})= $ \\  
$= s_{m,n}(x_{n+1},\ldots,x_{m+n},x_1,\ldots,x_n) = (x_1,\ldots,x_{m+n}) = 1_{\scriptAut(G^{m+n})},$
\item[(ii)] $c(m,n)c(n,m) = (s_{m,n},1_{G^{m+n}})(s_{n,m},1_{G^{m+n}}) = (s_{m,n}s_{n,m}, 1_{G^{m+n}}) = $\\ 
$= (1_{\scriptAut(G^{m+n})}, 1_{G^{m+n}}) = 1_{\scriptHol(G^{m+n})},$
\item[(iii)] if $1\le i\le m$, then $c(m,n)c(n,m)(x_i) = c(m,n)(x_{i+m}) = x_{i+m-m} = x_i;$
the case $m+1\le i\le m+n$ is similar. 
\end{enumerate}

\item
\begin{enumerate}
\item[(i)] $((c(p,m)\square 1_n)\circ c(m+n,p))(x_1,\ldots,x_{m+n+p}) = $\\
$= (s_{p,m} \times 1_{\scriptAut(G^n)}) s_{m+n,p}(x_1,\ldots,x_{m+n+p}) = $\\
$= (s_{p,m} \times 1_{\scriptAut(G^n)}) (x_{m+n+1}, \ldots, x_{m+n+p}, x_1, \ldots, x_{m+n}) = $ \\
$= (x_1, \ldots, x_m, x_{m+n+1}, \ldots, x_{m+n+p}, x_{m+1}, \ldots, x_{m+n}) = $ \\
$= (1_{\scriptAut(G^m)} \times s_{n,p})(x_1,\ldots,x_{m+n+p}) = 1_m \square c(n,p)(x_1,\ldots,x_{m+n+p}),$ 
\item[(ii)] $(c(p,m)\square 1_n)\circ c(m+n,p) = $ \\ 
$= ((s_{p,m}, 1_{G^{m+p}}) \square (1_{\scriptAut(G^n)},1_{G^n})) \circ (s_{m+n,p},1_{G^{m+n+p}}) = $\\ 
$= (s_{p,m} \times 1_{\scriptAut(G^n)}, 1_{G^{m+n+p}}) \circ (s_{m+n,p}, 1_{G^{m+n+p}}) = $\\
$= ((s_{p,m} \times 1_{\scriptAut(G^n)})s_{m+n,p}, 1_{G^{m+n+p}}) = 
(1_{\scriptAut(G^m)} \times s_{n,p}, 1_{G^{m+n+p}}) = $\\ 
$= (1_{\scriptAut(G^m)}, 1_{G^m}) \square (s_{n,p},1_{G^{n+p}}) = 
(1_m \square c(n,p)),$
\item[(iii)] if $1\le i \le m$, then $((c(p,m)\square 1_n)\circ c(m+n,p)) (x_i) = $\\
$= (c(p,m)\square 1_n)(x_{i+p})= x_i = (1_m \square c(n,p)) (x_i)$; \\ 
if $m+1\le i\le m+n$, then $((c(p,m)\square 1_n)\circ c(m+n,p)) (x_i) = $\\
$= (c(p,m)\square 1_n)(x_{i+p}) = x_{i+p} = (1_m \square c(n,p)) (x_i)$; \\
if $m+n+1\le i\le m+n+p$, then $((c(p,m)\square 1_n)\circ c(m+n,p)) (x_i) = $\\
$= (c(p,m)\square 1_n)(x_{i-m-n}) = x_{i-m-n+m} = x_{i-n} = (1_m \square c(n,p)) (x_i)$.
\end{enumerate}

\end{enumerate}
\end{proof}

These permutative categories give rise to infinite loop spaces which will be studied later.

\begin{remark}
The family $\left\{ \Hol\left(\free_n G\right), n\ge 0\right\}$ is not a permutative category because there is no ``natural'' map 
$\Hol\left(\free_m G\right) \times \Hol\left(\free_n G\right) \rightarrow \Hol\left(\free_{m+n} G\right)$.
\end{remark}

\chapter{Injectives in the category of groups}

\thispagestyle{kisa}

\label{injectives-chapter}

\section{Introduction}

This chapter is independent from the previous ones. In the next section we give a short proof of a theorem 
(\ref{injectives-theorem}) proved by S. Eilenberg and J. C. Moore in \cite{EM}. 

\begin{definition}
An object $I$ in a category $\C$ is called injective if for any monomorphism $K \rightarrow L$, and any map   
$K \rightarrow I$, there is a map $L \rightarrow I$ such that the diagram 

\centerline{\setlength{\unitlength}{2900sp}%
\begingroup\makeatletter\ifx\SetFigFont\undefined%
\gdef\SetFigFont#1#2#3#4#5{%
  \reset@font\fontsize{#1}{#2pt}%
  \fontfamily{#3}\fontseries{#4}\fontshape{#5}%
  \selectfont}%
\fi\endgroup%
\begin{picture}(3600,2055)(3601,-3161)
\thinlines
\special{ps: gsave 0 0 1 setrgbcolor}\put(4420,-1470){\line( 1, 0){800}}
\special{ps: grestore}\special{ps: gsave 0 0 1 setrgbcolor}\put(5755,-1470){\line( 1, 0){800}}
\special{ps: grestore}\special{ps: gsave 0 0 1 setrgbcolor}\put(5500,-1711){\line( 0,-1){800}}
\special{ps: grestore}\special{ps: gsave 0 0 1 setrgbcolor}\put(6635,-1710){\line(-1,-1){230}}
\special{ps: grestore}\special{ps: gsave 0 0 1 setrgbcolor}\put(6300,-2045){\line(-1,-1){230}}
\special{ps: grestore}
\put(5270,-1537){\makebox(0,0)[rm]{\smash{\SetFigFont{12}{15.0}{\rmdefault}{\mddefault}{\updefault}\special{ps:
gsave 0 0 1 setrgbcolor}{$\rightarrow$}\special{ps: grestore}}}}
\put(6605,-1537){\makebox(0,0)[rm]{\smash{\SetFigFont{12}{15.0}{\rmdefault}{\mddefault}{\updefault}\special{ps:
gsave 0 0 1 setrgbcolor}{$\rightarrow$}\special{ps: grestore}}}}
\put(5500,-2561){\makebox(0,0)[cb]{\smash{\SetFigFont{12}{15.0}{\rmdefault}{\mddefault}{\updefault}\special{ps:
gsave 0 0 1 setrgbcolor}{$\downarrow$}\special{ps: grestore}}}}
\put(5715,-2559){\makebox(0,0)[lb]{\smash{\SetFigFont{12}{15.0}{\rmdefault}{\mddefault}{\updefault}\special{ps:
gsave 0 0 1 setrgbcolor}{$\swarrow$}\special{ps: grestore}}}}
\put(4160,-1560){\makebox(0,0)[lb]{\smash{\SetFigFont{12}{15.0}{\rmdefault}{\mddefault}{\updefault}\special{ps:
gsave 0 0 1 setrgbcolor}{$1$}\special{ps: grestore}}}}
\put(5380,-1560){\makebox(0,0)[lb]{\smash{\SetFigFont{12}{16.0}{\rmdefault}{\mddefault}{\updefault}\special{ps:
gsave 0 0 1 setrgbcolor}{$K$}\special{ps: grestore}}}} 
\put(6700,-1560){\makebox(0,0)[lb]{\smash{\SetFigFont{12}{16.0}{\rmdefault}{\mddefault}{\updefault}\special{ps:
gsave 0 0 1 setrgbcolor}{$L$}\special{ps: grestore}}}}
\put(5450,-2960){\makebox(0,0)[lb]{\smash{\SetFigFont{12}{16.0}{\rmdefault}{\mddefault}{\updefault}\special{ps:
gsave 0 0 1 setrgbcolor}{$I$}\special{ps: grestore}}}}
\end{picture}}

\noindent commutes. 
\end{definition}

\begin{theorem}
\label{injectives-theorem} (S. Eilenberg, J. C. Moore)
The only injective object in the category of groups is the trivial group.
\end{theorem}

We will need the following lemma. It follows from a classical proof of the fact that the free group on two letters contains 
the free group on countably many letters. We recall a proof of the lemma because we do not 
know reference for the lemma in this exact wording.

\begin{lemma}
\label{free-groups-lemma} (Classical) 
Let $F[a,b]$ denote the free group on letters $a$ and $b$. Then the group homomorphism 
$$F[a,b] \rightarrow F[c,d]$$ 
given by 
$$\begin{array}{rcl} 
a & \mapsto & c \\ 
b & \mapsto & dcd^{-1}
\end{array}$$
is an injection.
\end{lemma}

\begin{figure}
\centerline{\epsfbox{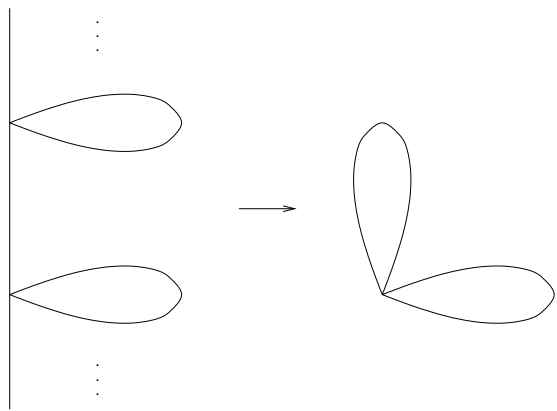}}
\vspace{-3mm}
\caption{Covering space of a bouquet of two circles} 
\label{covering}
\end{figure}

\begin{figure}
\centerline{\epsfbox{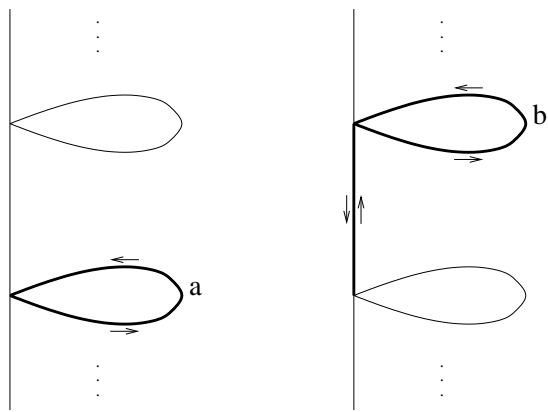}}
\vspace{-3mm}
\caption{Two loops whose classes generate subgroup $F[a,b]$ of the fundamental group of the covering space}
\label{loops}
\end{figure}

\begin{proof}
Consider the covering space of a bouquet of two circles shown in figure \ref{covering} on the next page. 
The covering space is 
homotopy equivalent to a bouquet of countably many circles, and hence its fundamental group is isomorphic to the 
free group on countably many generators. Classes of loops {\bf a} and {\bf b} shown in figure \ref{loops} generate a free subgroup 
$F[a,b]$. The fundamental group of the covering space injects into the fundamental group of the bouquet which is 
isomorphic to the free group on two generators, say, $F[c,d]$ (see figure 6.3). Let $i$ denote the 
restriction of this injection on $F[a,b]$, then we have 
$$\begin{array}{rcl}
F[a,b] & \buildrel i \over \hookrightarrow & F[c,d], \\
i(a)   & =                                 & c,      \\
i(b)   & =                                 & dcd^{-1}.
\end{array}$$
\end{proof}

\newpage
 
\vfill

\smallskip

\vfill 

\vfill 

\centerline{\epsfbox{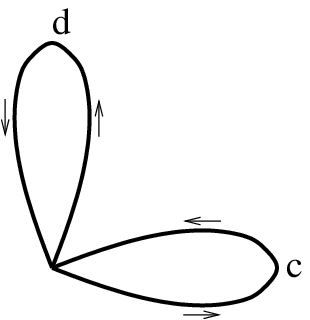}}

\noindent {Figure 6.3: \ The fundamental group of the bouquet of two circles is free on two generators}

\vfill

\newpage

\section{Proof of theorem \ref{injectives-theorem}}

\begin{proof}

Suppose $G$ is injective, and let $x\in G$ be any element. Let the homomorphism 
$$F[a,b] \buildrel f \over \rightarrow G$$ be given by  
$$\begin{array}{rcl}
f(a) & = & 1 \\
f(b) & = & x 
\end{array}$$
and let $F[a,b] \buildrel i \over \rightarrow F[c,d]$ be as in lemma \ref{free-groups-lemma}.
Then there exists a homomorphism 
$$F[c,d] \buildrel g \over \rightarrow G$$ 
such that the diagram 

\centerline{\setlength{\unitlength}{2900sp}%
\begingroup\makeatletter\ifx\SetFigFont\undefined%
\gdef\SetFigFont#1#2#3#4#5{%
  \reset@font\fontsize{#1}{#2pt}%
  \fontfamily{#3}\fontseries{#4}\fontshape{#5}%
  \selectfont}%
\fi\endgroup%
\begin{picture}(3600,2055)(3601,-3161)
\thinlines
\special{ps: gsave 0 0 1 setrgbcolor}\put(4400,-1470){\line( 1, 0){650}}
\special{ps: grestore}\special{ps: gsave 0 0 1 setrgbcolor}\put(5925,-1470){\line( 1, 0){650}}
\special{ps: grestore}\special{ps: gsave 0 0 1 setrgbcolor}\put(5500,-1711){\line( 0,-1){800}}
\special{ps: grestore}\special{ps: gsave 0 0 1 setrgbcolor}\put(6670,-1710){\line(-1,-1){230}}
\special{ps: grestore}\special{ps: gsave 0 0 1 setrgbcolor}\put(6335,-2045){\line(-1,-1){230}}
\special{ps: grestore}
\put(5090,-1537){\makebox(0,0)[rm]{\smash{\SetFigFont{12}{15.0}{\rmdefault}{\mddefault}{\updefault}\special{ps:
gsave 0 0 1 setrgbcolor}{$\rightarrow$}\special{ps: grestore}}}}
\put(6655,-1537){\makebox(0,0)[rm]{\smash{\SetFigFont{12}{15.0}{\rmdefault}{\mddefault}{\updefault}\special{ps:
gsave 0 0 1 setrgbcolor}{$\rightarrow$}\special{ps: grestore}}}}
\put(5500,-2561){\makebox(0,0)[cb]{\smash{\SetFigFont{12}{15.0}{\rmdefault}{\mddefault}{\updefault}\special{ps:
gsave 0 0 1 setrgbcolor}{$\downarrow$}\special{ps: grestore}}}}
\put(5750,-2559){\makebox(0,0)[lb]{\smash{\SetFigFont{12}{15.0}{\rmdefault}{\mddefault}{\updefault}\special{ps:
gsave 0 0 1 setrgbcolor}{$\swarrow$}\special{ps: grestore}}}}
\put(4160,-1560){\makebox(0,0)[lb]{\smash{\SetFigFont{12}{15.0}{\rmdefault}{\mddefault}{\updefault}\special{ps:
gsave 0 0 1 setrgbcolor}{$1$}\special{ps: grestore}}}}
\put(5500,-1560){\makebox(0,0)[cb]{\smash{\SetFigFont{12}{16.0}{\rmdefault}{\mddefault}{\updefault}\special{ps:
gsave 0 0 1 setrgbcolor}{$F[a,b]$}\special{ps: grestore}}}} 
\put(6710,-1560){\makebox(0,0)[lb]{\smash{\SetFigFont{12}{16.0}{\rmdefault}{\mddefault}{\updefault}\special{ps:
gsave 0 0 1 setrgbcolor}{$F[c,d]$}\special{ps: grestore}}}}
\put(5410,-2960){\makebox(0,0)[lb]{\smash{\SetFigFont{12}{16.0}{\rmdefault}{\mddefault}{\updefault}\special{ps:
gsave 0 0 1 setrgbcolor}{$G$}\special{ps: grestore}}}}
\put(6150,-1400){\makebox(0,0)[lb]{\smash{\SetFigFont{12}{16.0}{\rmdefault}{\mddefault}{\updefault}\special{ps:
gsave 0 0 1 setrgbcolor}{$i$}\special{ps: grestore}}}}
\put(5290,-2230){\makebox(0,0)[lb]{\smash{\SetFigFont{12}{16.0}{\rmdefault}{\mddefault}{\updefault}\special{ps:
gsave 0 0 1 setrgbcolor}{$f$}\special{ps: grestore}}}}
\put(6265,-2260){\makebox(0,0)[lb]{\smash{\SetFigFont{12}{16.0}{\rmdefault}{\mddefault}{\updefault}\special{ps:
gsave 0 0 1 setrgbcolor}{$g$}\special{ps: grestore}}}}
\end{picture}}
\noindent commutes. Then we have 
$$g(c)=g(i(a))=f(a)=1,$$
and 
$$x=f(b)=g(i(b))=g(dcd^{-1})=g(d)g(c)g(d^{-1})=g(d)1(g(d))^{-1}=1,$$
i.e. any element of $G$ is the identity. 
\end{proof}



\begin{thebibliography}{99999999}
\addcontentsline{toc}{chapter}{\bibname}
\thispagestyle{kisa}
\bibitem[AM]{AM}A. Adem and R. J. Milgram, Cohomology of finite groups, Springer-Verlag, 1994.
\bibitem[A]{A} G. Andrews, Number theory, Dover Publications, New York, 1994. 
\bibitem[BB]{BB}J. A. Beachy and W. D. Blair, Abstract algebra: Supplementary lecture notes, 2nd ed., Waveland Press, Illinois, 1996. 
\bibitem[BP]{BP}W. Browder and J. Pakianathan, ``Cohomology of uniformly powerful $p$-groups'', Trans. Amer. Math. Soc., 352 (2000)
2659-2688.
\bibitem[B]{B}K. Brown, Cohomology of groups, Graduate Texts in Mathematics 87, Springer-Verlag, 1982. 
\bibitem[CE]{CE}H. Cartan and S. Eilenberg, Homological algebra, Princeton University Press, 1956.
\bibitem[CCNPW]{AFG} J. H. Conway, R. T. Curtis, S. P. Norton, R. A. Parker, R. A. Wilson, Atlas of finite groups, Clarendon 
Press, Oxford, 1989.
\bibitem[D]{D}L. E. Dickson, ``A fundamental system of invariants of the general modular linear group with a solution of the 
form problem. Trans. Amer. Math. Soc., 12 (1911), 75-98. 
\bibitem[EM]{EM}S. Eilenberg and J. C. Moore, Foundations of relative homological algebra, Mem. Amer. Math. Soc., 55 (1965).
\bibitem[FP]{FP}Z. Fiedorowicz and S. Priddy, Homology of classical groups over finite fields and their associated infinite loop 
spaces, Springer-Verlag, 1978.  
\bibitem[FTY]{FTY}M. Furusawa, M. Tezuka, and N. Yagita, ``On the cohomology of classifying spaces of torus bundles and 
automorphic forms, J. London Math. Soc. (2) 37 (1988), 520-534.  
\bibitem[H]{H}J. Huebschmann, ``The mod $p$ cohomology rings of metacyclic groups'', Journal of Pure and Applied Algebra, 60 (1989), 
53-105. 
\bibitem[J1]{J1}C. Jensen, ``Homology of holomorphs of free groups'', submitted. 
\bibitem[J2]{J2}C. Jensen, ``Cohomology of $\Aut(F_n)$ in the $p$-rank two case'', Journal of Pure and Applied Algebra, 158 (2001), 
41-81. 
\bibitem[K]{K}A. G. Kurosh, The theory of groups, volume 1, Chelsea publishing company, New York, 1955.
\bibitem[M]{M}S. Mac Lane, Homology, Springer-Verlag, 1995. 
\bibitem[MP]{MP}J. Martino and S. Priddy, ``Classification of BG for groups with dihedral or quaternion Sylow 
2-groups'', Journal of Pure and Applied Algebra, 73 (1991), 13-21.
\bibitem[MS]{MS} P. A. Minh and P. Symonds, ``The cohomology of pro-$p$ groups with a powerfully embedded subgroup'', preprint,  
available at {\tt http://www.math.uga.edu/archive/minh-symonds.html}  
\bibitem[P]{P}J. Pakianathan, On the cohomology of a family of $p$-groups associated to Lie algebras, Princeton University thesis, 1997.
\bibitem[R]{R}D. Redmond, Number theory, Marcel Dekker, 1996.
\bibitem[S]{S}J.-P. Serre, Cohomologie Galoisienne, Lecture Notes in Mathematics 5, Springer-Verlag, 1964. 
\bibitem[W]{W}C. T. C. Wall, ``Resolutions for extensions of groups'', Proc. Cambridge Philos. Soc. 57 1961, 251-255.
\bibitem[We]{We} T. Weigel, Exp and log functors for the categories of powerful $p$-central groups and Lie algebras, 
Habilitationsschrift, Freiburg, 1994.
\end{thebibliography}
\end{document}